\documentclass[11pt]{article}

\usepackage{defs/packages}
\usepackage{defs/macros}
\usepackage{defs/diagrammatics}
\usepackage{defs/custom-biblatex-style}

\newcommand{\satex}[2][]{\vcenter{\hbox{\includegraphics[#1,scale=0.9]{satex/#2}}}}
\newcommand{\figfoam}[2][]{\vcenter{\hbox{\includegraphics[#1,scale=0.9]{foam/#2}}}}

\usepackage{tocloft}
\usepackage[a4paper, total={6in, 9in}]{geometry}

\counterwithin{figure}{section}

\usepackage{footnote}
\makesavenoteenv{tabular}
\usepackage{authblk}

\title{Rewriting modulo in diagrammatic algebras\\and application to categorification}
\author{Léo Schelstraete \orcidlink{0000-0001-7167-3964}\footnote{Max Planck Institute for Mathematics (Bonn, Germany). Contact: \texttt{leoschelstraete@gmail.com}}}
\date{\vspace{-6ex}}

\begin{document}

\makeatletter
\patchcmd{\@maketitle}{\vskip 2em}{\vskip -6em}{}{}
\makeatother
\maketitle

\begin{abstract}
  We develop a rewriting theory modulo suitable for higher linear structures. In particular, the theory is suited for diagrammatic algebras as they appear in categorification, representation theory and quantum topology.
  As an application, we use the theory to prove the basis conjecture of a certain super-2-category related to odd Khovanov homology.
  
  Our approach combines linear rewriting, higher rewriting and rewriting modulo. For diagrammatic algebras, the modulo rules typically capture a categorical property, such as pivotality.
  In the process, we revisit the foundations of these theories, including the notion of confluence.
  Other important tools include termination rules that depend on contexts, rewriting modulo invertible scalars, and a method to classify branchings modulo.
  
  This article includes an introduction to rewriting theory for non-experts.
\end{abstract}

\tableofcontents

\section{Introduction}

Monoidal categories, and more generally 2-categories, are ubiquitous in mathematics.
Their study categorifies classical algebra, exhibiting novel structures and phenomena \cite{Selinger_SurveyGraphicalLanguages_2011}; we refer to it as \emph{2-dimensional algebra}, or more loosely speaking, \emph{higher algebra}.
While symbols are the language of classical algebra, \emph{string diagrams} (dual to pasting diagrams) are the language of higher algebra.
For that reason, linear (strict) 2-categories presented using string diagrams are often called \emph{diagrammatic algebras}.

\begin{wrapfigure}[9]{r}{.38\textwidth}
  \captionsetup{font=small}
  \vspace*{-20pt}
  \begin{gather*}
    \TLtikz{
      \draw[line width=.1pt,dashed] (-.5,5) to (6.5,5);
      \TLstrand{1}{4}\TLstrand{2}{4}\TLstrand{3}{4}\TLcap{4}{4}
      \TLstrand{1}{3}\TLcup{2}{3}\TLstrand{4}{3}\TLstrand{5}{3}
      \TLstrand{1}{2}\TLcap{1.5}{2}\TLstrand{4}{2}\TLstrand{5}{2}
      \TLstrand{1}{1}\TLcup{1.5}{1}\TLdiag{3}{1}\TLantidiag{5}{1}
      \TLstrand{1}{0}\TLdiag{2}{0}\TLcap{4}{0}\TLstrand{6}{0}
      \TLcup{1}{-1}\TLstrand{4}{-1}\TLstrand{5}{-1}\TLstrand{6}{-1}
      \draw[line width=.1pt,dashed] (-.5,-1) to (6.5,-1);
    }[scale=.6]
    =
    {(q+q^{-1})}
    \TLtikz{
      \draw[line width=.1pt,dashed] (-.5,3) to (2.5,3);
      \TLcup{1}{2}
      \draw[TL_black_strand] (2,0) to (0,3);
      \TLcap{0}{0}
      \draw[line width=.1pt,dashed] (-.5,0) to (2.5,0);
    }
  \end{gather*}
  \vspace*{-8pt}
  \caption{A relation in the Temperley--Lieb category. Thanks to the diagrammatics, it is easily understood as ``evaluate the closed loop to $q+q^{-1}$, and apply a planar isotopy''.}
  \label{fig:intro_temperley-lieb}
\end{wrapfigure}
%
Diagrammatic algebras are particularly prevalent in representation theory and low-dimensional topology.
A classical example is the Temperley--Lieb category: it describes (some of) the representation theory of the quantum group $U_q(\mathfrak{sl}_2)$. It is also at the heart of the definition of the Jones polynomial \cite{Jones_PolynomialInvariantKnots_1985,Kauffman_StateModelsJones_1987}, an invariant of knots whose discovery birthed the field of quantum topology.
The morphisms of the Temperley--Lieb category are pictured as $\bC(q)$\nbd-linear combinations of certain diagrams (see \cref{fig:intro_temperley-lieb}).
This allows some of its defining relations to be interpreted as certain isotopies of diagrams---a typical feature of diagrammatic algebras.

Other examples of diagrammatic algebras in representation theory include:
\begin{itemize}
  \item diagrammatic algebras describing intertwiners of certain representations, such as the Tem\-per\-ley--Lieb category, $\fgl_n$-webs \cite{Kuperberg_SpidersRank2_1996,CKM_WebsQuantumSkew_2014}, or the Brauer category \cite{LZ_BrauerCategoryInvariant_2015};
  \item \emph{diagrammatic categorification}, such as the 2-Kac--Moody algebra $\cU(\fg)$ of Khovanov--Lauda \cite{KL_CategorificationQuantum$sln$_2010,KL_DiagrammaticApproachCategorification_2009} and Rouquier \cite{Rouquier_2KacMoodyAlgebras_2008} categorifying the quantum group $U_q(\fg)$ (where $\fg$ is a simple complex Lie algebra), or the Elias--Williamson category $\cH(W)$ \cite{EW_HodgeTheorySoergel_2014} categorifying the Hecke algebra $H_q(W)$ (where $W$ is a Coxeter group);
  \item \emph{diagrammatic supercategorification}, such as the 2-Kac--Moody superalgebra $\cS\cU(\fg)$ \cite{BE_SuperKacMoody_2017} conjecturally categorifying the covering quantum group $U_{q,\pi}(\fg)$ \cite{CHW_QuantumSupergroupsFoundations_2013}, or the graded-2-category $\gfoam_d$ of graded $\glt$-foams \cite{SV_OddKhovanovHomology_2023} (discussed below).
\end{itemize}
While the third family of examples is less classical than the other two, it is representative of a general trend: newly discovered diagrammatic algebras (here $\cS\cU(\fg)$) often arise as variations, or ``deformations'', of classical diagrammatic algebras (here $\cU(\fg)$), in the loose sense that the presentation of $\cS\cU(\fg)$ closely resembles that of $\cU(\fg)$.

Unfortunately, diagrammatic algebras are hard to study.
In particular, finding a hom-basis, i.e.\ a basis for each hom-space, is hard.
Often, a careful study of the presentation leads to a candidate hom-basis, shown to generate hom-spaces---the hard part is to show linear independence.
A classical solution is to find a concrete faithful representation of the diagrammatic algebra, extrapolating linear independence of the candidate hom-basis from the linear independence of its image.
However, in the higher context, finding such a concrete representation is a difficult task---assuming it even exists.
To give an example, the basis theorem for 2-Kac--Moody algebras $\cU(\fg)$ was originally proved in the restricted case $\fg=\fsl_n$ by Khovanov--Lauda \cite{KL_CategorificationQuantum$sln$_2010} using a representation on the cohomology of flag varieties. The general proof for any type $\fg$ was only sketched a decade later, by work of Webster \cite{Webster_UnfurlingKhovanovLaudaRouquierAlgebras_2018}.
Outside of the case $\fg=\fsl_2$ \cite{BK_OddGrassmannianBimodules_2026}, the analogue statement for 2\nbd-Kac--Moody superalgebras remains conjectural.

As one can expect, many subsequent questions about these diagrammatic algebras rely on knowing that hom-spaces have the expected dimension.
For instance, if the basis conjecture holds for 2-Kac--Moody superalgebras, it would follow that they do indeed categorify covering quantum groups---this remains a conjecture outside special cases \cite{EL_OddCategorification$U_qmathfraksl_2$_2016,KKO_SupercategorificationQuantumKacMoody_2013,KK_CategorificationHighestWeight_2012}.

The above discussion suggests the following:
\begin{quote}
  \textsc{Question:} \emph{Are there general techniques to study presentations of diagrammatic algebras, in particular to find hom-bases?}
\end{quote}
Better still, we may hope to implement these techniques on a computer.
In classical algebra, these questions and related ones have a long history, which includes Gröbner bases and Buchberger's algorithm for commutative algebras \cite{Buchberger_AlgorithmFindingBasis_2006}, Shirshov's work on PBW-bases in Lie algebras \cite{Shirshov_AlgorithmicProblemsLie_2009}, and Bokut's composition lemma \cite{Bokut_ImbeddingsSimpleAssociative_1976} and Bergmann's diamond lemma for associative algebras \cite{Bergman_DiamondLemmaRing_1978}.
Each of these approaches can be seen as an instance of (classical) \emph{rewriting theory}: the study of presented (classical) algebraic structures from an algorithmic perspective.
By its very nature, the rewriting approach is well-suited to be implemented on a computer. It also gives an abstract tool to find bases, and more generally to find explicit generators of higher relations, known as \emph{syzygies}; see \cite{GHM_ConvergentPresentationsPolygraphic_2019} for a modern perspective.

Recent decades have seen the emergence of \emph{higher rewriting theory}, an analogue of classical rewriting theory for 2-categories and $n$-categories; see \cref{subsec:hrt_state_of_the_art} for a discussion on the literature.
Unfortunately, the application of the theory to diagrammatic algebras has had limited success, despite some notable attempts \cite{%
  Alleaume_RewritingHigherDimensional_2018,%
  Dupont_RewritingModuloIsotopies_2021,%
  Dupont_RewritingModuloIsotopies_2022,%
  DEL_SuperRewritingTheory_2024%
} in that direction in recent years; we discuss them and relations with our work in \cref{subsubsec:intro_literature}.

\begin{table}
  \centering
  \begin{tabular}{@{}c@{\hspace{50pt}}c@{}}
    \emph{categorical property} & \emph{topological interpretation}
    \\\midrule
    interchange law & planar rectilinear isotopies
    \\
    pivotality & planar isotopies
    \\
    symmetric structure & rectilinear (non-planar) isotopies
    \\
    symmetric pivotality & (non-planar) isotopies
  \end{tabular}
  \captionsetup{justification=centering,width=1\textwidth}
  \caption{categorical properties of 2-categories and their topological interpretations as string diagrams; see
  \cite{%
    JS_GeometryTensorCalculus_1991,
    FY_BraidedCompactClosed_1989,
    JS_GeometryTensorCalculus_,
    KL_CoherenceCompactClosed_1980
  } for (partial) proofs and \cite{Selinger_SurveyGraphicalLanguages_2011} for a review.
  }
  \label{tab:cat_prop_topo_interpretation}
\end{table}

In this work, we lay down the foundations of a new rewriting theory:
\begin{quote}
  \textsc{Linear Gray (or Higher) rewriting theory:} \emph{A rewriting theory for linear higher categorical structures, sufficiently flexible to apply to diagrammatic algebras, and in particular prove basis theorems.}
\end{quote}
Some difficulties are theoretical in nature, as the interplay between linearity and higher rewriting leads to unexpected behaviours (see \cref{subsubsec:intro_independent_contextualized_branchings}).
This forces some conceptual shifts, such as the notion of tamed congruence discussed below.
Other difficulties are practical: indeed, diagrammatic algebras do not lend themselves easily to terminating algorithms. This leads to ideas such as rewriting modulo, Gray rewriting, or context-dependent rewriting.
We discuss these aspects below; a more in-depth overview is given in \cref{subsec:intro_example}. The reader new to rewriting theory may wish to first read the state of the art (\cref{subsec:hrt_state_of_the_art}).

\medbreak

\noindent$\bullet$ \textsc{Tamed congruence.}\quad Rather than confluence, our theory is built on the notion of \emph{$\succ$-tamed congruence} (\cref{defn:ARSM_tame_congruence}), for $\succ$ a partial order; see \cref{eq:intro_tamed_congruence}.
Related ideas can already be found in the work of Winkler and Buchberger \cite{WB_CriterionEliminatingUnnecessary_1986} and Bergman \cite{Bergman_DiamondLemmaRing_1978}, as well as in van Oostrom's notion of decreasingness \cite{vanOostrom_ConfluenceDecreasingDiagrams_1994}.

\begin{figure}
  \centering
  \begin{IEEEeqnarray*}{CcC}
    \tikzpic{
      \node (L) at (0,0) {$\bullet$};
      \node (T) at (1.5,1) {$\cdot$};
      \node (B) at (1.5,-1) {$\cdot$};
      \node (R) at (3,0) {$\cdot$};
      \draw[->] (L) to[out=60,in=180] node[above]{$f$} (T);
      \draw[->] (L) to[out=-60,in=180] node[below]{$g$} (B);
      \draw[->] (T) to[out=0,in=120] (R);
      \draw[->] (B) to[out=0,in=-120](R);
    }[scale=.7]
    &\mspace{120mu}&
    \tikzpic{
      \node (C) at (0,0) {$\bullet$};
      %
      \node[inner sep=0pt] (A1) at (2,1) {$\cdot$};
      \node[inner sep=0pt] (A2) at (2+.2,1-.4) {$\cdot$};
      \node[inner sep=0pt] (A3) at (2-.6,1-2*.4) {$\cdot$};
      \node[inner sep=0pt] (A4) at (2+1.2,1-3*.4) {$\cdot$};
      \node[inner sep=0pt] (A5) at (2-.2,1-4*.4) {$\cdot$};
      \node[inner sep=0pt] (A6) at (3,-1) {$\cdot$};
      \draw (A1) to (A2);
      \draw[
        dotted,
        dash pattern=on 3pt off 3pt, 
      ] (A2) to (A3);
      \draw[
        dotted,
        dash pattern=on 3pt off 3pt, 
      ] (A3) to (A4);
      \draw[
        dotted,
        dash pattern=on 3pt off 3pt, 
      ] (A4) to (A5);
      \draw (A5) to (A6);
      \draw[->] (C) to[out=60,in=180] node[above]{$f$} (A1);
      \draw[->] (C) to[out=-60,in=180] node[below]{$g$} (A6);
      \draw[dotted,thick] (.5,2) to (.5,-1.5);
      \node at (1,2) {$\bullet\succ$};
    }[scale=.7]
    \\*\nonumber
    \text{\small confluence}
    &&
    \text{\small $\succ$-tamed congruence}
  \end{IEEEeqnarray*}
  \caption{Schematic of confluence and tamed congruence. In the latter case, horizontal positions are used to suggest relative ordering with respect to $\succ$, reading from left to right.}
  \label{eq:intro_tamed_congruence}
\end{figure}

\medbreak

\noindent$\bullet$ \textsc{Rewriting modulo.}\quad Inspired by the work of Dupont \cite{Dupont_RewritingModuloIsotopies_2022}, we allow rewriting modulo another set of rules.
As pointed out above, certain defining relations of a diagrammatic algebra often capture a categorical property, which in turn admits a topological interpretation in terms of string diagrams (see \cref{tab:cat_prop_topo_interpretation}).
For instance, in a pivotal category string diagrams are best understood up to planar isotopies. When rewriting modulo, we partition the set of relations into two sets: \emph{oriented} relations, thought of as specific to the given diagrammatic algebra, and \emph{unoriented} relations, thought of as intrinsic to the underlying categorical structure.
When working modulo, branchings themselves need to be understood modulo.
In practice, this is done via the naturality conditions associated with the chosen categorical structure (see \cref{subsubsec:intro_classify_branching}).

We also allow the modulo data to contain scalars.
For instance, the 2-Kac--Moody superalgebra $\cS\cU(\fg)$ is a super-2-category \cite{BE_MonoidalSupercategories_2017}, in the sense that the interchange only holds up to sign:
\begin{equation*}
  \tikzpic{
    \draw (0,2) node[above] {\scriptsize $\stdmor'$} to
      node[black_dot,pos=.3] {}
      node[right,pos=.3] {\scriptsize $a$}
      (0,0) node[below] {\scriptsize $\stdmor$};
    \draw (1,2) node[above] {\scriptsize $\stdmorbis'$} to
      node[black_dot,pos=.7] {}
      node[right,pos=.7] {\scriptsize $b$}
      (1,0) node[below] {\scriptsize $\stdmorbis$};
  }[scale=0.6]
  \;=\;(-1)^{p(a)p(b)}
  \tikzpic{
    \draw (0,2) node[above] {\scriptsize $\stdmor'$} to
      node[black_dot,pos=.7] {}
      node[right,pos=.7] {\scriptsize $a$}
      (0,0) node[below] {\scriptsize $\stdmor$};
    \draw (1,2) node[above] {\scriptsize $\stdmorbis'$} to
      node[black_dot,pos=.3] {}
      node[right,pos=.3] {\scriptsize $b$}
      (1,0) node[below] {\scriptsize $\stdmorbis$};
  }[scale=0.6]
\end{equation*}
Here $p(a)$ and $p(b)$ denote the extra data of a parity for $a$ and $b$. In that situation, one would still like to think of string diagrams up to planar rectilinear isotopies, although this equivalence only holds up to sign.

\medbreak

\noindent$\bullet$ \textsc{Gray rewriting.}\quad In fact, our theory is not based on linear (strict) 2-categories, but on \emph{linear sesquicategories} (\cref{sec:higher_structure}), i.e., linear 2-categories without the interchange law.
This is both motivated by examples such as $\cS\cU(\fg)$ and by theoretical considerations; see \cref{subsubsec:intro_independent_contextualized_branchings}.
In analogy with the terminology of \cite{FM_RewritingGrayCategories_2022}, we call our theory \emph{linear Gray rewriting modulo}.

\medbreak

\noindent$\bullet$ \textsc{Context-dependent rewriting.}\quad As pointed out in \cite{Alleaume_RewritingHigherDimensional_2018}, terminating reduction algorithms in diagrammatic algebras often depend on context.
For instance, say we rewrite modulo the interchange law with the following rewriting rule:
\begin{gather*}
  \ssr\colon
  \tikzpic{
    \minicap[0][1]\ministrand[2][1]
    \minicup\ministrand[2][0]
  }[scale=.4]
  \;\to\;
  \tikzpic{
    \ministrand[1][1]\minicap[2][1]
    \ministrand[1][0]\minicup[2][0]
  }[scale=.4]\;.
\end{gather*}
If we allow the rewriting rule $\Gamma[\ssr]$ for every context $\Gamma$, our rewriting system does not terminate:
\begin{gather*}
  \tikzpic{
    \minicap[2][2]
    \minicap[0][1]\ministrand[2][1]\ministrand[3][1]
    \minicup\ministrand[2][0]\ministrand[3][0]
  }[scale=.4]
  \;\to\;
  \tikzpic{
    \minicap[0][2][{xscale=2}]
    \ministrand[0][1]\minicap[.5][1]\ministrand[2][1]
    \ministrand[0][0]\minicup[.5][0]\ministrand[2][0]
  }[scale=.4]
  \;\to\;
  \tikzpic{
    \minicap[0][2]
    \ministrand[0][1]\ministrand[1][1]\minicap[2][1]
    \ministrand[0][0]\ministrand[1][0]\minicup[2][0]
  }[scale=.4]
  \;\sim_{\mathrm{interchange}}\;
  \tikzpic{
    \minicap[2][2]
    \minicap[0][1]\ministrand[2][1]\ministrand[3][1]
    \minicup\ministrand[2][0]\ministrand[3][0]
  }[scale=.4]\;.
\end{gather*}
To recover termination, we should allow
\tikzpic{\minicap[2][1]\minicup[2][0]}[scale=.3]
to slide ``outside a cap'', but not ``inside a cap''; note that this condition depends on the context $\Gamma$.
To cover such situations, we introduce \emph{context-dependent rewriting}; see \cref{subsubsec:HLRSM_definition}.

\medbreak

Finally, we apply the theory to a diagrammatic algebra arising from categorification and quantum topology, the graded-2-category $\gfoam_d$ of graded $\glt$-foams (\cref{sec:rewriting_foam}):

\begin{bigtheorem*}[\cref{thm:foam_basis_theorem}]
  \label{thm:main_thm_foam}
  The graded-2-category $\gfoam_d$ has the expected hom-basis.
\end{bigtheorem*}

A graded-2-category \cite{SV_OddKhovanovHomology_2023} is analogous to a super-2-category, where the interchange law only holds up to scalar.
At the time of writing, our approach is the only known approach to this result; it gives the first application of higher rewriting theory to categorification and quantum topology.
The result is of independent interest, as the higher-representation-theoretic construction of odd Khovanov homology in \cite{SV_OddKhovanovHomology_2023} relies on it.

\subsubsection{Organisation}

The rest of this introduction consists of: a state of the art (\cref{subsec:hrt_state_of_the_art}), which also serves as an introduction to rewriting theory for non-experts; an elementary example highlighting some aspects of the theory (\cref{subsec:intro_example}); a comparison with the literature (\cref{subsubsec:intro_literature}); and some perspectives for future directions of research (\cref{subsec:perspectives}).

Apart from this introduction, the article has three sections.
The short \cref{sec:higher_structure} defines the necessary categorical structures at play; in particular, a suitable notion of presentation for graded-2-categories given by \emph{linear Gray polygraphs}.
\Cref{sec:foundation_rewriting} is the heart of the paper: it develops \emph{Linear Gray rewriting modulo}, starting from first principles.
Finally, we prove the basis theorem for graded $\glt$-foams in \cref{sec:rewriting_foam}, exhibiting all the techniques developed in \cref{sec:foundation_rewriting}.

More precisely, \cref{subsec:coherence_modulo,subsec:abstract_rewriting_modulo,subsec:higher_rewriting_modulo} develop higher (non-linear) rewriting modulo. After reading these subsections, the reader is equipped to read \cref{subsec:rewriting_foam_coherence_iso}, which deals with coherence of the modulo data in graded $\glt$-foams.
Higher linear rewriting modulo is then developed in \cref{subsec:linear_rewriting_modulo,subsec:higher_linear_rewriting_modulo} and applied in \cref{subsec:rewriting_foam_confluence_modulo_iso}.
Finally, \cref{subsec:rewriting_summary} gives a blueprint on how to apply higher linear rewriting modulo and \cref{subsec:notes} consists of further notes and further comparison with the literature (see also \cref{subsubsec:intro_literature}).

\subsubsection{Acknowledgments}

The author thanks Sigiswald Barbier, Jon Brundan, Ben Elias, Yves Guiraud, Louis-Hadrien Robert, Pedro Vaz and Emmanuel Wagner for their interest and for many helpful questions and comments, and an anonymous referee for their thorough reading and valuable suggestions.
Some string diagrams were done using \href{https://smimram.github.io/satex/}{SaTeX}.
The author was supported by the Fonds de la Recherche Scientifique--FNRS
under the Aspirant Fellowship FC 38559 and by the Max Planck Institute for Mathematics.

\subsection{State of the art}
\label{subsec:hrt_state_of_the_art}

\subsubsection{What is rewriting theory?}
\label{subsubsec:intro_ov_rw}

The word problem for monoids asks the following question: given a presented monoid $G$ with generators in the set $\sX$, is there an algorithm that decides whether two words with letters in $\sX$ are equal as elements of $G$?
While known to be undecidable in general \cite{Post_RecursiveUnsolvabilityProblem_1947,Markoff_ImpossibilityCertainAlgorithms_1947}, one can hope to solve the word problem in practical cases.
\emph{Rewriting theory} suggests the following method, which we illustrate with the symmetric group on three strands (understood as a monoid), defined using its Coxeter presentation:
\[
  \fG_3=\left\langle \sigma,\tau\mid \sigma\sigma=1,\tau\tau=1,\sigma\tau\sigma=\tau\sigma\tau\right\rangle.
\]
A \emph{rewriting system} consists of a choice of orientations on the defining relations. For instance:
\begin{equation}
  \label{eq:Sthree_rw}
  \sP_{\fG_3}\coloneqq (\sX,\sR),
  \quad\text{with}\quad
  \sX=\{\sigma,\tau\}
  \quad\an\quad
  \sR=\{\sigma\sigma\to 1,\tau\tau\to 1,\sigma\tau\sigma\to \tau\sigma\tau\}.
\end{equation}
Let $\sX^*$ denote the set of words in letters in $\sX$.
Any rewriting system defines a non-deterministic algorithm, where for each oriented relation ${A\to B}$ and words $x,y\in\sX^*$, the algorithm may perform the reduction $xAy\to xBy$, but not the reduction $xBy\to xAy$.
In the terminology of rewriting theory, $xAy\to xBy$ is called a \emph{rewriting step}, and a successive composition of rewriting steps, denoted $a\overset{*}{\to}b$, is called a \emph{rewriting sequence}.

A pair of co-initial rewriting sequences $(f\colon a\overset{*}{\to} b,g\colon a\overset{*}{\to} b')$ is called a \emph{branching}, and a pair of co-terminal rewriting sequences $(f'\colon b\overset{*}{\to} c,g'\colon b'\overset{*}{\to} c)$ is called a \emph{confluence}.
A branching that admits a confluence is said to be \emph{confluent}. This is illustrated in the diagrams below, where plain arrows (resp.\ dotted arrows) denote branchings (resp.\ confluences):
\begin{IEEEeqnarray*}{CcCcC}
  \begin{tikzcd}[ampersand replacement=\&,cramped,row sep=.7em]
    \& b \\
    a \& \\
    \& {b'}
    \arrow["*"{description},"f", curve={height=-9pt}, from=2-1, to=1-2,]
    \arrow["*"{description},"{g}"', curve={height=9pt}, from=2-1, to=3-2,]
  \end{tikzcd}
  &&
  \begin{tikzcd}[ampersand replacement=\&,cramped,row sep=.7em]
    \& b \\
    \&\& c \\
    \& {b'}
    \arrow["*"{description},"f'", curve={height=-9pt}, from=1-2, to=2-3,densely dotted,]
    \arrow["*"{description},"{g'}"'{shift={(-.05,.05)}}, curve={height=9pt}, from=3-2, to=2-3,densely dotted]
  \end{tikzcd}
  &&
  \begin{tikzcd}[ampersand replacement=\&,cramped,row sep=.7em]
    \& b \\
    a \&\& c \\
    \& {b'}
    \arrow["*"{description},"f", curve={height=-9pt}, from=2-1, to=1-2,]
    \arrow["*"{description},"{g}"', curve={height=9pt}, from=2-1, to=3-2,]
    \arrow["*"{description},"f'", curve={height=-9pt}, from=1-2, to=2-3,densely dotted,]
    \arrow["*"{description},"{g'}"'{shift={(-.05,.05)}}, curve={height=9pt}, from=3-2, to=2-3,densely dotted]
  \end{tikzcd}
  \\*
  \text{\small a branching}
  &\mspace{70mu}&
  \text{\small a confluence}
  &\mspace{70mu}&
  \text{\small a confluent branching}
\end{IEEEeqnarray*}
\Cref{fig:intro_branching} gives examples of confluent branchings in $\sP_{\fG_3}$.

To solve the word problem, we ask two key properties: \emph{termination}, which postulates that any rewriting sequence terminates, and \emph{confluence}, which postulates that every branching is confluent.
A rewriting system which is both terminating and confluent is called \emph{convergent}. In that case, the non-deterministic algorithm always produces an output, and this output is always the same. Given a word as input, we call the output its \emph{normal form}. One can readily show that under convergence, two words are equal in the associated monoid if and only if they have identical {normal forms}. Assuming one can tell whether a word is a normal form, this provides a solution to the word problem.

\begin{figure}
  \begin{gather*}
    \begin{tikzpicture}[scale=.5,xscale=.8,yscale=.8]
      \node (B) at (-3.5,0) {$\sigma\sigma\tau\sigma$};
      \node (B1) at (0,3) {$\tau\sigma$};
      \node (B2) at (0,-3) {$\sigma\tau\sigma\tau$};
      \node (C) at (5,0) {$\tau\sigma$};
      \node (C2) at (3.5,-2) {$\tau\sigma\tau\tau$};
      \draw[->,bend left=20pt] (B) to (B1);
      \draw[->,bend right=20pt] (B) to (B2);
      \draw[bend left=20pt,densely dotted,double distance=1.5pt] (B1) to (C);
      \draw[->,bend right=10pt,densely dotted] (B2) to (C2);
      \draw[->,bend right=10pt,densely dotted] (C2) to (C);
    \end{tikzpicture}
    \qquad
    \begin{tikzpicture}[scale=.5,xscale=.8,yscale=.8]
      \node (B) at (-3.5,0) {$\sigma\tau\sigma\sigma$};
      \node (B1) at (0,3) {$\sigma\tau$};
      \node (B2) at (0,-3) {$\tau\sigma\tau\sigma$};
      \node (C) at (5,0) {$\sigma\tau$};
      \node (C2) at (3.5,-2) {$\tau\tau\sigma\tau$};
      \draw[->,bend left=20pt] (B) to (B1);
      \draw[->,bend right=20pt] (B) to (B2);
      \draw[bend left=20pt,densely dotted,double distance=1.5pt] (B1) to (C);
      \draw[->,bend right=10pt,densely dotted] (B2) to (C2);
      \draw[->,bend right=10pt,densely dotted] (C2) to (C);
    \end{tikzpicture}
    \qquad
  \begin{tikzpicture}[scale=.5,xscale=.8,yscale=.8]
    \node (B) at (-3.5,0) {$\sigma\tau\sigma\tau\sigma$};
    \node (B1) at (0,3) {$\sigma\tau\tau\sigma\tau$};
    \node (B2) at (0,-3) {$\tau\sigma\tau\tau\sigma$};
    \node (C) at (5,0) {$\tau$};
    \node (C1) at (3.5,2) {$\sigma\sigma\tau$};
    \node (C2) at (3.5,-2) {$\tau\sigma\sigma$};
    \draw[->,bend left=20pt] (B) to (B1);
    \draw[->,bend right=20pt] (B) to (B2);
    \draw[->,bend left=10pt,densely dotted] (B1) to (C1);
    \draw[->,bend left=10pt,densely dotted] (C1) to (C);
    \draw[->,bend right=10pt,densely dotted] (B2) to (C2);
    \draw[->,bend right=10pt,densely dotted] (C2) to (C);
  \end{tikzpicture}
  \end{gather*}

  \caption{Critical branchings of $\sP_{\fG_3}$.}
  \label{fig:intro_branching}
\end{figure}

It remains to show that $\sP_{\fG_3}$ is convergent. Termination is not hard: it can be shown using a suitable partial order. On the other hand, confluence must in principle be checked for \emph{every} branching. A branching $(f,g)$ for which $f$ and $g$ are rewriting \emph{steps} (in contrast to rewriting \emph{sequences}) is said to be \emph{local}.
For instance, the three branchings in \cref{fig:intro_branching}, depicted in plain arrows, are local.
\emph{Newman's lemma} (\cref{lem:ARSM_newmann_lemma}) states that assuming termination, confluence follows from confluence of local branchings.

\begin{figure}
  \begin{IEEEeqnarray*}{CcC}
    \begin{tikzpicture}[scale=.5,xscale=.8,yscale=.8]
      \node (B) at (-5,0) {$(\sigma\sigma)(\tau\tau)$};
      \node (B1) at (0,3) {$(1)(\tau\tau)$};
      \node (B2) at (0,-3) {$(\sigma\sigma)(1)$};
      \node (C) at (5,0) {$(1)(1)$};
      \draw[->,bend left=20pt] (B) to (B1);
      \draw[->,bend right=20pt] (B) to (B2);
      \draw[->,bend left=20pt,densely dotted] (B1) to (C);
      \draw[->,bend right=20pt,densely dotted] (B2) to (C);
    \end{tikzpicture}
    &\mspace{80mu}&
    \begin{tikzpicture}[scale=.5,xscale=1.2,yscale=.8]
      \node (B) at (-3.5,0) {$x(\sigma\tau\sigma\sigma)y$};
      \node (B1) at (0,3) {$x(\sigma\tau)y$};
      \node (B2) at (0,-3) {$x(\tau\sigma\tau\sigma)y$};
      \node (C) at (5,0) {$x(\sigma\tau)y$};
      \node (C2) at (4,-2.5) {$x(\tau\tau\sigma\tau)y$};
      \draw[->,bend left=20pt] (B) to (B1);
      \draw[->,bend right=20pt] (B) to (B2);
      \draw[bend left=20pt,densely dotted,double distance=1.5pt] (B1) to (C);
      \draw[->,bend right=10pt,densely dotted] (B2) to (C2);
      \draw[->,bend right=10pt,densely dotted] (C2) to (C);
    \end{tikzpicture}
  \end{IEEEeqnarray*}
  \caption{An independent and a contextualized branching in $\sP_{\fG_3}$, for any words $x,y\in\sX^*$.}
  \label{fig:intro_indep_context_branching_monoid}
\end{figure}

Still, many local branchings remain, as the ones described in \cref{fig:intro_indep_context_branching_monoid}.
Each of them admits a canonical confluence.
The first branching is an \emph{independent branching}: intuitively, it consists of two rewriting steps that do not interact with each other.
The second branching is a \emph{contextualized branching}: it is of the form $x(f,g)y$, for $(f,g)$ the first local branching in \cref{fig:intro_branching}, and the confluence of $x(f,g)y$ is canonically induced from the confluence of $(f,g)$.
A \emph{critical branching} is a local branching which is neither an independent branching nor a contextualized branching.
As suggested by our discussion, confluence of local branchings follows from confluence of critical branchings.
The latter constitute the minimal amount of computations one needs to perform; the rest follows from general considerations.
The reader may convince themself that \cref{fig:intro_branching} gives a complete list of critical branchings in $\sP_{\fG_3}$. It follows that the rewriting system $\sP_{\fG_3}$ is both terminating and locally confluent, and hence convergent by Newman's lemma.

\subsubsection{Polygraphs}
\label{subsubsec:intro_ov_polygraph}

What is a presentation of a (strict and small) $n$-category? As a first example, consider again the rewriting system $\sP_{\fG_3}=(\sX,\sR)$ defined in \eqref{eq:Sthree_rw}.
Viewing a monoid as a category with a single object $\{*\}$, the set $\sX$ consists of generating 1-cells with source and target the object $\{*\}$. We encapsulate the latter fact with (trivial) source and target maps $s_0,t_0\colon\sX\to\{*\}$.
In this perspective, the set of words $\sX^*$ is the free category generated by $\sX$.
Similarly, the set $\sR$ consists of generating 2-cells, and we define maps $s_1,t_1\colon \sR\to\sX^*$ setting $s_1(r)=A$ and $t_1(r)=B$ for each oriented relation $r\colon A\to B$.
Reformulated in this way, $\sP_{\fG_3}$ defines the data of a \emph{2-polygraph}.
More generally, an $n$-category can be presented by an \emph{$(n+1)$-polygraph}, with generating $(k+1)$-cells $\sP_{k+1}$ defined on the free $k$-category $\sP_k^*$ generated by the lower cells:
\begin{IEEEeqnarray*}{cCc}
  \begin{tikzcd}[ampersand replacement=\&,cramped]
    {\{*\}^*} \& \sX^*\\
    {\{*\}} \& \sX \& {\sR}
    \arrow[from=2-1, to=1-1,equals]
    \arrow[shift right, from=2-2, to=1-1,"s_0"']
    \arrow[shift left, from=2-2, to=1-1,"t_0"]
    \arrow[from=2-2, to=1-2]
    \arrow[shift right, from=2-3,to=1-2,"s_1"']
    \arrow[shift left, from=2-3, to=1-2,"t_1"]
  \end{tikzcd}
  &\mspace{50mu}&
  \begin{tikzcd}[ampersand replacement=\&,cramped]
    {\sP_0^*} \& \ldots \& {\sP_{n-1}^*} \& {\sP_{n}^*} \\
    {\sP_0} \& \ldots \& {\sP_{n-1}} \& {\sP_n} \& {\sP_{n+1}}
    \arrow[from=2-1, to=1-1,equals]
    \arrow[shift right, from=2-2, to=1-1,"s_0"']
    \arrow[shift left, from=2-2, to=1-1,"t_0"]
    \arrow[shift right, from=2-3, to=1-2,"s_{n-2}"'{xshift=-5pt}]
    \arrow[shift left, from=2-3, to=1-2,"t_{n-2}"{xshift=7pt}]
    \arrow[from=2-3, to=1-3]
    \arrow[shift right, from=2-4, to=1-3,"s_{n-1}"'{xshift=-5pt}]
    \arrow[shift left, from=2-4, to=1-3,"t_{n-1}"{xshift=7pt}]
    \arrow[from=2-4, to=1-4]
    \arrow[shift right, from=2-5, to=1-4,"s_n"']
    \arrow[shift left, from=2-5, to=1-4,"t_n"]
  \end{tikzcd}
  \\*[1ex]
  \text{\small the 2-polygraph $\sP_{\fG_3}$ presenting $\fG_3$}
  &&
  \text{\small an $(n+1)$-polygraph $\sP$ presenting an $n$-category}
\end{IEEEeqnarray*}
The $n$-category presented by an $(n+1)$-polygraph $\sP$ is obtained by quotienting $\sP_n^*$ by the relation $s_n(r)=t_n(r)$ for each $r\in\sP_{n+1}$.

Polygraphs were first introduced by Street \cite{Street_LimitsIndexedCategoryvalued_1976}, under the name of \emph{computads}; the term \emph{signatures} also appears in the literature.
Polygraphs were independently introduced by Burroni \cite{Burroni_HigherdimensionalWordProblems_1993} to study generalizations of the word problem---the terminology is by now standard in the rewriting literature.

\subsubsection{Higher rewriting}
\label{subsubsec:intro_ov_higher_rw}

The rewriting theory associated with an $(n+1)$-polygraph is called \emph{$(n+1)$-dimensional rewriting theory}.
Higher rewriting theory has attracted increasing interest in the last two decades
\cite{
  Lafont_AlgebraicTheoryBoolean_2003,
  Guiraud_TerminationOrdersThreedimensional_2006,
  GM_HigherdimensionalCategoriesFinite_2009,
  GM_HigherdimensionalNormalisationStrategies_2012,
  Mimram_3dimensionalRewritingTheory_2014,
  GM_PolygraphsFiniteDerivation_2018,
}, with as a major application the construction of a homology theory for $\omega$-categories, based on polygraphic resolutions
\cite{
  Mimram_ComputingCriticalPairs_2010,%
  Metayer_ResolutionsPolygraphs_2003,%
  Guetta_HomologyCategoriesPolygraphic_2021,%
  LMW_FolkModelStructure_2010,%
  LM_PolygraphicResolutionsHomology_2009,%
}; see the recent monograph on the subject \cite{ABG+_PolygraphsRewritingHigher_2025}.
The 1-dimensional case corresponds to rewriting in sets; it is also known as \emph{abstract rewriting}. As we have seen, 2-dimensional rewriting theory corresponds to rewriting in categories, and monoids in particular. Correspondingly, three-dimensional rewriting is rewriting in 2-categories.

Recall the notion of a contextualized branching from \cref{fig:intro_indep_context_branching_monoid}.
Expressed in the terminology of 2-polygraphs, contextualization is the idea that given a 2-cell $r\colon A\to B$ and two 1-cells $x$ and $y$ suitably composable with the source and target of $A$ (or $B$), there exists a 2-cell
\[x[r]y\colon xAy\to xBy,\]
defined using left and right compositions.
When doing rewriting in 2-categories, a relation $r\colon A\to B$ is a generating 3-cell, and contextualization amounts to first composing horizontally with 1-cells $x$ and $y$, and then vertically with 2-cells $\alpha$ and $\beta$:
\begin{gather*}
  \Gamma[r]\coloneqq
  \satex[scale=1.3]{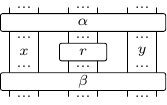}
  \colon\quad
  \satex[scale=1.3]{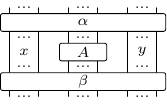}
  \longrightarrow
  \satex[scale=1.3]{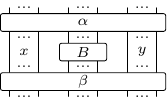}.
\end{gather*}
Here $\Gamma$ is the ``diagram with a hole'' defined by $x$, $y$, $\alpha$ and $\beta$.
Contextualization is one of the main features that differentiate abstract rewriting from higher rewriting.

\subsubsection{Linear rewriting}
\label{subsubsec:intro_ov_linear_rw}

Linear rewriting is the algorithmic study of presented modules.
Given a commutative ring $\Bbbk$, a $\Bbbk$-module $\lM$ is presented by a set $\cB$, together with a set $\lR$ of oriented relations on $\langle\cB\rangle_\Bbbk$, the free $\Bbbk$-module generated by $\cB$.
This is encapsulated by the data of a \emph{linear 1-polygraph}:
\[\begin{tikzcd}[cramped]
  \langle\cB\rangle_\Bbbk &  \lR
  \arrow["{t}", shift left, from=1-2, to=1-1]
  \arrow["{s}"', shift right, from=1-2, to=1-1]
\end{tikzcd}.\]
Similarly, \emph{linear 2-polygraphs} present linear categories (and in particular associative algebras) and linear $(n+1)$-polygraphs present linear $n$-categories.
The associated linear $(n+1)$-dimensional rewriting theory was studied by Guiraud, Hoffbeck and Malbos in the case $n=2$ \cite{GHM_ConvergentPresentationsPolygraphic_2019}, and by Alleaume in the case $n=3$ \cite{Alleaume_RewritingHigherDimensional_2018}; see \cref{subsubsec:intro_literature} for how our approach differs.

The most classical setting for linear rewriting is commutative algebras. For that reason, elements of $\cB$ are called \emph{monomials}. An oriented relation 
$r\colon s(r)\to t(r)\in\lR$ is assumed to be of the form ``rewrite a monomial into a linear combination of monomials'', that is:
\begin{gather}
  \label{eq:left-monomial}
  r\colon b\to_\lR \lambda_1b_1+\ldots+\lambda_nb_n\qquad \text{for }\lambda_1,\ldots,\lambda_n\in\Bbbk,b,b_1,\ldots,b_n\in\cB
\end{gather}
We say that relations in $r$ are \emph{left-monomial}.

A generic rewriting step is of the form
\[\lambda r+v\colon\lambda s(r)+v\dashrightarrow_\lR\lambda t(r) +v\quad \lambda\in\Bbbk\setminus\{0\},v\in\langle\cB\rangle_\Bbbk,\] 
and it is said to be \emph{positive} if the monomial $s(r)\in\cB$ does not appear in the linear decomposition of $v$.
For example, if $\cB=\{a,b,c\}$ and $r\colon a\to b+c\in\lR$ is an oriented relation, then both
\[2a+b=a+(a+b)\dashrightarrow_\lR (b+c)+(a+b)\quad\an\quad 2a+b\to_\lR 2(b+c)+b\]
are rewriting steps, but only the latter is positive; note the use of dashed and plain arrows to distinguish the two.
To avoid rewriting loops, one must restrict to positive rewriting steps:
\begin{gather*}
  0=a-a\dashrightarrow_\lR (b+c)-a\dashrightarrow_\lR (b+c)-(b+c)=0.
\end{gather*}
This restriction constitutes the main difficulty of the linear setting; see \cref{subsubsec:intro_independent_contextualized_branchings}.

We denote $\lR^+$ the set of positive rewriting steps. Similarly to the abstract setting, we can reduce the study of confluence to the study of \emph{local} confluence.
In fact, we can further reduce to local \emph{monomial} confluence, where a \emph{monomial branching} is a branching whose source is a monomial:

\begin{lemma}[linear Newman's lemma]
  \label{lem:intro_linear_newmann}
  Let $\lS=(\cB;\lR)$ be a linear 1-polygraph.
  If $\lR^+$ terminates, then confluence of monomial $\lR^+$-branchings implies confluence. 
\end{lemma}

\subsubsection{Rewriting modulo}
\label{subsubsec:intro_ov_rw_modulo}

Enforcing all relations to be oriented can be too restrictive. Instead, one may wish to rewrite with a set of oriented relations $\aR$\footnote{We use blackboard font to refer to abstract rewriting; this is unrelated to the set of real numbers.}, modulo another set of \emph{un}oriented relations $\aE$.
More precisely, the working data of \emph{abstract rewriting modulo} is given by two 1-polygraphs
\[\begin{tikzcd}[cramped]
  X & \aR
  \arrow["{t}", shift left, from=1-2, to=1-1]
  \arrow["{s}"', shift right, from=1-2, to=1-1]
\end{tikzcd}
\quad\an\quad
\begin{tikzcd}[cramped]
  X & \aE
  \arrow["{t}", shift left, from=1-2, to=1-1]
  \arrow["{s}"', shift right, from=1-2, to=1-1]
\end{tikzcd},\]
defined on the same underlying set $X$.
Denote by $\aE^\top=(\aE\cup\aE^{-1})^*$ the free groupoid generated by $\aE$. Intuitively, relations in $\aE^\top$ are \emph{un}oriented rewriting sequences in $\aE$.
In this context, a \emph{rewriting step modulo} is a composition $e'\circ r\circ e$ with $r\in\aR$ and $e,e'\in \aE^\top$:
\[
\begin{tikzcd}
  \bullet\arrow[r,"e","{\aE}"{subscript},snakecd]& \bullet\arrow[r,"r","{\aR}"{subscript}]& \bullet\arrow[r,"e'","{\aE}"{subscript},snakecd]& \bullet
\end{tikzcd}
\]
In other words, in between a rewriting step in $\aR$ one can apply an arbitrary number of relations in $\aE$, in any direction. The data $\aS=(X;\aR,\aE)$ defines an \emph{abstract rewriting system modulo}, and a rewriting step modulo as above is called an \emph{$\aS$-rewriting step}.

Typically, the modulo data $\aE$ will consist of relations thought as being ``structural'', in the sense of being part of some underlying algebraic structure.
For instance, rewriting in commutative algebras is implicitly rewriting in associative algebras modulo commutativity ($\aE$ consists of relations $xy\to yx$ for all monomials $x$ and $y$).

Rewriting modulo allows an inductive approach to the word problem. Indeed, instead of trying to fit the relations $\aR\sqcup\aE$ into a single convergent rewriting system, one can instead show that $\aR$ is convergent ``modulo $\aE$'' on one hand, and that $\aE$ is convergent on the other hand.

Different variants of rewriting modulo have been developed, often with more restrictive modulo rules than the ones described above
\cite{%
  JL_ChurchRosserPropertiesNormal_2012,
  JK_CompletionSetRules_1986,
  Huet_ConfluentReductionsAbstract_1977,
  PS_CompleteSetsReductions_1981,
  Marche_NormalizedRewritingUnified_1998,
  Viry_RewritingModuloRewrite_1995,
}.
Rewriting modulo is used in \cite{CDM_ConfluenceAlgebraicRewriting_2022} to study confluence in Lawvere theories.
A higher analogue to rewriting modulo was introduced in \cite{DM_CoherentConfluenceModulo_2022}, using the formalism of double categories.
In the linear setting, Dupont extended Alleaume's approach \cite{Alleaume_RewritingHigherDimensional_2018} modulo in order to rewrite modulo pivotality in 2-categories \cite{Dupont_RewritingModuloIsotopies_2022}. Based on his theory, he proposed an approach to the basis problem in 2-Kac--Moody algebras in simply-laced cases \cite{Dupont_RewritingModuloIsotopies_2021}.
An approach to rewriting modulo in super-2-categories was also proposed in \cite{DEL_SuperRewritingTheory_2024}, motivated by the study of the 2-Kac--Moody superalgebra in the case $\fg=\slt$; see \cref{subsubsec:intro_literature} for how our approach differs from \cite{Dupont_RewritingModuloIsotopies_2022,Dupont_RewritingModuloIsotopies_2021,DEL_SuperRewritingTheory_2024}.

\subsubsection{Gray rewriting}
\label{subsubsec:intro_ov_gray_rw}

Starting with $n=3$, the strict and weak notions of an $n$-category start to diverge: while a bicategory is alway equivalent to a 2-category, not every tricategory is equivalent to a 3-category. However, every tricategory is equivalent to a \emph{Gray category} \cite{GPS_CoherenceTricategories_1995}.
In Gray categories, associativity and unitality hold strictly, but the interchange law for 2-morphisms only holds up to certain coherent 3-morphisms, called \emph{interchangers}; see \cref{fig:intro_gray}.

\begin{figure}[t]
  \centering
  
  \begin{equation*}
    \begin{array}{ccccc}
      \tikzpic{
        \draw (0,2) node[above] {\scriptsize $\stdmor'$} to
          node[black_dot,pos=.3] {}
          node[right,pos=.3] {\scriptsize $a$}
          (0,0) node[below] {\scriptsize $\stdmor$};
        \draw (1,2) node[above] {\scriptsize $\stdmorbis'$} to
          node[black_dot,pos=.7] {}
          node[right,pos=.7] {\scriptsize $b$}
          (1,0) node[below] {\scriptsize $\stdmorbis$};
      }[scale=0.6]
      \;=\;
      \tikzpic{
        \draw (0,2) node[above] {\scriptsize $\stdmor'$} to
          node[black_dot,pos=.7] {}
          node[right,pos=.7] {\scriptsize $a$}
          (0,0) node[below] {\scriptsize $\stdmor$};
        \draw (1,2) node[above] {\scriptsize $\stdmorbis'$} to
          node[black_dot,pos=.3] {}
          node[right,pos=.3] {\scriptsize $b$}
          (1,0) node[below] {\scriptsize $\stdmorbis$};
      }[scale=0.6]
      &&
      \tikzpic{
        \draw (0,2) node[above] {\scriptsize $\stdmor'$} to
          node[black_dot,pos=.3] {}
          node[right,pos=.3] {\scriptsize $a$}
          (0,0) node[below] {\scriptsize $\stdmor$};
        \draw (1,2) node[above] {\scriptsize $\stdmorbis'$} to
          node[black_dot,pos=.7] {}
          node[right,pos=.7] {\scriptsize $b$}
          (1,0) node[below] {\scriptsize $\stdmorbis$};
      }[scale=0.6]
      \;\Rrightarrow\;
      \tikzpic{
        \draw (0,2) node[above] {\scriptsize $\stdmor'$} to
          node[black_dot,pos=.7] {}
          node[right,pos=.7] {\scriptsize $a$}
          (0,0) node[below] {\scriptsize $\stdmor$};
        \draw (1,2) node[above] {\scriptsize $\stdmorbis'$} to
          node[black_dot,pos=.3] {}
          node[right,pos=.3] {\scriptsize $b$}
          (1,0) node[below] {\scriptsize $\stdmorbis$};
      }[scale=0.6]
      &&
      \tikzpic{
        \draw (0,2) node[above] {\scriptsize $\stdmor'$} to
          node[black_dot,pos=.3] {}
          node[right,pos=.3] {\scriptsize $a$}
          (0,0) node[below] {\scriptsize $\stdmor$};
        \draw (1,2) node[above] {\scriptsize $\stdmorbis'$} to
          node[black_dot,pos=.7] {}
          node[right,pos=.7] {\scriptsize $b$}
          (1,0) node[below] {\scriptsize $\stdmorbis$};
      }[scale=0.6]
      \;?\;
      \tikzpic{
        \draw (0,2) node[above] {\scriptsize $\stdmor'$} to
          node[black_dot,pos=.7] {}
          node[right,pos=.7] {\scriptsize $a$}
          (0,0) node[below] {\scriptsize $\stdmor$};
        \draw (1,2) node[above] {\scriptsize $\stdmorbis'$} to
          node[black_dot,pos=.3] {}
          node[right,pos=.3] {\scriptsize $b$}
          (1,0) node[below] {\scriptsize $\stdmorbis$};
      }[scale=0.6]
      \\
      \text{3-categories}
      &\mspace{30mu}&
      \text{Gray categories}
      &\mspace{30mu}&
      \text{3-sesquicategories}
    \end{array}
  \end{equation*}

  \caption{In a 3-category, a Gray category or a 3-sesquicategory, the interchange law for 2-morphisms respectively holds strictly, holds weakly via interchangers, or does not hold (a priori).}
  \label{fig:intro_gray}
\end{figure}

In \cite{FM_RewritingGrayCategories_2022}, Forest and Mimram initiated the study of rewriting theory for {weak} $n$-categories, starting with Gray categories.
Their rewriting theory is based on certain higher categorical structures akin to Gray categories, but which do not contain any coherence data for interchange laws; they are known as \emph{$n$-sesquicategories} \cite{Street_CategoricalStructures_1996,FM_RewritingGrayCategories_2022,Araujo_SimpleStringDiagrams_2022}, also called \emph{$n$-precategories} in \cite{FM_RewritingGrayCategories_2022}.

\subsection{An elementary example}
\label{subsec:intro_example}

We illustrate some aspects of our theory using an elementary example.
Further aspects can be found in the main text; in particular, the example does not use context-dependent rewriting.

\subsubsection{Basis from convergence modulo}
\label{subsubsec:intro_linear_rw_mod}

In \cref{subsec:linear_rewriting_modulo}, we combine linear rewriting (\cref{subsubsec:intro_ov_linear_rw}) and rewriting modulo (\cref{subsubsec:intro_ov_rw_modulo}), starting with the notion of a \emph{linear rewriting system modulo} ${\lS=(\cB;\lR,\lE)}$ (\cref{defn:linear-rewriting-system-modulo}).
In this modulo setting, we have an analogue notion of {positive} rewriting steps $\lS^+$ and of the linear Newman's lemma (\cref{lem:intro_linear_newmann}); see \cref{thm:LRSM_tamed_newmann_lemma}.
Denote $\NF_\lS$ the $\Bbbk$-module of normal forms for $\lS^+$, that is, the $\Bbbk$-module consisting of elements on which $\lS^+$ terminates.

\begin{rewritingfact}[\LRSMbasisfromconvergencetheorem{}]
  \label{rewritingfact:basis}
  Let $\lM$ be a $\Bbbk$-module presented by a linear rewriting system modulo $\lS=(\cB;\lR,\lE)$.
  If $\lS^+$ is convergent modulo, then the canonical linear map
  \[\NF_\lS/\brak{\lE}_\Bbbk\to \lM\]
  is an isomorphism.
  In particular, if $B$ is basis for $\NF_\lS/\brak{\lE}_\Bbbk$, then it is a basis for $\lM$.
\end{rewritingfact}

\subsubsection{Superadjunctions}
\label{subsubsec:intro_linear_gray_polygraph}

\begin{bigdefinition*}[\cref{defn:presentation_graded_two_cat}]
  A \emph{presentation} of a super-2-category is a linear Gray polygraph.
\end{bigdefinition*}

As an elementary example, let us define $\miniP$, the ``$\bZ$-linear Gray polygraph of superadjunction''.
We first define its set of objects $\miniP_0=\{*\}$ (a single object), its set of generating 1-cell $\miniP_1=\left\{\;\tikzpic{\ministrand}[yscale=.4]\;\right\}$ (a single 1-cell) and its set of generating 2-cells $\miniP_2=\left\{\;\tikzpic{\minicap}[scale=.5]\,,\,\tikzpic{\minicup}[scale=.5]\;\right\}$ (unit and counit, read from bottom to top).
In addition, the two 2-cells are given a parity
\[p(\;\tikzpic{\minicap}[scale=.5]\;)=p(\;\tikzpic{\minicup}[scale=.5]\;)=1\in\bZ/2\bZ,\]
which extends additively to generic diagrams.
Finally, we define a set of generating rewriting steps $\miniP_3=\miniR_3\sqcup\miniE_3$ where
\begin{gather*}
  \newcommand{\tpscl}{.5}
  \miniR_3=
  \left\{\;
  \tikzpic{
    \ministrand[0][1]\minicap[1][1]
    \minicup\ministrand[2][0]
  }[scale=\tpscl]
  \;\to\;
  \tikzpic{
    \ministrand[0][1]
    \ministrand
  }[scale=\tpscl]
  \;,\;
  \tikzpic{
    \minicap[0][1]\ministrand[2][1]
    \ministrand\minicup[1][0]
  }[scale=\tpscl]
  \;\to\;-\;
  \tikzpic{
    \ministrand[0][1]
    \ministrand
  }[scale=\tpscl]
  \;\right\}
  \an
  \miniE_3=
  \left\{\;
  \begin{array}{c}
  \tikzpic{
    \draw (0,2) to
      node[black_dot,pos=.3] {}
      node[right,pos=.3] {\scriptsize $a$}
      (0,0);
    \draw (1,2) to
      node[black_dot,pos=.7] {}
      node[right,pos=.7] {\scriptsize $b$}
      (1,0);
  }[scale=0.6]
  \to (-1)^{p(a)p(b)}
  \tikzpic{
    \draw (0,2) to
      node[black_dot,pos=.7] {}
      node[right,pos=.7] {\scriptsize $a$}
      (0,0);
    \draw (1,2)  to
      node[black_dot,pos=.3] {}
      node[right,pos=.3] {\scriptsize $b$}
      (1,0);
  }[scale=0.6]
  \\
  \text{\scriptsize for all $a,b\in(\miniP_2)^*$}
  \end{array}
  \;\right\}.
\end{gather*}
We denote $\miniS=(\miniP_0,\miniP_1,\miniP_2;\miniR_3,\miniE_3)$ this data; it defines a \emph{Gray rewriting system modulo}, where $\miniE_3$ is the modulo data.

Denote $\miniC$ the super-2-category presented by $\miniP$.
Superadjunction is the super analogue of the classical notion of adjunction: accordingly, we expect hom-spaces in $\miniC$ to have the same dimension as its classical analogue.
Thanks to \cref{rewritingfact:basis}, we can use rewriting to show that fact. (Of course, one can give a much simpler proof, but this is not our point.)
The main difficulty is to show confluence.

\begin{figure}
  \begin{gather*}
    \label{eq:critical_branching_diag_intro}
    \begin{tikzcd}[ampersand replacement=\&,column sep=large,row sep=0]
      \tikzpic{
        \minicap[2][2]
        \minicap[0][1]\ministrand[2][1]\ministrand[3][1]
        \ministrand\minicup[1][0]\ministrand[3][0]
      }[scale=.4]
      \\
      \&
      -\;\tikzpic{\minicap}[scale=.5]
      \\
      -\;\tikzpic{
        \minicap[0][2]
        \ministrand[0][1]\ministrand[1][1]\minicap[2][1]
        \ministrand\minicup[1][0]\ministrand[3][0]
      }[scale=.4]
      \arrow[from=1-1,to=2-2]
      \arrow[from=3-1,to=2-2]
      \arrow[from=1-1,to=3-1,snakecd]
    \end{tikzcd}
    \qquad\an\qquad
    \begin{tikzcd}[ampersand replacement=\&,column sep=large,row sep=0]
      \tikzpic{
        \ministrand[0][2]\minicap[1][2]\ministrand[3][2]
        \minicup[0][1]\ministrand[2][1]\ministrand[3][1]
        \minicup[2][0]
      }[scale=.4]
      \\
      \&
      \tikzpic{\minicup}[scale=.5]
      \\
      -\;\tikzpic{
        \ministrand[0][2]\minicap[1][2]\ministrand[3][2]
        \ministrand[0][1]\ministrand[1][1]\minicup[2][1]
        \minicup[0][0]
      }[scale=.4]
      \arrow[from=1-1,to=2-2]
      \arrow[from=3-1,to=2-2]
      \arrow[from=1-1,to=3-1,snakecd]
    \end{tikzcd}
  \end{gather*}
  \caption{A set of critical branchings for $\miniP$}
  \label{fig:intro_crit_branchings_miniP}
\end{figure}

Two examples of branchings are given in \cref{fig:intro_crit_branchings_miniP}.
In each case, the first branch applies one of the zigzag rewriting steps, and the other branch uses the interchange law and then applies the other zigzag rewriting step. To ensure confluence, exactly one zigzag relation must carry a sign.
Intuitively, the branchings in \cref{fig:intro_crit_branchings_miniP} are ``critical branchings'' akin to the critical branchings of $\sP_{\fG_3}$ (\cref{fig:intro_branching}), in the sense that their confluence should guarantee confluence of the rewriting system.
Most of the theory is dedicated to justifying this intuition.

\subsubsection{Independent and contextualized branchings}
\label{subsubsec:intro_independent_contextualized_branchings}

Recall from \cref{subsubsec:intro_ov_rw} the notions of independent and contextualized branchings. We have analogous notions in higher rewriting, as already discussed in \cref{subsubsec:intro_ov_higher_rw}; see also \cref{subsubsec:HRSM_independent_branching} and \cref{subsubsec:HRSM_contextualization}.

Consider a contextualized branching $\Gamma[f,g]$; as in classical rewriting, we would like to say that if $(f,g)$ is positively confluent, then so is $\Gamma[f,g]$.
However, \emph{this is not true in general}.
Indeed, positivity is in general \emph{not} preserved by contexts:
\begin{quote}
  \textsc{Warning:} \emph{the contextualization of a positive rewriting step needs not be positive!}
\end{quote}
Consider our running example $\miniS$ and denote $\sim_{\miniE}$ the equivalence relation induced by $\miniE$, that is, equivalence via the superinterchange law; we call it \emph{$\miniE$-congruence}.
Then:
\begin{gather*}
  \newcommand{\tpscl}{.5}
  \tikzpic{
    \minicap[0][1]\ministrand[2][1]\ministrand[3][1]
    \minicup\ministrand[2][0]\ministrand[3][0]
  }[scale=\tpscl]
  \;\not\sim_{\miniE}\;
  \tikzpic{
    \ministrand[0][1]\ministrand[1][1]\minicap[2][1]
    \ministrand[0][0]\ministrand[1][0]\minicup[2][0]
  }[scale=\tpscl]
  \quad\text{ while }\quad
  \tikzpic{
    \minicap[2][2]
    \minicap[0][1]\ministrand[2][1]\ministrand[3][1]
    \minicup\ministrand[2][0]\ministrand[3][0]
  }[scale=\tpscl][(0,1*\tpscl)]
  \;\sim_{\miniE}\;
  \tikzpic{
    \minicap[0][2]
    \ministrand[0][1]\ministrand[1][1]\minicap[2][1]
    \ministrand[0][0]\ministrand[1][0]\minicup[2][0]
  }[scale=\tpscl][(0,1*\tpscl)]\;.
\end{gather*}
More formally, we have two 2-cells $v,w$ and a context $\Gamma$ such that $v\not\sim_{\miniE} w$, while ${\Gamma[v]\sim_{\miniE} \Gamma[w]}$.
In other words, \emph{contextualization does not act freely}.
In particular, if $\lambda r+v$ is a rewriting step, it may be that $\Gamma[s(r)]$ \emph{is} $\miniE$-congruent (up to sign) to a monomial in the linear decomposition of $\Gamma[v]$, even if $s(r)$ \emph{is not} $\miniE$-congruent (up to sign) to a monomial in the linear decomposition of $v$.
That is, it may be that $\lambda\Gamma[r]+\Gamma[v]$ is not positive, even if $\lambda r+v$ is positive.

This gives a theoretical motivation for using 2-sesquicategories as a foundation of our theory, even when rewriting in 2-categories: contextualization \emph{is} free on 2-sesquicategories, and it is the explicit addition of a modulo rule that prevents freeness.
Note that contextualization is free both for associative algebras and associative algebras modulo commutativity; this perhaps explains why the problem had not been recognised before (to the author's knowledge).

The same issue occurs when considering the confluence of independent branchings.
To solve that issue, we use the alternative notion of tamed congruence.

\subsubsection{Tamed congruence}
\label{subsubsec:intro_tamed_congruence}

Recall the linear Newman's lemma (\cref{lem:intro_linear_newmann}).
An analogous statement can be given in the modulo setting.
However, we shall need a more general version, based on the notion of tamed congruence (see \cref{eq:intro_tamed_congruence} and \cref{defn:ARSM_tame_congruence} in the text):

\begin{rewritingfact}[\LRSMtamednewmannlemma{}]
  \label{rewritingfact:tamed_newmann}
  Let $\lS=(\cB;\lR,\lE)$ be a \LRSM{} and $\succ$ a partial order on $\cB$ satisfying some conditions. If $\succ$ is well-founded and every monomial $\lS^+$\nbd-branching is $\succ$\nbd-tamely $\lS$\nbd-congruent, then $\lS^+$ is convergent.
\end{rewritingfact}

In practice, the partial order $\succ$ is the well-founded order used to show termination; in our running example $\miniS$, the partial order $\succ$ compares the number of generating 2-cells.

\medbreak

We can now deal with contextualized branchings:

\begin{rewritingfact}[\HLRSMcontextualizationlemma{}]
  \label{rewritingfact:contextualization_lemma}
  Let $\sS=(\sR,\sE)$ be a higher linear rewriting system modulo satisfying some conditions.
  Let $(f,g)$ be a monomial $\sS^+$\nbd-branching and $\Gamma$ a context.
  If $(f,g)$ is $\sS^+$\nbd-confluent, then $\Gamma[f,g]$ is $\succ$\nbd-tamely $\sS$\nbd-congruent.
\end{rewritingfact}

and independent branchings:

\begin{rewritingfact}[\cref{lem:HLRSM_congruence_independent_branching}]
  \label{rewritingfact:HLRSM_congruence_independent_branching}
  Let $\sS=(\sR,\sE)$ be a higher linear rewriting system modulo satisfying some conditions.
  Every independent $\sS^+$\nbd-branching is $\succ$\nbd-tamely $\sS$\nbd-congruent.
\end{rewritingfact}

Tamed congruence is also better behaved than confluence.
Namely, it is transitive, in the sense that if the branchings $(f_1,f_2)$ and $(f_2,f_3)$ are tamed congruent, then so is $(f_1,f_3)$ (\cref{lem:tamed_is_transitive}).
As we will see in \cref{sec:rewriting_foam}, transitivity allows to simplify the analysis of branchings.

\subsubsection{How to classify branchings?}
\label{subsubsec:intro_classify_branching}

The arguments of this subsection are abstract; thus we fix an \ARSM{} $\aS=(X;\aR,\aE)$.
Recall that the equivalence relation induced by $\aE$ on $X$ is called $\aE$-congruence; correspondingly, an unoriented composition of relations in $\aE$ is called an \emph{$\aE$-congruence}.

Working modulo induces a huge number of branchings: indeed, in between the two branches of a branching, one can apply an arbitrary number of modulo relations.
To classify branchings, we should not only understand elements of $X$ modulo, but also rewriting sequences modulo.
We say that two rewriting sequences $f$ and $f'$ are \emph{$\aE$-congruent} if there exists $\aE$-congruences $e_s$ and $e_t$ such that:
\begin{gather*}
  \begin{tikzcd}[ampersand replacement=\&]
    \cdot \& \cdot \\
    {\cdot} \& {\cdot}
    \arrow[""{name=0, anchor=center, inner sep=0},"*"{description},"f"{yshift=.5ex},"\aS"{subscript}, from=1-1, to=1-2]
    \arrow[""{name=1, anchor=center, inner sep=0},"*"{description},"f'"'{yshift=-.5ex},"\aS"{subscript}, from=2-1, to=2-2]
    \arrow["{e_s}"',snakecd, no head, from=1-1, to=2-1]
    \arrow["{e_t}",snakecd, no head, from=1-2, to=2-2]
  \end{tikzcd}
\end{gather*}
Two branchings $(f,g)$ and $(f',g')$ are \emph{branchwise $\aE$-congruent} if $f$ (resp.\ $g$) is {$\aE$-congruent} to $f'$ (resp.\ $g'$).
It is not difficult to see that:

\begin{rewritingfact}[\ARSMbranchwiseconfluencelemma{}]
  \label{rewritingfact:branchise_confluence_lemma}
  Let $\aS=(\aR,\aE)$ be an \ARSM{}.
  If $(f,g)$ and $(f',g')$ are branchwise ${\aE}$\nbd-cong\-ruent branchings, then $(f,g)$ is confluent if and only if $(f',g')$ is.
\end{rewritingfact}

Hence, the study of confluence can be done up to branchwise $\aE$-congruence.
Consider again our example $\miniS$.
What does it mean to understand the rewriting steps in $\miniR$ up to $\miniE$\nbd-cong\-ruence? The interchange law preserves the set of generating 2-cells; more precisely, any sequence of interchangers between two diagrams $s$ and $t$ induces a canonical bijection between the generating 2-cells in $s$ and the generating 2-cells in $t$.
Hence it makes sense to speak about \emph{the} cap and \emph{the} cup associated with a rewriting step in $\miniS$:

\begin{lemma}[characterization of rewriting steps in $\miniS$]
  \label{lem:intro_characterization_miniS}
  If two rewritings steps in $\miniS$ apply to the same cup and cap, then they are $\miniE$-congruent.
\end{lemma}

To show \cref{lem:intro_characterization_miniS}, we decompose $\miniE$-congruences as compositions of interchangers and proceed inductively.
Under the hood, we are using the naturality of interchangers in a Gray category; see \cref{subsubsec:scalar_HRSM}.
The situation is similar for other modulo data; see \cref{sec:rewriting_foam} for an example with pivotality.
As a rule of thumb, one should always include naturality axioms in the modulo.
Using the characterization given in \cref{lem:intro_characterization_miniS}, topological arguments can be used to deduce the following:

\begin{lemma}[classification of branchings in $\miniS$]
  \label{lem:intro_classification_miniS}
  Every monomial $\miniS$-branching is branchwise $\miniE$-congruent either to an independent branching, or to a contextualization of one the branchings depicted in \cref{fig:intro_crit_branchings_miniP}.
\end{lemma}

Combining \cref{lem:intro_classification_miniS} with \cref{rewritingfact:contextualization_lemma} (\HLRSMcontextualizationlemma{}), \cref{rewritingfact:HLRSM_congruence_independent_branching} (tamed congruence of independent branchings), \cref{rewritingfact:branchise_confluence_lemma} (\ARSMbranchwiseconfluencelemma{}) and \cref{rewritingfact:tamed_newmann} (\LRSMtamednewmannlemma{}), we conclude that $\miniS^+$ is confluent.

\begin{remark}
  Contrary to classical rewriting, critical branchings are not uniquely defined; hence we refer to branchings in \cref{fig:intro_crit_branchings_miniP} as ``a set of critical branchings'', rather than \emph{the} critical branchings, as we did in \cref{subsubsec:intro_ov_rw}.
\end{remark}

\subsection{Comparison with the literature}
\label{subsubsec:intro_literature}

Our rewriting theory has three main features: it can be \emph{higher}, it can be \emph{linear} and it can be \emph{modulo}.
In other words, it allows to rewrite modulo in 2-sesquicategories and in linear 2-sesquicategories.
We call it \emph{linear higher rewriting modulo}, or \emph{linear Gray rewriting modulo} if interchangers are in the modulo data. We compare with the literature, saying ``strict rewriting'' to emphasize that a rewriting theory does not allow modulo.

Our approach builds on the work of Guiraud--Hoffbeck--Malbos \cite{GHM_ConvergentPresentationsPolygraphic_2019} and their approach to linear (strict) rewriting in associative algebras, and the work of Forest--Mimram \cite{FM_RewritingGrayCategories_2022} and their approach to higher (strict) rewriting in 2-sesquicategories.
In that sense, our work is an extension of \cite{FM_RewritingGrayCategories_2022} to linear rewriting and rewriting modulo; or an extension of \cite{GHM_ConvergentPresentationsPolygraphic_2019} to (Gray) higher rewriting and rewriting modulo.

Alleaume \cite{Alleaume_RewritingHigherDimensional_2018}, Dupont \cite{Dupont_RewritingModuloIsotopies_2022,Dupont_RewritingModuloIsotopies_2021} and Dupont--Ebert--Lauda \cite{DEL_SuperRewritingTheory_2024} respectively developed a (strict) rewriting theory for linear 2-categories, a rewriting modulo theory for linear 2-categories, and a rewriting modulo theory for super-2-categories.
We have important debts to their work: the striking connection between rewriting theory and diagrammatic algebras in the first place, and the importance of rewriting modulo.
However, our approach is distinct in many respects, both theoretically and practically, and even when restricted to their respective setting.
Let us highlight two of them:
\begin{itemize}
  \item Theoretically, our theory is based on tamed congruence. Its importance, even for strict rewriting in linear 2-categories, is highlighted in \cref{subsubsec:intro_independent_contextualized_branchings}; see also \cref{footnote:gap_critical_confluence}.
  \item Practically, we spend a fair amount of time formalizing how to understand branchings modulo (see e.g.\ \cref{subsubsec:intro_classify_branching}).
  The ability to confidently classify branchings is a cornerstone of the theory, since if some forgotten critical branchings are not confluent, it prevents confluence altogether; see also \cref{rem:mistake_article_dupont}.
\end{itemize}
We stress further differences in the text; see \cref{subsec:notes}.

\subsection{Perspectives}
\label{subsec:perspectives}

A rewriting approach can be demanding: finding a practical convergent rewriting system may require a lot of trial and error, and classifying critical branchings can be laborious.
Moreover:
\begin{quote}
  \textsc{Obstruction to rewriting theory:} there is no guarantee that a convergent rewriting system exists, and if so, that it is sufficiently reasonable to be used in practice.
\end{quote}
However, once established, the rewriting perspective provides a rich understanding of the combinatorial structure of the presentation.
We give some future directions of research below.
As more and more examples are studied, we hope that the theory of linear Gray rewriting modulo will evolve into a standard set of tools.

\subsubsection{Computer implementation}

Assuming one has a candidate convergent presentation (see the obstruction above), the rewriting approach is relatively straightforward, at least in principle: enumerate critical branchings and show that each of them is confluent.
We expect both of these processes to be implementable on a computer, at least when working modulo (graded) interchangers.
This should vastly expand what is meant by a ``sufficiently reasonable'' convergent presentation.
This direction of research should relate with current developments in applied category theory, such as the graphical proof-assistants \emph{Globular} and \emph{homotopy.io}
\cite{
  BKV_GlobularOnlineProof_2018,
  Dorn_Associative$n$categories_2023,
  RV_HighlevelMethodsHomotopy_2019,
  BV_DataStructuresQuasistrict_2017,
}.

In some cases, even finding a convergent presentation can be made automatic. For commutative algebras, this is known as Buchberger's algorithm; in generality, this is known as the Knuth--Bendix completion.
One may hope that a similar algorithm can be implemented for certain classes of diagrammatic algebras.

\subsubsection{Classification}

Typically, once a rewriting approach is established for a diagrammatic algebra, it also applies to many variations of this diagrammatic algebra.
Consider the rewriting proof given in \cref{subsec:intro_example} for our example $\miniS$: the sign in the relations $\miniR$ was only relevant once we checked confluence of critical branchings (\cref{fig:intro_crit_branchings_miniP}).
Similarly, the rewriting proof given in \cref{sec:rewriting_foam} for graded $\glt$\nbd-foams tells us how scalars can be chosen in order to get the same basis theorem, leading to precisely two choices, the super-2-category $\gfoam_d$ and $\gfoam_d'$ (\cref{subsec:addendum_deformation_foam}).

Similar ideas have appeared before in the literature \cite{Barbier_DiagramCategoriesBrauer_2024,Elias_DiamondLemmaHecketype_2022}, although not explicitly using rewriting techniques.
We expect rewriting methods to yield classifications of certain types of diagrammatic algebras.

\subsubsection{Higher structures}

The application of rewriting theory goes beyond finding bases.
By considering how rewriting sequences form unoriented cycles, one obtains relations \emph{between} relations, or \emph{coherence data}; see \cref{subsec:coherence_modulo}.
Critical branchings of a convergent rewriting system lead to an explicit description of this coherence data.
For instance, coherent presentations of Artin monoids can be obtained via rewriting theory, leading to a new proof of Deligne's theorem on categorical actions of Artin monoids \cite{GGM_CoherentPresentationsArtin_2015}.
In the linear setting, this coherence data is known as \emph{syzygies}; applications include computation of homological invariants or study of Koszulness \cite{GHM_ConvergentPresentationsPolygraphic_2019}.
In this work, we only discuss coherence in \cref{subsec:coherence_modulo}; however, we expect that the coherence results of \cite{FM_RewritingGrayCategories_2022,GHM_ConvergentPresentationsPolygraphic_2019} can be adapted to our setting. We leave this for future work.

In recent years, stable homotopy theory and $\infty$-categories have become more and more prevalent in link homologies and higher representation theory
\cite{
  LS_KhovanovStableHomotopy_2014,
  HKK_FieldTheoriesStable_2016,
  LLS_KhovanovHomotopyType_2020,
  SSS_OddKhovanovHomotopy_2020,
  DGL+_Spectral2actionsFoams_2025,
  Liu_BraidingComplexOriented_2024,
  LMR+_BraidedMonoidal$infty2$category_2024,
  MR_HigherRepresentationsCornered_2020,
}.
At present however, the complexity of the $\infty$-setting remains a major obstacle to exploration beyond the simplest cases.
Explicit presentations of the coherence data, obtained via rewriting methods, could be an important step forward.

\section{Linear Gray polygraphs}
\label{sec:higher_structure}

This section introduces the necessary categorical structures to present graded-2-categories and define their rewriting theory.
This can be understood as a linear analogue to the work of Forest and Mimram \cite[section 2 and 3]{FM_RewritingGrayCategories_2022}; equivalently, as a generalization of the work of Alleaume on linear $n$\nbd-poly\-graphs \cite{Alleaume_RewritingHigherDimensional_2018} to allow weak interchangers.
For an introduction to the ideas of this section, see \cref{subsubsec:intro_ov_polygraph,subsubsec:intro_ov_gray_rw,subsubsec:intro_linear_gray_polygraph}.

The notion of an $n$\nbd-sesqui\-ca\-te\-gory was first defined by Street \cite{Street_CategoricalStructures_1996} in the case $n=2$. The general case was introduced by Forest--Mimram \cite{FM_RewritingGrayCategories_2022} (following the general theory of Weber \cite{Weber_FreeProductsHigher_2013}) under the name of ``$n$\nbd-precategories'' and independently by Araújo \cite[section~1.6]{Araujo_SimpleStringDiagrams_2022} under the name of ``$n$\nbd-sesqui\-ca\-te\-gories''.
Although we shall follow Forest and Mimram's presentation, we choose Araújo's terminology to avoid confusion with existing notions of $n$\nbd-precategories in the literature.
Enriched category theory provides yet another defining approach to $n$\nbd-sesqui\-ca\-te\-gories; see \cite[section~2.4 and Appendix~A]{FM_RewritingGrayCategories_2022}.

We refer to \cite{FM_RewritingGrayCategories_2022} for the definitions of $n$\nbd-sesquicategories and Gray categories. We also omit a formal definition of linear $n$\nbd-sesquicategories, as it is not strictly necessary for our purpose; a precise definition can be found in \cite{Schelstraete_OddKhovanovHomology_2024}.

\medbreak

\Cref{subsec:gray_polygraphs} discusses the non-linear setting; it is essentially a review of \cite{FM_RewritingGrayCategories_2022}, except for scalar 3-sesquipolygraphs (\cref{subsubsec:scalar_3sesquipolygraphs}). The material of this subsection suffices to read \cref{subsec:coherence_modulo,subsec:abstract_rewriting_modulo,subsec:higher_rewriting_modulo} and its application to foam isotopies in \cref{subsec:rewriting_foam_coherence_iso}.
The generalization to the linear setting is given in \cref{subsubsec:linear_gray_polygraph}.
Every categorical structure is assumed to be small.

\subsection{Gray presentations of 2-categories}
\label{subsec:gray_polygraphs}

\subsubsection{Globular sets}
\label{subsubsec:n-globular_sets}

Let $n\in\bN=\{0,1,\ldots\}$.
An \defemph{$n$\nbd-globular set} $\cC$ is a diagram of sets and functions as follows:
\begin{equation*}
  \begin{tikzcd}
    {\cC_0} & {\cC_1} & \ldots & {\cC_n}
    \arrow["{t_0}", shift left, from=1-2, to=1-1]
    \arrow["{s_0}"', shift right, from=1-2, to=1-1]
    \arrow["{s_1}"', shift right, from=1-3, to=1-2]
    \arrow["{t_1}", shift left, from=1-3, to=1-2]
    \arrow["{s_{n-1}}"', shift right, from=1-4, to=1-3]
    \arrow["{t_{n-1}}", shift left, from=1-4, to=1-3]
  \end{tikzcd}
\end{equation*}
such that $s_j\circ s_{j+1}=s_j\circ t_{j+1}$ and $t_j\circ s_{j+1}=t_j\circ t_{j+1}$ for each $0\leq j < n$. The maps $s_j$ and $t_j$ are respectively called \defemph{source maps} and \defemph{target maps}.
An element $u\in \cC_j$ is called a \defemph{$j$\nbd-cell}, with $s_{j-1}(u)$ and $t_{j-1}(u)$ respectively its \defemph{source} and \defemph{target}, which we sometimes simply denote by $s(u)$ and $t(u)$.
We refer to $j$ as the \defemph{dimension} of $u$.
For $0\leq i <j$, we define the \defemph{$i$\nbd-source} of $u$ to be
\[s_i(u) = s_i\circ \left(\text{\small any suitable composition of source and target maps}\right)(u),\]
where the choice in the bracket does not matter thanks to the properties of source and target maps. The subscript indicates that $s_i(u)$ is an $i$\nbd-cell. We define the \defemph{$i$\nbd-target} $t_i(u)$ similarly.
A \emph{morphism of globular sets $f\colon\cC\to\cD$} is a family $f_i\colon\cC_i\to\cD_i$ of functions that commute with the source and target maps. It is an isomorphism if each function $f_i$ is a bijection.

Given an $n$\nbd-globular set $\cC$, a \defemph{$0$\nbd-sphere} is an ordered pair of $0$\nbd-cells in $\cC$.
For $0<i\leq n$, an \defemph{$i$\nbd-sphere} is an ordered pair $(f,g)$ of $i$\nbd-cells such that $s(f)=s(g)$ and $t(f)=t(g)$. 
For $0\leq i<k\leq n$ and an $i$\nbd-sphere $(f,g)$, we set
\begin{equation*}
\cC_k(f,g) = \{u\in \cC_k\mid s_i(u)=f,t_i(u)=g\}.
\end{equation*}
The source and target maps restrict to maps
$\begin{tikzcd}[cramped]
  \cC_k(f,g) &  \cC_{k+1}(f,g)
  \arrow["{t_k}", shift left, from=1-2, to=1-1]
  \arrow["{s_k}"', shift right, from=1-2, to=1-1]
\end{tikzcd}$.

A \defemph{globular extension of $\cC$} is an $(n+1)$\nbd-globular set $\sP$ such that $\sP_k=\cC_k$ for $0\leq k\leq n$:
\begin{equation*}
    \begin{tikzpicture}
      \node[anchor=west] at (0,0) {
        {\begin{tikzcd}[cramped]
          {\cC_0} & {\cC_1} & \ldots & {\cC_n} & {\sP_{n+1}}
          \arrow["{t_0}", shift left, from=1-2, to=1-1]
          \arrow["{s_0}"', shift right, from=1-2, to=1-1]
          \arrow["{s_1}"', shift right, from=1-3, to=1-2]
          \arrow["{t_1}", shift left, from=1-3, to=1-2]
          \arrow["{s_{n-1}}"', shift right, from=1-4, to=1-3]
          \arrow["{t_{n-1}}", shift left, from=1-4, to=1-3]
          \arrow["{s_{n}}"', shift right, from=1-5, to=1-4]
          \arrow["{t_{n}}", shift left, from=1-5, to=1-4]
        \end{tikzcd}}
      };
      \draw[line width = 1pt,decoration={brace,mirror},decorate] (0,-.5) to node[below=3pt]{$\cC$} (5.5,-.5);
    \end{tikzpicture}
    \vspace*{-.5em}
\end{equation*}
For $k\in\bN$ such that $0\leq k\leq n$, the \defemph{$k$\nbd-restriction of $\cC$} is the subglobular set
\begin{equation*}
  \cC_{\leq k}\coloneqq
  \begin{tikzcd}[cramped]
    {\cC_0} & {\cC_1} & \ldots & {\cC_k}
    \arrow["{t_0}", shift left, from=1-2, to=1-1]
    \arrow["{s_0}"', shift right, from=1-2, to=1-1]
    \arrow["{s_1}"', shift right, from=1-3, to=1-2]
    \arrow["{t_1}", shift left, from=1-3, to=1-2]
    \arrow["{s_{k-1}}"', shift right, from=1-4, to=1-3]
    \arrow["{t_{k-1}}", shift left, from=1-4, to=1-3]
  \end{tikzcd}
\end{equation*}

\subsubsection{2-polygraphs}
\label{subsubsec:2polygraphs}

A \emph{0-polygraph} is a set $\sP_0$. The \emph{free 0-category generated by $\sP_0$} is the set $\sP_0^*=\sP_0$ itself.
A \emph{1-polygraph} is a 1-globular set
\[\begin{tikzcd}[cramped]
  \sP_0 &  \sP_1
  \arrow["{t_0}", shift left, from=1-2, to=1-1]
  \arrow["{s_0}"', shift right, from=1-2, to=1-1]
\end{tikzcd}.\]
The \emph{free 1-category generated by $(\sP_0,\sP_1)$}, denoted $\sP_1^*$, is the category whose objects are elements of $\cC_0$ and morphisms are formal compositions
\[\stdmor\coloneqq \stdmor_d\starop_0 \stdmor_{d-1}\starop_0 \ldots\starop_0 \stdmor_1,\]
where $\stdmor_i\in\sP_1$ and $s_0(\stdmor_{i+1})= t_0(\stdmor_i)$ for all $1\leq i\leq d-1$.
It inherits a source $s_0(\stdmor)=s_0(\stdmor_1)$ and target $t_0(\stdmor)=t_0(\stdmor_d)$.
Composition $\starop_0$ is given by concatenation.

\begin{example}
  \label{ex:polygraphA}
  Recall $\miniP$ from \cref{subsec:intro_example}.
  Its underlying 1-polygraph is $\miniP_0=\{*\}$ and $\miniP_1=\left\{\;\tikzpic{\ministrand}[yscale=.4]\;\right\}$, with $s_0(\;\tikzpic{\ministrand}[yscale=.4]\;)=t_0(\;\tikzpic{\ministrand}[yscale=.4]\;)=*$. The free 1-category is
\[(\miniP_1)^*=\{\;\tikzpic{\ministrand}[yscale=.4]\;\ldots\;\tikzpic{\ministrand}[yscale=.4]\mid n\in\bN\;\},\]
where $n$ denotes the number of strands.
\end{example}

A \emph{2-polygraph} is a globular extension $\sP_2$ of $\sP_1^*$:
\begin{gather*}
  \begin{tikzcd}[ampersand replacement=\&,cramped]
    {\sP_0^*} \& {\sP_1^*} \\
    {\sP_0} \& {\sP_1} \& {\sP_2}
    \arrow[from=2-1, to=1-1,equals]
    \arrow[shift right, from=2-2, to=1-1,"s_0"']
    \arrow[shift left, from=2-2, to=1-1,"t_0"{xshift=-1pt}]
    \arrow[from=2-2, to=1-2]
    \arrow[shift right, from=2-3, to=1-2,"s_1"']
    \arrow[shift left, from=2-3, to=1-2,"t_1"{xshift=-1pt}]
  \end{tikzcd}
\end{gather*}
Given a 2-polygraph, a \emph{whiskered 2-generator} is a formal composition of the form $\stdmor\starop_0 \stddiag\starop_0 \stdmorbis$, where $\stdmor,\stdmorbis\in\sP_1^*$ and $\stddiag\in\sP_2$, and such that $s_0(\stdmor)=t_0(\stddiag)$ and $s_0(\stddiag)=t_0(\stdmor)$.
We depict it as:
\begin{gather*}
  \stdmor\starop_{0}\stddiag\starop_{0}\stdmorbis
  \;=\;
  \newcommand{\tempid}[1]{\draw (#1,0) to (#1,3);}
  \tikzpic{
    \tempid{0}\tempid{1}\tempid{2}\tempid{3}\tempid{4}\tempid{5}
    \diagdots{.5}{0}\diagdots{2.5}{0}\diagdots{4.5}{0}
    \diagdots{.5}{3}\diagdots{2.5}{3}\diagdots{4.5}{3}
    \node at (.5,1.5) {$\stdmor$};
    \node at (4.5,1.5) {$\stdmorbis$};
    \node[rounded corners,fill=white,draw=black,minimum width=1cm] at (2.5,1.5) {$\stddiag$};
  }[scale=.5,yscale=.7]\;.
\end{gather*}
It inherits a 1-source $s_1(\stdmor\starop_{0}\stddiag\starop_{0}\stdmorbis)=\stdmor\starop_{0}s_1(\stddiag)\starop_{0}\stdmorbis$ and a 1-target $t_1(\stdmor\starop_{0}\stddiag\starop_{0}\stdmorbis)=\stdmor\starop_{0}t_1(\stddiag)\starop_{0}\stdmorbis$.
We define a globular extension $\sP_2^*$ of $\sP_1^*$ as the set of formal compositions
\[D\coloneqq D_d\starop_1 D_{d-1}\starop_1\ldots\starop_1 D_1,\]
where $D_i$ is a whiskered 2-generator and $s_1(\stddiag_{i+1})= t_1(\stddiag_i)$ for all $1\leq i\leq d-1$.
The compositions $\starop_0$ and $\starop_1$ endow $\sP_2^*$ with the structure of a 2-sesquicategory, called the \emph{free 2-sesquicategory generated by $(\sP_0,\sP_1,\sP_2)$}.

\begin{example}
  \label{ex:polygraphB}
  Following up on \cref{ex:polygraphA}:
  \begin{gather*}
    \miniP_2=\left\{\;\tikzpic{\minicap}[scale=.7]\,,\,\tikzpic{\minicup}[scale=.7]\;\right\},
  \end{gather*}
  where the source and target are read as the bottom and top boundary, respectively.
  The free 2-sesquicategory $\miniP_2^*$ consists of diagrams with cups and caps, up to rectilinear isotopies that preserve the relative vertical positions of cups and caps.
\end{example}

\subsubsection{Contexts}
\label{subsubsec:contexts}

Let $(\sP_0,\sP_1,\sP_2)$ be a 2-polygraph.
Let $\square=(\stdmor_s,\stdmor_t)$ be a 1-sphere, that is, $\stdmor_s,\stdmor_t\in\sP_1^*$ are such that $s_0(\stdmor_s)=s_0(\stdmor_t)$ and $t_0(\stdmor_s)=t_0(\stdmor_t)$. A \emph{context in $\sP$ with boundary $\square$} is a formal composition
\begin{gather*}
  \Gamma
  \;\coloneqq\;
  D_1\starop_1(\stdmorbis_1\starop_0\square\starop_0 \stdmorbis_0)\starop_1D_0
  \;=\;
  \newcommand{\tempid}[1]{\draw (#1,-2) to (#1,5);}
  \tikzpic{
    \tempid{0}\tempid{1}\tempid{2}\tempid{3}\tempid{4}\tempid{5}
    \diagdots{.5}{.5}\diagdots{2.5}{.5}\diagdots{4.5}{.5}
    \diagdots{.5}{-1.5}\diagdots{2.5}{-1.5}\diagdots{4.5}{-1.5}
    \diagdots{.5}{2.5}\diagdots{2.5}{2.5}\diagdots{4.5}{2.5}
    \diagdots{.5}{4.5}\diagdots{2.5}{4.5}\diagdots{4.5}{4.5}
    \node at (.5,1.5) {$\stdmorbis_1$};
    \node at (4.5,1.5) {$\stdmorbis_0$};
    \node[rounded corners,fill=white,draw=black,minimum width=1cm,minimum height=.5cm] at (2.5,1.5) {$\square$};
    \node[rounded corners,fill=white,draw=black,minimum width=3cm] at (2.5,-.5) {$D_0$};
    \node[rounded corners,fill=white,draw=black,minimum width=3cm] at (2.5,3.5) {$D_1$};
  }[scale=.5,yscale=.7]\;,
\end{gather*}
where $\stdmorbis_0,\stdmorbis_1\in\sP_1^*$ and $D_0,D_1\in\sP_2^*$ are suitably composable.
Given a 2-cell $A\in\sP_2^*$ whose 1-boundary is $\square$, we write $\Gamma[A]$, the \emph{contextualization of $A$}, the 2-cell
\[\Gamma[A]=D_1\starop_1(\stdmorbis_1\starop_0 A\starop_0 \stdmorbis_0)\starop_1D_0.\]

\subsubsection{3-sesquipolygraphs}
\label{subsubsec:3sesquipolygraphs}

Let $(\sP_0,\sP_1,\sP_2)$ be a 2-polygraph.
A \emph{3-sesquipolygraph} is a globular extension $\sP_3$ of $\sP_2^*$:
\begin{gather*}
  \begin{tikzcd}[ampersand replacement=\&,cramped]
    {\sP_0^*} \& {\sP_1^*} \& {\sP_2^*}\\
    {\sP_0} \& {\sP_1} \& {\sP_2} \& {\sP_3}
    \arrow[shift right, from=2-4, to=1-3,"s_2"']
    \arrow[shift left, from=2-4, to=1-3,"t_2"{xshift=-1pt}]
    \arrow[from=2-1, to=1-1,equals]
    \arrow[shift right, from=2-2, to=1-1,"s_0"']
    \arrow[shift left, from=2-2, to=1-1,"t_0"{xshift=-1pt}]
    \arrow[from=2-2, to=1-2]
    \arrow[shift right, from=2-3, to=1-2,"s_1"']
    \arrow[shift left, from=2-3, to=1-2,"t_1"{xshift=-1pt}]
    \arrow[from=2-3, to=1-3]
  \end{tikzcd}
\end{gather*}

Let $\sP=(\sP_0,\sP_1,\sP_2,\sP_3)$ be a 3-sesquipolygraph.
Let $f\colon s(f)\to t(f)$ be a 3-cell in $\sP_3$.
Given a context $\Gamma$ as in \cref{subsubsec:contexts}, the \emph{contextualization of $f$} is the formal 3-cell $\Gamma[f]\colon\Gamma[s(f)]\to\Gamma[t(f)]$, depicted as
\begin{gather*}
  \Gamma[f]
  \;=\;
  \newcommand{\tempid}[1]{\draw (#1,-2) to (#1,5);}
  \tikzpic{
    \tempid{0}\tempid{1}\tempid{2}\tempid{3}\tempid{4}\tempid{5}
    \diagdots{.5}{.5}\diagdots{2.5}{.5}\diagdots{4.5}{.5}
    \diagdots{.5}{-1.5}\diagdots{2.5}{-1.5}\diagdots{4.5}{-1.5}
    \diagdots{.5}{2.5}\diagdots{2.5}{2.5}\diagdots{4.5}{2.5}
    \diagdots{.5}{4.5}\diagdots{2.5}{4.5}\diagdots{4.5}{4.5}
    \node at (.5,1.5) {$\stdmorbis_1$};
    \node at (4.5,1.5) {$\stdmorbis_0$};
    \node[rounded corners,fill=white,draw=black,minimum width=1cm,minimum height=.5cm] at (2.5,1.5) {$s(f)$};
    \node[rounded corners,fill=white,draw=black,minimum width=3cm] at (2.5,-.5) {$D_0$};
    \node[rounded corners,fill=white,draw=black,minimum width=3cm] at (2.5,3.5) {$D_1$};
  }[scale=.5,yscale=.7]
  \quad\to\quad
  \tikzpic{
    \tempid{0}\tempid{1}\tempid{2}\tempid{3}\tempid{4}\tempid{5}
    \diagdots{.5}{.5}\diagdots{2.5}{.5}\diagdots{4.5}{.5}
    \diagdots{.5}{-1.5}\diagdots{2.5}{-1.5}\diagdots{4.5}{-1.5}
    \diagdots{.5}{2.5}\diagdots{2.5}{2.5}\diagdots{4.5}{2.5}
    \diagdots{.5}{4.5}\diagdots{2.5}{4.5}\diagdots{4.5}{4.5}
    \node at (.5,1.5) {$\stdmorbis_1$};
    \node at (4.5,1.5) {$\stdmorbis_0$};
    \node[rounded corners,fill=white,draw=black,minimum width=1cm,minimum height=.5cm] at (2.5,1.5) {$t(f)$};
    \node[rounded corners,fill=white,draw=black,minimum width=3cm] at (2.5,-.5) {$D_0$};
    \node[rounded corners,fill=white,draw=black,minimum width=3cm] at (2.5,3.5) {$D_1$};
  }[scale=.5,yscale=.7]\;.
\end{gather*}
We denote $\Cont(\sP)$ the set of contextualized 3-cells $\Gamma[f]$, for all contexts $\Gamma$ and 3-cells $f\in\sP_3$.

\subsubsection{Scalar 3-sesquipolygraphs}
\label{subsubsec:scalar_3sesquipolygraphs}

Let $\ring$ be a commutative ring.

A \emph{scalar 3-sesquipolygraph} is a 3-sesquipolygraph $\sP=(\sP_0,\sP_1,\sP_2,\sP_3)$ together with a function $\scl\colon\sP_3\to\ring$.
If the image of $\scl$ consists of invertible scalars $\ring^\times\subset\ring$, we say that $\sP$ is \emph{scalar-invertible}.
We extend $\scl$ to $\Cont(\sP)$ by setting $\scl(\Gamma[f])=\scl(f)$.

\begin{example}
  \newcommand{\tpscl}{.5}
  \label{ex:polygraphC}
  Recall \cref{ex:polygraphB}. We define a scalar 3-sesquipolygraph on $\miniP_2$ as follows:
  \begin{gather*}
  \scl(\miniP)_3=
  \left\{\;
  \begin{array}{ccccc}
  \tikzpic{
    \ministrand[0][1]\minicap[1][1]
    \minicup\ministrand[2][0]
  }[scale=\tpscl]
  \;\Rrightarrow\;
  \tikzpic{
    \ministrand[0][1]
    \ministrand
  }[scale=\tpscl]
  \;&,&\;
  \tikzpic{
    \minicap[0][1]\ministrand[2][1]
    \ministrand\minicup[1][0]
  }[scale=\tpscl]
  \;\overset{-1}{\Rrightarrow}\;
  \tikzpic{
    \ministrand[0][1]
    \ministrand
  }[scale=\tpscl]
  \;&,&\;
  \tikzpic{
    \draw (0,2) to
      node[black_dot,pos=.3] {}
      node[right,pos=.3] {\scriptsize $a$}
      (0,0);
    \draw (1,2) to
      node[black_dot,pos=.7] {}
      node[right,pos=.7] {\scriptsize $b$}
      (1,0);
  }[scale=0.6]
  \overset{(-1)^{p(a)p(b)}}{\Rrightarrow}
  \tikzpic{
    \draw (0,2) to
      node[black_dot,pos=.7] {}
      node[right,pos=.7] {\scriptsize $a$}
      (0,0);
    \draw (1,2)  to
      node[black_dot,pos=.3] {}
      node[right,pos=.3] {\scriptsize $b$}
      (1,0);
  }[scale=0.6]
  \\
  &&&&
  \text{\scriptsize for all $a,b\in(\miniP_2)^*$}
  \end{array}
  \;\right\}.
\end{gather*}
Here we depict the scalar associated to each 3-cell in $\miniP_3$ above the arrow.
\end{example}

\subsubsection{Gray polygraphs}
\label{subsubsec:gray_polygraph}

Let $\sQ$ be a 2\nbd-sesqui\-poly\-graph. The \defemph{3\nbd-sesqui\-poly\-graph of interchangers} is the 3\nbd-sesqui\-poly\-graph $\sQ\Gray$ such that $\sQ\Gray_{\leq 2}=\sQ$ and $\sQ\Gray_3$ consists of \defemph{interchange generators}, defined for each 0-compo\-sable $A,\chi,B$ with $A\colon \xi\Rightarrow \xi'$, $\chi\in \sP_1^*$, $B\colon \zeta\Rightarrow \zeta'\in\sP_2$, as the 3-cell
\[
  X_{A,\chi,B}\colon (A\starop_0 \chi\starop_0 \zeta')\starop_1 (\xi\starop_0 \chi\starop_0 B) \Rrightarrow (\xi'\starop_0 \chi\starop_0 B)\starop_1(A\starop_0\chi\starop_0 \zeta),
\]
pictured as:
\[
  X_{A,\chi,B}\colon
  \newcommand{\tempid}[1]{\draw (#1,0) to (#1,4);}
  \tikzpic{
    \tempid{0}\tempid{1}\tempid{2}\tempid{3}\tempid{4}\tempid{5}
    \diagdots{.5}{0}\diagdots{2.5}{0}\diagdots{4.5}{0}
    \diagdots{.5}{4}\diagdots{2.5}{4}\diagdots{4.5}{4}
    \node[below] at (.5,0) {\scriptsize $\xi$};
    \node[below] at (2.5,0) {\scriptsize $\chi$};
    \node[below] at (4.5,0) {\scriptsize $\zeta$};
    \node[above] at (.5,4) {\scriptsize $\xi'$};
    \node[above] at (2.5,4) {\scriptsize $\chi$};
    \node[above] at (4.5,4) {\scriptsize $\zeta'$};
    \node[rounded corners,fill=white,draw=black,minimum width=1cm,minimum height=.3cm] at (.5,3) {\scriptsize $A$};
    \node[rounded corners,fill=white,draw=black,minimum width=1cm,minimum height=.3cm] at (4.5,1) {\scriptsize $B$};
  }[scale=.5,yscale=.7]
  \Rrightarrow\mspace{10mu}
  \tikzpic{
    \tempid{0}\tempid{1}\tempid{2}\tempid{3}\tempid{4}\tempid{5}
    \diagdots{.5}{0}\diagdots{2.5}{0}\diagdots{4.5}{0}
    \diagdots{.5}{4}\diagdots{2.5}{4}\diagdots{4.5}{4}
    \node[below] at (.5,0) {\scriptsize $\xi$};
    \node[below] at (2.5,0) {\scriptsize $\chi$};
    \node[below] at (4.5,0) {\scriptsize $\zeta$};
    \node[above] at (.5,4) {\scriptsize $\xi'$};
    \node[above] at (2.5,4) {\scriptsize $\chi$};
    \node[above] at (4.5,4) {\scriptsize $\zeta'$};
    \node[rounded corners,fill=white,draw=black,minimum width=1cm,minimum height=.3cm] at (.5,1) {\scriptsize $A$};
    \node[rounded corners,fill=white,draw=black,minimum width=1cm,minimum height=.3cm] at (4.5,3) {\scriptsize $B$};
  }[scale=.5,yscale=.7]
\]

\begin{definition}
  A \defemph{Gray polygraph} is a 3\nbd-sesqui\-poly\-graph $\sP$ such that
  \[{\sP_{\leq 2}\Gray\subset \sP}.\]
  In other words, $\sP$ contains its own 3\nbd-sesqui\-poly\-graph of interchangers.
\end{definition}

\begin{example}
  \label{ex:polygraphD}
  The scalar 3-sesquipolygraph from \cref{ex:polygraphC} is a Gray polygraph.
\end{example}

Let $\sP$ be a 3-sesquipolygraph. Let $\sim_\sP$ be the smallest equivalence relation on $\sP_2^*$ such that $\Gamma[s(f)]\sim_\sP\Gamma[t(f)]$ for each $\Gamma[f]\in\Cont(\sP)$.
We write $[\sP]_\sim$ the 2-sesquicategory obtained by quotienting $\sP_2^*$ by $\sim_\sP$.
If $\sP$ is a Gray polygraph, then $[\sP]_\sim$ is a 2\nbd-cate\-gory.
This leads to the following definition:

\begin{definition}
  A \emph{presentation of a 2\nbd-cate\-gory $\cC$} is the data of a Gray polygraph $\sP$ such that $[\sP]_\sim$ is isomorphic to $\cC$.
\end{definition}

\subsection{Linear Gray presentations of graded-2-categories}
\label{subsec:linear_gray_polygraphs}

\begin{notation}
  \label{not:sesquicat_notation}
  We fix throughout the section an abelian group $G$, a commutative ring $\ring$ and a $\bZ$\nbd-bilinear map $\mu\colon G\times G\to\ring^\times$. The word ``graded'' always refers to $G$\nbd-graded, and ``linear'' to $\ring$\nbd-linear. Given a homogeneous element $v$, we write $\deg(v)$ its grading.
\end{notation}

\subsubsection{Linear 3-sesquipolygraphs}

A \emph{graded 2-polygraph} is a 2-polygraph $(\sP_0,\sP_1,\sP_2)$ such that $\sP_2$ is a $G$-graded set, that is, endowed with a function $\deg\colon\sP_2\to G$.
The degree function $\deg$ extends to the free 2-sesquicategory $\sP_2^*$ as $\deg(\stdmor\starop_{0}\stddiag\starop_{0}\stdmorbis)=\deg(\stddiag)$ and
\[
\deg(D_d\starop_1 D_{d-1}\starop_1\ldots\starop_1 D_1)
=
\deg(D_d)+\deg(D_{d-1})+\ldots +\deg(D_1).
\]

Let $(\sP_0,\sP_1,\sP_2)$ be a graded 2-polygraph.
We write $\sP_2^l$ the globular extension of $\sP_1^*$ such that for each 1-sphere $(\stdmor_s,\stdmor_t)$, we set $\sP_2^l(\stdmor_s,\stdmor_t)$ to be the free $G$-graded $\ring$-module generated by $\sP_2^*(\stdmor_s,\stdmor_t)$.
A \emph{(graded) linear 3-sesquipolygraph} is a globular extension $\sP_3$ of $\sP_2^l$:
\begin{gather*}
  \begin{tikzcd}[ampersand replacement=\&,cramped]
    {\sP_0^*} \& {\sP_1^*} \& {\sP_2^l}\\
    {\sP_0} \& {\sP_1} \& {\sP_2} \& {\sP_3}
    \arrow[shift right, from=2-4, to=1-3,"s_2"']
    \arrow[shift left, from=2-4, to=1-3,"t_2"{xshift=-1pt}]
    \arrow[from=2-1, to=1-1,equals]
    \arrow[shift right, from=2-2, to=1-1,"s_0"']
    \arrow[shift left, from=2-2, to=1-1,"t_0"{xshift=-1pt}]
    \arrow[from=2-2, to=1-2]
    \arrow[shift right, from=2-3, to=1-2,"s_1"']
    \arrow[shift left, from=2-3, to=1-2,"t_1"{xshift=-1pt}]
    \arrow[from=2-3, to=1-3]
  \end{tikzcd}
\end{gather*}
such that for each $f\in\sP_3$, we have $\deg(s(f))=\deg(t(f))$.

\subsubsection{Monomial linear 3-sesquipolygraphs}
\label{subsubsec:monomial_linear_nsesquipoly}

A linear $3$\nbd-sesqui\-poly\-graph $\sP$ is \defemph{monomial} if:
\[
  s_2(r)\in\sP_{2}^*
  \;\an\;
  t_{2}(r)\in\ring\sP_{2}^*\text{ for all }r\in\sP_{3},
\]
where $\ring\sP_{2}^*$ is the subset of $\sP_{2}^l$ consisting of vectors of the form $\lambda b$ for $\lambda\in\ring$ and $b\in\sP_{2}^*$. We further say it is \emph{monomial-invertible} if $t_{2}(r)\in\ring^\times\sP_{2}^*$.
We have the following canonical bijection (recall \cref{subsubsec:scalar_3sesquipolygraphs}):
\[\begin{tikzcd}[ampersand replacement=\&,cramped]
  {\left\{\;
  \substack{
    \text{\normalsize scalar(-invertible)}\\[1ex]
    \text{\normalsize  $3$\nbd-sesqui\-poly\-graph}
  }\;\right\}}
  \&
  {\left\{\;
  \substack{
    \text{\normalsize monomial(-invertible)}\\[1ex]
    \text{\normalsize linear $3$\nbd-sesqui\-poly\-graph}
  }\;\right\}}
  \arrow[shift left=1, from=1-2, to=1-1,"\scl"]
  \arrow[shift left=1, from=1-1, to=1-2,"\lin"]
\end{tikzcd}\]
Given a scalar $3$\nbd-sesquipoly\-graph $\sP$, its \defemph{linearization} is the monomial linear $3$\nbd-sesqui\-poly\-graph $\mathrm{lin}(\sP)$ defined as $\mathrm{lin}(\sP)_{\leq 2}=\sP_{\leq 2}$ and $\mathrm{lin}(\sP)_{3}=\sP_{3}$ as sets, but with the following source and target maps:
\begin{gather*}
  \mathrm{lin}(s_{2})(r) = s_{2}(r)\quad\an\quad \mathrm{lin}(t_{2})(r) = \mathrm{scl}(r)t_{2}(r).
\end{gather*}
The inverse $\scl\coloneqq\lin^{-1}$ is defined analogously.

\begin{example}
  In \cref{subsec:intro_example}, we defined the scalar 3-sesquipolygraph $\scl(\miniP)$ of \cref{ex:polygraphC} as a monomial linear 3-sesquipolygraph $\miniP$.
\end{example}

\subsubsection{Linear Gray polygraphs}

\label{subsubsec:linear_gray_polygraph}

Let $\sQ$ be a graded 2\nbd-sesqui\-poly\-graph.
Recall the 3-sesquipolygraph of interchangers $\sQ\Gray$ from \cref{subsubsec:gray_polygraph}.
We equip $\sQ\Gray$ with the function
\begin{align*}
  \sQ\Gray &\to \ring^\times,\\
  X_{\alpha,g,\beta}&\mapsto \mu(\deg\alpha,\deg\beta).
\end{align*}
This makes $\sQ\Gray$ a scalar-invertible 3\nbd-sesqui\-poly\-graph, called the \defemph{3\nbd-sesqui\-poly\-graph of $(G,\mu)$-graded interchangers}.
Given the previous section, we may also view $\sQ\Gray$ as a monomial-invertible linear $3$\nbd-sesquipo\-ly\-graph.

\begin{definition}
  \label{defn:linear_Gray_polygraph}
  A \defemph{$(G,\mu)$\nbd-linear Gray polygraph} is a graded linear 3\nbd-sesqui\-polygraph $\sP$ such that $\sP$ contains its own monomial graded 3\nbd-sesqui\-polygraph of $(G,\mu)$\nbd-graded interchangers.
\end{definition}

A relation $\sim$ defined on a $\ring$-module is \emph{linear} if $v\sim w$ implies $\lambda v + u\sim\lambda w + u$ for every scalar $\lambda\in\ring$ and vector $u$.
For $\sP$ a (graded) linear 3-sesquipolygraph, denote $\sim_\sP$ the smallest linear equivalence relation such that $\Gamma[s(f)]\sim_\sP\Gamma[t(f)]$ for each $\Gamma[f]\in\Cont(\sP)$. We write $[\sP]_\sim$ the (graded) linear 2-sesquicategory obtained by quotienting $\sP_2^l$ by $\sim$; more precisely, by quotienting $\sP_2^l(\xi_s,\xi_t)$ for every 1-sphere 
$(\xi_s,\xi_t)$.
If $\sP$ is a $(G,\mu)$\nbd-linear Gray polygraph, then $[\sP]_\sim$ is a $(G,\mu)$-graded-2\nbd-cate\-gory \cite{SV_OddKhovanovHomology_2023}.
That is, a structure akin to a linear 2-category, but where the interchange law is twisted:
\begin{equation}
  \label{eq:graded_interchange_law}
  \tikzpic{
    \draw (0,2) to
      node[black_dot,pos=.3] {}
      node[right,pos=.3] {\scriptsize $a$}
      (0,0);
    \draw (1,2) to
      node[black_dot,pos=.7] {}
      node[right,pos=.7] {\scriptsize $b$}
      (1,0);
  }[scale=0.6]
  \;=\;\mu(\deg a,\deg b)
  \tikzpic{
    \draw (0,2) to
      node[black_dot,pos=.7] {}
      node[right,pos=.7] {\scriptsize $a$}
      (0,0);
    \draw (1,2)  to
      node[black_dot,pos=.3] {}
      node[right,pos=.3] {\scriptsize $b$}
      (1,0);
  }[scale=0.6]
\end{equation}
This leads to the following definition:

\begin{definition}
  \label{defn:presentation_graded_two_cat}
  A \emph{presentation of a $(G,\mu)$\nbd-graded-2\nbd-cate\-gory $\cC$} is the data of a $(G,\mu)$\nbd-linear Gray polygraph $\sP$ such that $[\sP]_\sim$ is isomorphic to $\cC$.
\end{definition}

\begin{remark}
  If $G=\{*\}$ is trivial, a $(\{*\},\id)$\nbd-graded-2\nbd-cate\-gory is just a linear 2\nbd-cate\-gory. Hence, a $(\{*\},\id)$\nbd-linear Gray polygraph, and in particular a $(\{*\},\id)$\nbd-scalar Gray polygraph, defines a presentation of a linear 2\nbd-cate\-gory.
\end{remark}

\begin{example}
  \label{ex:polygraphE}
  Consider the scalar Gray polygraph $\scl(\miniP)$ from \cref{ex:polygraphC}.
  Its linearization $\miniP$ (see \cref{subsubsec:monomial_linear_nsesquipoly}) presents the graded-2-category $\miniC$ discussed in \cref{subsec:intro_example}.
\end{example}

\section{Linear Gray rewriting modulo}
\label{sec:foundation_rewriting}

This section introduces \emph{linear Gray rewriting modulo}, the rewriting theory of graded-2-categories.
It has three main features: it allows \emph{modulo}, it is \emph{linear}, and it is \emph{higher}.
Each of these features requires its own treatment, in addition to a description of how they combine.

\Cref{subsec:coherence_modulo,subsec:abstract_rewriting_modulo,subsec:higher_rewriting_modulo} introduce coherence modulo, abstract rewriting modulo and higher rewriting modulo, respectively.
After reading these sections, the reader is equipped to read the application to foam isotopies given in \cref{subsec:rewriting_foam_coherence_iso}.
\Cref{subsec:linear_rewriting_modulo,subsec:higher_linear_rewriting_modulo} discuss linear rewriting modulo and higher linear rewriting modulo, respectively; these sections are the technical core of the paper.
We sometimes use graded $\glt$-foams (\cref{sec:rewriting_foam}) as a source of (counter)examples; the reader is then referred to \cref{subsec:review_foams} for the definitions.

\begin{notation}
  We fix the same notations $G$, $\Bbbk$ and $\mu$ as in \cref{not:sesquicat_notation} throughout the section.
  We use blackboard fonts (e.g.\ $\aP$) for abstract notions, typewriter fonts (e.g.\ $\lP$) for linear notions and sans serif fonts (e.g.\ $\sP$) for (linear) higher notions.
\end{notation}

\subsection{Coherence modulo from convergence modulo}
\label{subsec:coherence_modulo}

Given a groupoid $\cP$, one may be interested in understanding its set of connected components $\pi_0(\cP)$,
or its $\pi_1(\cP)$, that is, its \emph{coherence}.
In this section, we describe how to use convergence to answer both questions.
We fix the following data:
\bigbreak
\hfill%
\begin{tabular}{c}
  a groupoid $\cP$,
  a category $\cS$ such that $\cS^\top=\cP$,\\
  and a wide\footnote{Containing all objects; equivalently, containing all identities.} subgroupoid $\cE\subset\cS$.
\end{tabular}%
\hfill \customlabel{data_convergence_modulo}{$(*)$}
\bigbreak
\noindent
Here $\cS^\top$ is the localisation of the category $\cS$.
Morphisms in $\cS$ are depicted with plain arrow $x\to y$ (or $\begin{tikzcd}[cramped]
  x\rar & y
\end{tikzcd}$).
Morphisms in $\cE$ are depicted with unoriented wiggly lines $x\sim y$ (or $\begin{tikzcd}[cramped]
  x\rar[snakecd] & y
\end{tikzcd}$).
Unspecified morphisms $x\to_\cS y$ (resp.\ $x\sim_\cE y$) indicate the mere existence of a morphism in $\cS$ (resp.\ in $\cE$) between $x$ and $y$.
Note the use of subscripts to specify the category, used extensively in what follows.

\begin{definition}
  \label{defn:general_coherence_modulo}
  The groupoid $\cP$ is said to be \emph{coherent modulo $\cE$} if every endomorphism is parallel to a morphism in $\cE$, up to conjugation.
  That is, for every endomorphism $f\colon x\to_\cP x$, there exist morphisms $g\colon x\to_\cP y$ and $e\colon y\sim_\cE y$:
  \begin{gather*}
    \forall\;
    \tikzpic{
      \node (A) at (.5,0) {$x$};
      \draw[->] (A) to[out=-120,in=0] (0,-.5) to[out=180,in=-90] (-.5,0) node[left]{$f$} to[out=90,in=180] (0,.5) to[out=0,in=120] (A);
    }[scale=.8],
    \qquad
    \exists\;
    \tikzpic{
      \node (B) at (-.5,0) {$y$};
      \node (A) at (-2,0) {$x$};
      \draw[snakecd] (B) to[out=-60,in=180] (0,-.5) to[out=0,in=-90] (.5,0) node[right]{$e$} to[out=90,in=0] (0,.5) to[out=180,in=60] (B);
      \draw[->] (-2+.2,.1) to node[above]{$g$} (-.5-.2,.1);
      \draw[<-] (-2+.2,-.1) to node[below]{$g^{-1}$} (-.5-.2,-.1);
    }\;.
  \end{gather*}
\end{definition}

\begin{lemma}[transitivity of coherence modulo]
  \label{lem:transitivity_coherence_modulo}
  Let $\cC$ be a groupoid and $\cE\subset \cD\subset \cC$ wide subgroupoids of $\cC$. If $\cC$ is coherent modulo $\cD$ and $\cD$ is coherent modulo $\cE$, then $\cC$ is coherent modulo $\cE$.
  \hfill\qed
\end{lemma}

An \emph{$\cS$\nbd-branching} is a pair $(f,f')$ with $f\colon x\to_\cS y$ and $f'\colon x\to_\cS y'$.
An \emph{$\cS$\nbd-confluence} is a pair $(g,g')$ with $g\colon y\to_\cS z$ and $g'\colon y'\to_\cS z$.
A branching $(f,f')$ is said to be \emph{confluent} if there exists a confluence $(g,g')$ such that $g\circ f\equiv g'\circ f'$:
\begin{IEEEeqnarray*}{CcCcC}
  \begin{tikzcd}[ampersand replacement=\&,cramped,row sep=.7em]
    \& y \\
    x \& \\
    \& {y'}
    \arrow["f", curve={height=-9pt}, from=2-1, to=1-2]
    \arrow["{f'}"', curve={height=9pt}, from=2-1, to=3-2]
  \end{tikzcd}
  &&
  \begin{tikzcd}[ampersand replacement=\&,cramped,row sep=.7em]
    \& y \\
    \&\& z \\
    \& {y'}
    \arrow["g", curve={height=-9pt}, from=1-2, to=2-3]
    \arrow["{g'}"'{shift={(-.05,.05)}}, curve={height=9pt}, from=3-2, to=2-3]
  \end{tikzcd}
  &&
  \begin{tikzcd}[ampersand replacement=\&,cramped,row sep=.7em]
    \& y \\
    x \&\& z \\
    \& {y'}
    \arrow["g", curve={height=-9pt}, from=1-2, to=2-3]
    \arrow[draw=none, from=1-2, to=3-2]
    \arrow["f", curve={height=-9pt}, from=2-1, to=1-2]
    \arrow["{f'}"', curve={height=9pt}, from=2-1, to=3-2]
    \arrow["{g'}"'{shift={(-.05,.05)}}, curve={height=9pt}, from=3-2, to=2-3]
  \end{tikzcd}
  \\*
  \text{\scriptsize an $\cS$\nbd-branching}
  &\mspace{70mu}&
  \text{\scriptsize an $\cS$\nbd-confluence}
  &\mspace{70mu}&
  \substack{\text{an $\cS$\nbd-branching}\\\text{which is confluent}}
\end{IEEEeqnarray*}
We say that $\cS$ is \emph{confluent} if every branching is confluent.

We say that $\cS$ is \emph{terminating modulo $\cE$} if every infinite sequence $(f_n)_{n\in\bN}$ of morphisms in $\cS$ with ${t(f_n)=s(f_{n+1})}$ eventually terminates in morphisms in $\cE$.
We say that $\cS$ is \emph{convergent modulo $\cE$} if it is both confluent and terminating.
An object $y\in\cS$ for which ${y\to_\cS z}$ implies ${y\sim_\cE z}$ is called an \emph{$\cS$\nbd-normal form modulo $\cE$}. We denote $\NF_\cS$ the set of $\cS$\nbd-normal forms modulo $\cE$. If $x\to_\cS y$ and $y$ is an $\cS$\nbd-normal form modulo $\cE$, we say that \emph{$y$ is an $\cS$\nbd-normal form modulo $\cE$ for $x$}.

When $\cE$ is clear from the context, we sometimes omit ``modulo $\cE$'' and simply say ``termination'', ``convergence'' and ``normal form''.

\medbreak

The following two propositions explain how convergence modulo leads to a description of $\pi_0(\cS^{\top})$ and $\pi_1(\cS^{\top})$, respectively.

\begin{proposition}
  \label{prop:pi_zero_from_normal_forms}
  If $\cS$ is convergent modulo $\cE$, then the canonical mapping $\NF_\cS\to \pi_0(\cS^\top)$ sending an $\cS$\nbd-normal form to its connected component induces a bijection
  \[\begin{tikzcd}[cramped]
    \faktor{\NF_\cS}{\pi_0(\cE)}\rar["\sim"] &\pi_0(\cS^\top)
  \end{tikzcd},\]
  where the quotient identifies objects belonging to the same connected component of $\cE$.
\end{proposition}

\begin{proof}
  Since $\cS$ is terminating modulo $\cE$, the mapping ${\NF_\cS\to \pi_0(\cS^\top)}$ is surjective. Assume then that $x$ and $y$ are $\cS$\nbd-normal forms belonging to the same connected component of $\cS^\top$.
  We have the following (see e.g.\ \cite[1.3.18]{ABG+_PolygraphsRewritingHigher_2025}):
  \begin{quote}
    \emph{Church--Rosser property:} If $\cS$ is confluent, then every morphism $f\in\cS^\top$ can be decomposed as $f= h^{-1}\circ g$ for some $f\colon x\to_\cS z$ and $g\colon y\to_\cS z$ in $\cS$.
  \end{quote}
  Since $x$ and $y$ are $\cS$\nbd-normal forms, $f$ and $g$ must be in $\cE$, so that $x$ and $y$ belong to the same connected component of $\cE$.
\end{proof}

\begin{proposition}
  \label{prop:coherence_from_convergence}
  If $\cS$ is convergent modulo $\cE$, then $\cS^\top$ is coherent modulo $\cE$.
\end{proposition}

\begin{proof}
  Thanks to the Church--Rosser property (see the proof of \cref{prop:pi_zero_from_normal_forms}), for any endomorphism ${f\colon x\to_{\cS^\top}x}$ in $\cS^\top$ there exist $f_1,f_2\colon x\to_\cS z$ morphisms in $\cS$ such that $f= f_1^{-1}\circ f_2$.
  Since $\cS$ is terminating modulo $\cE$, there exists an $\cS$\nbd-normal form $y$ and a morphism $h\colon z\to_\cS y$. Because $\cS$ is confluent, the branching $(h\circ f_1,h\circ f_2)$ admits a confluence. Since $y$ is an $\cS$\nbd-normal form, this confluence is in $\cE$, and there exists $e\in\cE$ such that $e\circ h\circ f_1= h\circ f_2$. Setting $g\coloneqq h\circ f_1$ gives $e\circ g= g\circ f$, which concludes.
\end{proof}

An \emph{abstract equivalence} is an equivalence relation $\equiv$ on the set of morphisms, such that if $f\equiv f'$, then (i) $s(f)=s(f')$ and $t(f)=t(f')$, and (ii) if $g$ and $h$ are suitably composable morphisms, then $g\circ f\circ h\equiv g\circ f'\circ h$.
The notion of coherence modulo $\cE$ generalizes to \emph{$\equiv$-coherence modulo $\cE$} by further imposing that $f\equiv g^{-1}\circ e\circ g$ in the notation of \cref{defn:general_coherence_modulo}; similarly, confluence generalizes to \emph{$\equiv$-confluence}.
The results of this section extend verbatim to these generalized notions.

Beyond the discrete equivalence $\equiv_{\textrm{dis.}}$ (if two morphisms have the same source and target, they are equivalent), the scalar equivalence $\equiv_{\textrm{scl.}}$ (two morphisms $f$ and $g$ with the same source and target are equivalent if and only if $\scl(f)=\scl(g)$) will also be relevant to our purpose, in particular when analyzing scalars in the modulo data (see \cref{subsec:rewriting_foam_coherence_iso}).
In this case, the notion of $\equiv_{\textrm{scl.}}$-coherence modulo amounts to the following:

\begin{definition}
  \label{defn:scalar-coherence}
  Let $\cP$ be a scalar groupoid with scalars $\scl\colon\cP_1\to\Bbbk$ and $\cE\subset\cP$ a subgroupoid.
  We say that $\cP$ is \emph{scalar-coherent} if for every endomorphism $f\colon x\to_\cP x$, there exist morphisms $g\colon x\to_\cP y$ and $e\colon y\sim_\cE y$ such that $\scl(f)= \scl(e)$.
\end{definition}

\subsection{Abstract rewriting modulo}
\label{subsec:abstract_rewriting_modulo}

Recall that a 1\nbd-poly\-graph is a 1-globular set; in the context of rewriting, a 1-polygraph is called an \emph{abstract rewriting system (\ARS{})}.
Unpacking the definition, an \ARS{} is the data of a set $X$ of 0-cells, called \emph{elements}, and a set $\aP$ of 1-cells, called \emph{rewriting steps}, equipped with \emph{source} and \emph{target maps}:
\[\begin{tikzcd}[cramped]
  X &  \aP
  \arrow["{t}", shift left, from=1-2, to=1-1]
  \arrow["{s}"', shift right, from=1-2, to=1-1]
\end{tikzcd},\]
Note that we abuse notation and denote $\aP$ both the 1\nbd-poly\-graph and the set of rewriting steps.

Recall the notion of the free category $\aP^*$ generated by a 1\nbd-poly\-graph $\aP$ (\cref{subsubsec:2polygraphs}). We denote $\aP^\top$ the localisation of $\aP^*$.
A morphism in $\aP^*$ (resp.\ in $\aP^\top$) is called a \emph{rewriting sequence} (resp.\ a \emph{congruence}). If a rewriting sequence (resp.\ a congruence) decomposes in $n$ $\aP$\nbd-rewriting steps, we call the number $n$ its \emph{length}.
We write $x\to_\aP y$ to refer to an unspecified rewriting step with source $x$ and target $y$ (note that there could be more than one rewriting step between given source and target), or to indicate the existence of such a rewriting step.
Similarly, we write $x\overset{*}{\to}_\aP y$ (resp.\ $x\sim_\aP y$) to denote a rewriting sequence (resp.\ a congruence), and say that $x$ \emph{rewrites into} (resp.\ \emph{is congruent to}) $y$.
\begin{IEEEeqnarray*}{CcCcC}
  x\to_\aP y && x\overset{*}{\to}_\aP y && x\sim_\aP y\\*
  \text{\scriptsize rewriting step}
  &\mspace{70mu}&
  \text{\scriptsize rewriting sequence}
  &\mspace{70mu}&
  \text{\scriptsize congruence}
\end{IEEEeqnarray*}
Following these notations, $x\overset{*}{\to}_\aP y$ (resp.\ $x\to_{\aP^*} y$) coincides with $x\sim_\aP y$ (resp.\ $x\to_{\aP^\top} y$).

\begin{definition}
  \label{defn:ARSM}
  An \emph{abstract rewriting system modulo (\ARSM{})} is the data $\aS={(\aR,\aE)}$ of two 1\nbd-poly\-graphs $\aR\coloneqq(X,\aR)$ and $\aE=(X,\aE)$.
\end{definition}

Given $\aS=(\aR,\aE)$ as above, we define the following 1\nbd-poly\-graph:
\[
  {}_\aE\aR_\aE\coloneqq \aE^\top\times_X\aR\times_X\aE^\top.
\]
In other words, a rewriting step in ${}_\aE\aR_\aE$ is a triple $(e,r,e')\in \aE^\top\times \aR\times \aE^\top$ with $t(e)=s(r)$ and $t(r)=s(e')$. The source and target maps are defined as $s(e,r,e')=s(e)$ and $t(e,r,e')=t(e')$.
An \emph{$\aS$\nbd-rewriting sequence} is either an $\aE$\nbd-congruence (in which case it has \emph{length zero}) or an ${}_\aE\aR_\aE$\nbd-rewriting sequence (in which case it has the same length as the ${}_\aE\aR_\aE$\nbd-rewriting sequence, i.e.\ the number of $\aR$-rewriting steps).%
\footnote{
  \label{footnote:ARSM-distinction-dupont-S-rewriting-sequence}
  Our terminology differs from \cite{Dupont_RewritingModuloIsotopies_2022,Dupont_RewritingModuloIsotopies_2021}, where $\aS$\nbd-rewriting sequences coincide with ${}_\aE\aR_\aE$\nbd-rewriting sequences. In particular, our notions of branching and confluence do not explicitly depend on the modulo data, while the notion of termination does; this is the converse of \cite{Dupont_RewritingModuloIsotopies_2022,Dupont_RewritingModuloIsotopies_2021}.
}
We denote $\aS^*$ the set of $\aS$\nbd-rewriting sequences. Note that $({}_\aE\aR_\aE)^*\cup\aE^\top=\aS^*$.
We similarly denote ${x\to_\aS y}$ (resp.\ ${x\overset{*}{\to}_\aS y}$, resp.\ ${x\sim_\aS y}$) an $\aS$\nbd-rewriting step (resp.\ an $\aS$\nbd-rewriting sequence, resp.\ an $\aS$\nbd-congruence).
 
An \ARSM{} inherits all the notions and results introduced in the previous section, setting \ref{data_convergence_modulo} as
\[
  \cS\coloneqq \aS^*
  \quad\an\quad
  \cE\coloneqq\aE^\top. 
\]
For instance, we say that $\aS$ is \emph{terminating modulo $\aE$} (or simply \emph{terminating}) if $\cS$ is terminating modulo $\cE$; equivalently, $\aS$ is terminating if there is no infinite sequence $\{f_n\}_{n\in\bN}$ such that $f_n$ is an $\aS$\nbd-rewriting step and $t(f_n)=s(f_{n+1})$.
(Note that with our conventions, one \emph{cannot} replace ``step'' by ``sequence'' in the above characterization.)

In contrast with the previous section, an \ARSM{} provides a notion of \emph{locality}.
A \emph{local $\aS$\nbd-branching} is a pair $(f,g)$ with $f\colon x\to_\aS y$ and $g\colon x\to_\aS y'$ two rewriting steps in $\aS$:
\begin{IEEEeqnarray*}{CcCcC}
  \begin{tikzcd}[ampersand replacement=\&,cramped,row sep=.7em]
    \& y \\
    x \& \\
    \& {y'}
    \arrow["*"{description},"f", curve={height=-9pt}, from=2-1, to=1-2,"\aS"{subscript}]
    \arrow["*"{description},"{g}"', curve={height=9pt}, from=2-1, to=3-2,"\aS"{subscript}]
  \end{tikzcd}
  &&
  \begin{tikzcd}[ampersand replacement=\&,cramped,row sep=.7em]
    \& y \\
    x \& \\
    \& {y'}
    \arrow["f", curve={height=-9pt}, from=2-1, to=1-2,"\aS"{subscript}]
    \arrow["{g}"', curve={height=9pt}, from=2-1, to=3-2,"\aS"{subscript}]
  \end{tikzcd}
  &&
  \begin{tikzcd}[ampersand replacement=\&]
    x \dar[snakecd,"e"']\arrow[r,"\aR"{subscript},"f"] \&y
    \\
    x' \rar["\aR"{subscript},"g"] \&y'
  \end{tikzcd}
  \\*
  \text{\scriptsize an $\aS$\nbd-branching}
  &\mspace{70mu}&
  \text{\scriptsize a local $\aS$\nbd-branching}
  &\mspace{70mu}&
  \text{\scriptsize an $\aS$\nbd-local triple}
\end{IEEEeqnarray*}
We say that $\aS$ is \emph{locally confluent} if all local branchings are confluent.
An \emph{$\aS$\nbd-local triple} is a triple $[f,e,g]$ with $\aR$\nbd-rewriting steps $f$ and $g$ and $\aE$\nbd-congruence $e$, such that $s(f)=s(e)$ and $t(e)=s(g)$. Note that every $\aS$\nbd-local triple defines a local $\aS$\nbd-branching $(f,g\circ e)$, and that if every $\aS$\nbd-local triple is confluent, then $\aS$ is confluent.
Local $\aS$\nbd-branchings and $\aS$\nbd-local triples should be thought of as essentially identical notions, one being more suited for general statements while the other being better suited for explicit computations.

The next subsections describe how confluence can be achieved from a local analysis.

\subsubsection{Tamed congruence }
\label{subsubsec:ARSM_tamed_newmann_lemma}

Newman's lemma \cite{Newman_TheoriesCombinatorialDefinition_1942} readily extends to rewriting modulo:

\begin{lemma}[Newman's lemma modulo]
  \label{lem:ARSM_newmann_lemma}
  Let $\aS=(\aR,\aE)$ be an \ARSM{}. If $\aS$ is terminating modulo $\aE$ and locally confluent, then it is confluent.
\hfill\qed
\end{lemma}

However, local confluence does not always suffice, especially in the linear context. In this subsection, we introduce the weaker notion of \emph{$\succ$\nbd-tamed congruence}, for which an analogue of Newman's lemma still holds.

\begin{definition}
  \label{defn:ARS_esp_order}
  Let $(X,\aE)$ be an \ARS{}. An \emph{(abstract) esp-order} $\succ$ is a strict partial order on $X$ which is \emph{$\aE$-invariant}, that is, such that $(x'\sim_\aE x\an x\succ y\an y\sim_\aE y')$ implies $(x'\succ y')$.
\end{definition}

The prefix ``esp'' stands for ``$\aE$-invariant'', ``Strict'' and ``Partial''.
If $M\subset X$ is a set of elements in $X$, we write $x\succ M$ if $x\succ y$ for all $y\in M$.
If $f=f_n\circ\ldots\circ f_1$ is a sequence of composable arrows on $X$, we write $x\succ f$ to mean $x\succ\{s(f_1),t(f_1),\ldots,t(f_n)\}$.

\begin{definition}
  \label{defn:ARSM_tame_congruence}
  Let $\aS=(\aR,\aE)$ be an \ARSM{} on the set $X$ and $\succ$ an esp-order on $X$.
  An $\aS$\nbd-branching $(f,g)$ of source $\bullet$ is said to be \emph{$\succ$\nbd-tamely congruent} (resp.\ \emph{$\succ$\nbd-tamely confluent}) if there exists a congruence $h$ (resp.\ con\-fluence $(f',g')$) such that $\bullet\succ h$ (resp.\ $\bullet\succ f'^{-1}\circ g'$).
\end{definition}
 
In particular, $\succ$\nbd-tameness implies $\bullet\succ t(f)$ and $\bullet\succ t(g)$.
Here is a schematic for a $\succ$\nbd-tamed congruence, where horizontal positions are used to suggest relative orderings with respect to $\succ$:
\begin{gather*}
  \begin{tikzpicture}[scale=.7]
    \node (C) at (0,0) {$\bullet$};
    \coordinate (A1) at (2,1);
    \coordinate (A2) at (2+.2,1-.4);
    \coordinate (A3) at (2-.6,1-2*.4);
    \coordinate (A4) at (2+1.2,1-3*.4);
    \coordinate (A5) at (2-.2,1-4*.4);
    \coordinate (A6) at (3,-1);
    \draw (A1) to (A2);
    \draw (A2) to (A3);
    \draw[
      dotted,
      dash pattern=on 3pt off 3pt, 
    ] (A3) to (A4);
    \draw (A4) to (A5);
    \draw (A5) to (A6);
    \draw[->] (C) to[out=60,in=180] node[above]{$f$} (A1);
    \draw[->] (C) to[out=-60,in=180] node[below]{$g$} (A6);
    \draw[dotted,thick] (.5,2) to (.5,-1.5);
    \node at (1,1.8) {$\bullet\succ$};
  \end{tikzpicture}
\end{gather*}
In hindsight, tamed congruence is reminiscent of older ideas, as similar notions appeared already in the work of Bergman \cite{Bergman_DiamondLemmaRing_1978} (``ambiguities resolvable relative to an order'') and in the work of Winkler and Buchberger \cite{WB_CriterionEliminatingUnnecessary_1986}\footnote{We thank Vincent van Oostrom for pointing out this reference.} (``generalized Newman's lemma'').
It is also reminiscent of the notion of \emph{confluence by decreasingness} as introduced by van Oostrom \cite{vanOostrom_ConfluenceDecreasingDiagrams_1994} and extended modulo in \cite{Fv_ProofOrdersDecreasing_2013,Felgenhauer_ConfluenceTermRewriting_2015}.
Tamed congruence is a very special case, choosing the partial order $(h_1\succ h_2)\Leftrightarrow(s(h_1)\succ s(h_2)\an s(h_1)\succ t(h_2))$ on $\aS$-rewriting steps (see \cite[Definition~2]{vanOostrom_ConfluenceDecreasingDiagrams_1994} and \cite[Theorem~31]{Fv_ProofOrdersDecreasing_2013}).
We do not know whether the comparison between tamed congruence and confluence by decreasingness extends further, starting with transitivity (\cref{lem:tamed_is_transitive}).

\begin{definition}
  \label{defn:ARSM_compatible_preorder}
  Let $\aS$ be an \ARSM{}. An esp-order $\succ$ on $X$ is said to be \emph{compatible with $\aS$} if $x\to_\aS y$ implies $x\succ y$.
\end{definition}

A compatible esp-order $\succ$ controls termination: minimal elements for $\succ$ are $\aS$\nbd-normal forms and if $\succ$ is well-founded, then $\aS$ terminates modulo $\aE$.
Moreover, if $\succ$ is compatible then confluence implies $\succ$\nbd-tamed congruence.

\begin{lemma}[tamed Newman's lemma modulo]
  \label{lem:ARSM_tamed_newmann_lemma}
  Let $\aS=(\aR,\aE)$ be an \ARSM{} and $\succ$ a compatible esp-order. If $\succ$ is well-founded and every local $\aS$\nbd-branching is $\succ$\nbd-tamely congruent, then $\aS$ is confluent.
\end{lemma}

\begin{proof}
  Since $\succ$ is well-founded, we can proceed by induction on $\succ$ (see e.g.\ \cite[section 1.3.9]{ABG+_PolygraphsRewritingHigher_2025}) 
  to show that the following property holds for every $x\in X$:
  \begin{quote}
    $P(x)$: every $\aS$\nbd-branching with source $y$ such that $x\succ y$ is $\aS$\nbd-confluent.
  \end{quote}
  We only sketch the remaining of the proof. Let $x\in X$ and assume $P(y)$ holds for all $x\succ y$. Using induction, one first argues that any tamely congruent local branching with source $x$ is in fact confluent.
  Then if $(f,g)$ is a branching with source $x$, decompose $f$ and $g$ into $\aS$\nbd-rewriting steps $f_m\circ\ldots\circ f_1$ and $g_n\circ\ldots\circ g_1$, apply the previous sentence to the local branching $(f_1,g_1)$ to get a confluence, and complete the rest of the diagram using induction:
  \[\begin{tikzcd}[ampersand replacement=\&,row sep=1em,column sep=4em]
    \& \cdot \& \cdot \& \cdot \\
    \cdot \& \cdot \& \cdot \& \cdot \& \cdot \\
    \cdot \& \cdot \& \cdot \& \cdot \& \cdot \\
    \& \cdot \& \cdot \& \cdot
    \arrow[equals, from=1-2, to=2-2]
    \arrow[equals, from=3-2, to=4-2]
    \arrow[from=2-1, to=3-1,equals]
    \arrow[from=2-3, to=3-3,equals]
    \arrow[from=3-4, to=4-4,equals]
    \arrow[from=1-4, to=2-4,equals]
    \arrow[from=2-5, to=3-5,equals]
    \arrow["{*}"{description},"\aS"{subscript}, from=2-3, to=2-4]
    \arrow["{*}"{description},"\aS"{subscript}, from=3-3, to=3-4]
    \arrow["{*}"{description},"\aS"{subscript}, from=2-4, to=2-5]
    \arrow["{*}"{description},"\aS"{subscript}, from=3-4, to=3-5]
    \arrow["{*}"{description},"\aS"{subscript}, from=2-2, to=2-3]
    \arrow["{*}"{description},"\aS"{subscript}, from=3-2, to=3-3]
    \arrow["{*}"{description},"\aS"{subscript}, from=1-3, to=1-4]
    \arrow["{*}"{description},"\aS"{subscript}, from=4-3, to=4-4]
    \arrow["{\text{induction}}"{description}, draw=none, from=1-2, to=2-4]
    \arrow["{\text{induction}}"{description}, draw=none, from=4-2, to=3-4]
    \arrow["{\text{induction}}"{description}, draw=none, from=2-3, to=3-5]
    \arrow["{\text{first part of the argument}}"{description}, draw=none, from=2-1, to=3-3]
    \arrow["{f_1}", from=2-1, to=2-2]
    \arrow["{g_1}"', from=3-1, to=3-2]
    \arrow["{f_m\circ\ldots\circ f_2}", from=1-2, to=1-3]
    \arrow["{g_n\circ\ldots\circ g_2}"', from=4-2, to=4-3]
  \end{tikzcd}\]
  This concludes.
\end{proof}

A first reason to prefer tamed congruence to confluence is the following:

\begin{lemma}[transitivity of tamed congruence]
  \label{lem:tamed_is_transitive}
  Let $\aS$ be an \ARSM{} and $\succ$ an esp-order on $X$.
  Let $(f_1,f_2)$ and $(f_2,f_3)$ be two $\aS$-branchings.
  If $(f_1,f_2)$ and $(f_2,f_3)$ are $\succ$-tamely $\aS$-congruent, then $(f_1,f_3)$ is $\succ$-tamely $\aS$-congruent.
  \hfill\qed
\end{lemma}

Note that in contrast, confluence is \emph{not} a transitive property.
As we will see in \cref{sec:rewriting_foam}, transitivity is of great use in the context of higher rewriting.
In hindsight, one can find hints of transitivity in related earlier notions (see the discussion after \cref{defn:ARSM_tame_congruence}); see e.g.\ the ``triangle lemma'' as reviewed in \cite{BD_AlgebraicOperads_2016}.

\subsubsection{Branchwise congruences}
\label{subsubsec:ARSM_E_congruence}

In principle, understanding branchings when rewriting modulo is a difficult task: given that the length of the $\aE$\nbd-congruence in a local triple is arbitrary, the number of local triples is in general infinite.
To circumvent this problem, rewriting steps need to be understood \emph{up to $\aE$\nbd-congruence}, similarly to how elements in $X$ are understood up to $\aE$\nbd-congruence.
In practice, one has canonical ways to do so, coming from naturality axioms with regard to the modulo: naturality of interchangers when working modulo interchange (see \cref{subsubsec:scalar_HRSM}), naturality of the pivotal structure when working up to isotopies (see \cref{lem:foam_pivotal_naturalities}), and so on. We now formalize this situation.

\medbreak

Let $\aS=(\aR,\aE)$ be an \ARSM{}.
Two $\aS$\nbd-rewriting sequences $f$ and $g$ are said to be $\aE$\nbd-\emph{congruent} if there exist $\aE$\nbd-congruences $e_s\colon s(f)\sim_\aE s(g)$ and $e_t\colon t(f)\sim_\aE t(g)$:
\begin{IEEEeqnarray*}{CcCcC}
  \begin{tikzcd}[ampersand replacement=\&]
    \cdot \& \cdot \\
    {\cdot} \& {\cdot}
    \arrow[""{name=0, anchor=center, inner sep=0},"*"{description},"f"{yshift=.5ex},"\aS"{subscript}, from=1-1, to=1-2]
    \arrow[""{name=1, anchor=center, inner sep=0},"*"{description},"g"'{yshift=-.5ex},"\aS"{subscript}, from=2-1, to=2-2]
    \arrow["{e_s}"',snakecd, no head, from=1-1, to=2-1]
    \arrow["{e_t}",snakecd, no head, from=1-2, to=2-2]
    \arrow["\scriptstyle"{marking,allow upside down}, draw=none, from=2-1, to=1-2]
  \end{tikzcd}
  &&
  \begin{tikzcd}[ampersand replacement=\&,cramped,row sep=.7em]
    \& \cdot \\
    \cdot \& \\
    \& {\cdot}
    \arrow["*"{description},"f", curve={height=-9pt}, from=2-1, to=1-2,"\aS"{subscript,pos=.95}]
    \arrow["*"{description},"{g}"', curve={height=9pt}, from=2-1, to=3-2,"\aS"{subscript,pos=.91},""{name=0, anchor=center, inner sep=0,pos=.4}]
    \arrow["{e}", from=1-2, to=3-2,snakecd]
    \arrow[from=1-2,to=0,draw=none,"\scriptstyle"{marking,allow upside down}]
  \end{tikzcd}
  &&
\begin{tikzcd}[ampersand replacement=\&,cramped]
	\& \cdot \\
	\cdot \& \cdot \\
	\cdot \& \cdot \\
	\& \cdot
	\arrow["*"{description},""{name=0, anchor=center, inner sep=0}, "{g}", curve={height=-6pt}, from=2-1, to=1-2,"\aS"{subscript,pos=.85}]
	\arrow["*"{description},"f"'{yshift=-.5ex}, from=2-1, to=2-2,"\aS"{subscript}]
	\arrow["*"{description},"f'"{yshift=.5ex}, from=3-1, to=3-2,"\aS"{subscript}]
	\arrow["*"{description},""{name=1, anchor=center, inner sep=0}, "{g'}"', curve={height=6pt}, from=3-1, to=4-2,"\aS"{subscript,pos=1}]
	\arrow[from=2-1, to=3-1,equals]
	\arrow[from=1-2, to=2-2,snakecd]
	\arrow[from=4-2, to=3-2,snakecd]
	\arrow["\scriptstyle"{marking,allow upside down}, draw=none, from=0, to=2-2]
	\arrow["\scriptstyle"{marking,allow upside down}, draw=none, from=1, to=3-2]
\end{tikzcd}
  \\*
  \substack{%
    \text{\scriptsize two $\aE$\nbd-congruent}\\
    \text{\scriptsize $\aS$\nbd-rewriting sequences}
  }
  &\mspace{70mu}&
  \substack{%
    \text{\scriptsize a $\aE$\nbd-congruent}\\
    \text{\scriptsize $\aS$\nbd-branching}
  }
  &\mspace{70mu}&
  \substack{%
    \text{\scriptsize two branchwise $\aE$\nbd-congruent}\\
    \text{\scriptsize $\aS$\nbd-branchings}
  }
\end{IEEEeqnarray*}
An $\aS$\nbd-branching $(f,g)$ is said to be \emph{$\aE$\nbd-con\-gruent} if there exists an $\aE$\nbd-con\-gruent $e\colon t(f)\sim_\aE t(g)$.
Two $\aS$\nbd-branchings $(f,f')$ and $(g,g')$ are said to be \emph{branchwise $\aE$\nbd-congruent} if they have the same source and $(f,g)$ (resp.\ $(f',g')$) is $\aE$\nbd-congruent.

The following lemma states that ``confluence is preserved under branchwise $\aE$\nbd-congruence'':

\begin{lemma}
  \label{lem:ARSM_branchwise_confluence_lemma}
  
  Let $\aS=(\aR,\aE)$ be an \ARSM{}.
  If $(f,g)$ and $(f',g')$ are branchwise ${\aE}$\nbd-cong\-ruent $\aS$\nbd-branchings, then $(f,g)$ is confluent if and only if $(f',g')$ is.
  \hfill\qed
\end{lemma}

Branchwise $\aE$\nbd-congruence generalizes to branchwise $\succ$\nbd-tamed congruence:

\begin{definition}
  \label{defn:branchwise_confluence_congruence}
  Let $\aS=(\aR,\aE)$ be an \ARSM{} and $\succ$ an esp-order on $X$.
  Two $\aS$\nbd-branchings $(f,g)$ and $(f',g')$ are \emph{branchwise $\succ$\nbd-tamely cong\-ruent} if they have the same source and the branchings $(f,f')$ and $(g,g')$ are respectively $\succ$\nbd-tamely congruent.
\end{definition}

Contrary to branchwise $\aE$\nbd-congruence, working up to branchwise confluence does not preserve confluence. However, the following lemma states that ``tamed congruence is preserved under branchwise tamed congruence''. Having applications in mind, we state it with respect to a sub-\ARSM{}.

\begin{lemma}
  [\textsc{Branchwise Tamed Congruence Lemma}]
  \label{lem:ARSM_branchwise_tamed_congruence_lemma}
  Let $\aS=(\aR,\aE)$ be an \ARSM{} and $\succ$ an esp-order on $X$ compatible with $\aS$.
  Let also $\aT\subset\aS$ be a sub-\ARSM{}.
  If $(f,g)$ and $(f',g')$ are branchwise $\succ$\nbd-tamely ${\aT}$\nbd-con\-gruent $\aS$\nbd-branchings, then $(f,g)$ is ${\aT}$\nbd-tamely ${\aT}$\nbd-congruent if and only if $(f',g')$ is.\hfill\qed
\end{lemma}

The proof of the above lemma fits into one picture, and essentially comes down to transitivity (\cref{lem:tamed_is_transitive}):
\begin{center}
  \begin{tikzcd}[ampersand replacement=\&,cramped,anchor=south]
    \cdot \& \cdot \& \cdot \& \cdot \\
    \& \cdot \& \cdot
    \arrow["{f'}"', from=1-2, to=1-1,"\aS"'{subscript}]
    \arrow["{f}", from=1-2, to=2-2,"\aS"'{subscript}]
    \arrow[Rightarrow, no head, from=1-3, to=1-2]
    \arrow["{g'}", from=1-3, to=1-4,"\aS"{subscript}]
    \arrow["{g}"', from=1-3, to=2-3,"\aS"{subscript}]
    \arrow["{*}"{description}, tail reversed, from=1-4, to=2-3,"\aT"'{subscript,pos=.5}]
    \arrow["{*}"{description}, tail reversed, from=2-2, to=1-1,"\aT"'{subscript,pos=.5}]
    \arrow["{*}"{description}, tail reversed, from=2-3, to=2-2,"\aT"'{subscript,pos=.5}]
  \end{tikzcd}
\end{center}

\subsection{Higher rewriting modulo}
\label{subsec:higher_rewriting_modulo}

\subsubsection{Higher rewriting system modulo}
\label{subsubsec:HRSM_definition}

Recall the notion of 3\nbd-sesqui\-poly\-graph (\cref{subsubsec:3sesquipolygraphs}).

\begin{definition}
  \label{defn:HRSM}
  A \emph{higher rewriting system modulo} (\HRSM{}) $\sS=(\sR,\sE)$ is the data of two 3\nbd-sesqui\-poly\-graphs $\sR$ and $\sE$ with the same underlying 2\nbd-sesqui\-poly\-graph $\sR_{2}=\sE_2$.
\end{definition}

If $\sE_3=\emptyset$, we simply say \emph{higher rewriting system}.
Recall the notion of a Gray polygraph \cref{subsubsec:gray_polygraph}. In most cases, interchangers are part of the modulo data.
This leads to the following notion:

\begin{definition}
  \label{defn:GRSM}
  A \emph{Gray rewriting system modulo} (\GRSM{}) is the data of a \HRSM{} $\sS=(\sR,\sE)$ such that $\sE$ is a Gray polygraph.
\end{definition}

Recall the notion of contexts (\cref{subsubsec:contexts}). Given a \HRSM{} ${\sS=(\sR,\sE)}$ every choice of 1-sphere $\square$ in $\sR_1^*$ defines an \ARSM{} $\aS(\square)=(\aR(\square),\aE(\square))$ on the underlying set $X(\square)$, where
\begin{gather*}
  X(\square)=\sR_2^*(\square)=\sE_2^*(\square),
  \\
  \aR(\square)=\Cont(\sR_3)(\square)
  \quad\an\quad
  \aE(\square)=\Cont(\sE_3)(\square).
\end{gather*}
Moreover, every context $\Gamma$ on $\square$ defines a morphism of \ARSM{}
\begin{equation}
  \label{eq:context_as_morphism}
  \Gamma\colon\aS(\square)\to \aS\big(s_1(\Gamma),t_1(\Gamma)\big).
\end{equation}
(A morphism of \ARSM{} is a pair of morphisms of 1-globular sets, the latter denoting a pair of functions commuting with the source and target maps; see \cref{subsubsec:n-globular_sets}.)
We say that a branching $(f,g)$ is a \emph{contextualization} of another branching $(f',g')$ whenever $\Gamma[f',g']=(f,g)$.

We think of a \HRSM{} as a category of \ARSM{}s, with morphisms given by contexts. In principle, one could deal with each \ARSM{} $\aS(\square)$ independently, using the tools of abstract rewriting theory modulo (\cref{subsec:abstract_rewriting_modulo}).
In practice, one would like to leverage the fact that these \ARSM{}s gather together in a category, and relate to one another (and themself) via contextualization.
In that sense, higher rewriting theory (modulo) is nothing else than the study of morphisms of abstract rewriting systems (modulo).

\medbreak

In many practical situations, termination can only be obtained with context-dependent termination rule, where $r$ being a rewriting rule does not imply that $\Gamma[r]$ is.
For that reason, we introduce the following notion:

\begin{definition}
  Let $\sS=(\sR,\sE)$ be a \HRSM{}.
  An \emph{abstract sub-system $\aT$} is the data of a family of subsets $\aT(\square)\subset\sR_3^*(\square)$ for each 1-sphere $\square$ in $\sR_1^*$.
\end{definition}

We think of $\aT$ as a family of sub-\ARSM{}s $\aT(\square)\subset\aS(\square)$ modulo $\sE(\square)$ and on the set $X(\square)$, and write $\aT\subset\sS$ to emphasize this point.
We say that $\aT$ is \emph{confluent} if every sub-\ARSM{} $\aT(\square)$ is confluent, and similarly for other notions.
One can always consider $\sS$ as its own abstract sub-system. We refer to the case $\aT=\sS$ as \emph{context-agnostic}, and to the case  $\aT\subsetneq\sS$ as \emph{context-dependent}.

\medbreak

Abstract notions of orders generalize verbatim:

\begin{definition}
  \label{defn:HRSM_compatible_preorder}
  Let $\sS=(\sR,\sE)$ be a \HRSM{} and $\aT\subset\sS$ an abstract sub-system.
  An \emph{esp-order} is an esp-order on every \ARSM{} $\aS(\square)$, in the sense of \cref{defn:ARS_esp_order}.
  An esp-order is said to be \emph{compatible with $\aT$} if it is compatible with every sub-\ARSM{} $\aT(\square)$, in the sense of \cref{defn:ARSM_compatible_preorder}.
\end{definition}

\subsubsection{Contextualization}
\label{subsubsec:HRSM_contextualization}
\begin{lemma}
  \label{lem:HRSM_contextualization_branching_lemma}
  Let $\sS=(\sR,\sE)$ be a \HRSM{} and $\aT\subset\sS$ an abstract sub-system.
  Let $(f,g)$ be a local $\aT$\nbd-branching and $\Gamma$ a context. Assume $(f,g)$ admits a $\aT$\nbd-confluence $(f',g')$ such that $\Gamma[f',g']$ belongs to $\aT$. Then $\Gamma[f,g]$ is $\aT$\nbd-confluent.
  \hfill\qed
\end{lemma}

In particular, in the context-agnostic case:

\begin{lemma}
  \label{lem:HRSM_contextualization_branching_lemma_agnostic}
  Let $\sS=(\sR,\sE)$ be a \HRSM{}.
  Let $(f,g)$ be a local $\sS$\nbd-branching and $\Gamma$ a context. If $(f,g)$ is $\sS$\nbd-confluent, then $\Gamma[f,g]$ is $\sS$\nbd-confluent.
  \hfill\qed
\end{lemma}

\subsubsection{Independent branching}
\label{subsubsec:HRSM_independent_branching}

The higher structure of a higher rewriting system $\sP$ induces canonical types of local $\sP$\nbd-bran\-chings. A local $\sP$\nbd-branching is \emph{independent} if it is of the form $(f,g)=(\ssf\starop_1\psi,\phi\starop_1 \ssg)$, for $\sP$\nbd-rewriting steps $\ssf\colon \phi\to\phi'$ and $\ssg\colon\psi\to\psi'$:
\[\begin{tikzcd}[column sep=large,row sep=-1.2em]
  &
  \satex{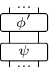}
  \\
  \satex{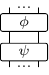}
  &
  \\
  &
  \satex{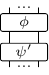}
  \tar[from=2-1,to=1-2,"\ssf\starop_1\psi"{yshift=3pt},curve={height=-20pt}]
  \tar[from=2-1,to=3-2,"\phi\starop_1 \ssg"'{yshift=-3pt},curve={height=20pt}]
\end{tikzcd}\]
When working modulo, that is with a \HRSM{} $\sS=(\sR,\sE)$, independent branchings refer to independent $\sR$\nbd-branchings. 
Note that the property of being an independent branching is preserved by context: if $(f,g)$ is an independent branching, so is $\Gamma[f,g]$.
A local $\sS$\nbd-branching (or an $\sS$\nbd-local triple) is said to be \emph{overlapping} if it is not branchwise $\sE$\nbd-congruent to an independent $\sR$\nbd-branching.

An independent branching always have a canonical $\sS$\nbd-confluence, given by
\[(g',f')\coloneqq(\phi'\starop_1 \ssg,\ssf\starop_1\psi').\]
The following is tautological, thanks to the independence axiom:

\begin{lemma}
  Let $\sS=(\sR,\sE)$ be a \HRSM{} and $\aT\subset\sS$ an abstract sub-system.
  Let $(f,g)$ is an independent $\aT$\nbd-branching. If its canonical $\aT$\nbd-confluence $(g',f')$ is in $\aT$, then it defines a ${\aT}$\nbd-conf\-luence for $(f,g)$.
  \hfill\qed
\end{lemma}

In particular, in the context-agnostic case:

\begin{lemma}
  \label{lem:HRSM_independent_branching_agnostic}
  Let $\sS=(\sR,\sE)$ be a \HRSM{}. Every independent $\sS$\nbd-bran\-ching is $\sS$\nbd-conf\-luent.
  \hfill\qed
\end{lemma}

\subsubsection{Higher Newman's lemma}

In the context-agnostic case, the above results can be gathered in a single black-box:

\begin{lemma}[higher Newman's lemma]
  \label{lem:HRSM_higher_newmann_lemma}
  Let $\sS=(\sR,\sE)$ be a \HRSM{}. If $\sS$ is terminating modulo $\sE$ and if every overlapping $\sS$\nbd-branching is, up to branchwise $\sE$\nbd-congruence, a contextualization of a $\sS$\nbd-confluent branching, then $\sS$ is convergent modulo $\sE$.
\end{lemma}

\begin{proof}
  This follows from the (abstract) Newman's lemma (\cref{lem:ARSM_newmann_lemma}), in combination with the branchwise confluence lemma (\cref{lem:ARSM_branchwise_confluence_lemma}), confluence of contextualizations (\cref{lem:HRSM_contextualization_branching_lemma_agnostic}) and confluence of independent branchings (\cref{lem:HRSM_independent_branching_agnostic}).
\end{proof}

\subsubsection{Independent rewriting}
\label{subsubsec:HRSM_independent_rewriting}

In the higher setting, one often wants to rewrite part of a diagram away from a given branching. This gives the notion of \emph{independent rewriting}.
Given rewriting steps $\ssf\colon\phi\to\phi_1$, $\ssg\colon\phi\to\phi_2$ and $\ssh\colon\psi\to\psi'$, setting
\[
  f=\ssf\starop_1\psi,
  \quad 
  g=\ssg\starop_1\psi
  \quad\an\quad
  h=\phi\starop_1\ssh
\]
defines three pairs of branchings $(f,g)$, $(f,h)$ and $(h,g)$, the latter two being independent branchings.
Define also
\[
  f'=\ssf\starop_1\psi',
  \quad 
  g'=\ssg\starop_1\psi',
  \quad
  h_1=\phi_1\starop_1\ssh
  \quad\an\quad
  h_2=\phi_2\starop_1\ssh.
\]
We say that the branching $(f,g)$ \emph{rewrites} into the branching $(f',g')$ via the triple $(h_1,h,h_2)$ (and similarly if the vertical positions of $f,g$ and $h$ are swapped), as pictured below:
\[\begin{tikzcd}[column sep=large,row sep=0em]
  &
  \satex{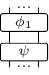}
  &&
  \satex{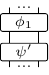}
  \\
  \satex{indep-m-l}
  &&
  \satex{indep-m-r}
  \\
  &
  \satex{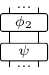}
  &&
  \satex{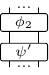}
  \tar[from=2-1,to=1-2,"f"]
  \tar[from=2-1,to=3-2,"g"']
  \tar[from=2-1,to=2-3,"h"{yshift=2pt}]
  \tar[from=2-3,to=1-4,"f'"]
  \tar[from=2-3,to=3-4,"g'"']
  \tar[from=1-2,to=1-4,"h_1"{yshift=2pt}]
  \tar[from=3-2,to=3-4,"h_2"{yshift=2pt}]
\end{tikzcd}\]

We have the following lemma:
 
\begin{lemma}
  \label{lem:HRSM_independent_rewriting_context_dependent}
  Let $\sS=(\sR,\sE)$ be a \HRSM{}, $\aT\subset\sS$ an abstract sub-system and $\succ$ an esp-order compatible with $\aT$.
  In the situation above, assume $(f,g)$ is a $\aT$\nbd-branching.
  If $h$, $h_1$ and $h_2$ are in $\aT$ and if $(f',g')$ is $\succ$\nbd-tamed $\aT$\nbd-congruent, then so is $(f,g)$.
  \hfill\qed
\end{lemma}

In particular, in the context-agnostic case:

\begin{lemma}
  \label{lem:HRSM_independent_rewriting_context_agnostic}
  Let $\sS=(\sR,\sE)$ be a \HRSM{} and $\succ$ an esp-order compatible with $\sS$. In the situation above, if $(f',g')$ is $\succ$\nbd-tamed $\sS$\nbd-congruent, then so is $(f,g)$.
  \hfill\qed
\end{lemma}

\subsubsection{Scalar higher rewriting system modulo}
\label{subsubsec:scalar_HRSM}

In the linear context, one sometimes encounters linear rewriting systems where every rewriting step is of the form $a\to \lambda b$ for $a,b\in\cB$ basis elements and $\lambda\in\Bbbk^\times$ an invertible scalar; typically, when studying the coherence of a modulo data. We call this situation \emph{monomial-invertible} (see \cref{subsubsec:LRSM}).
Instead of viewing $a\to \lambda b$ as a linear rule, we can view it as an abstract rule $a\overset{\lambda}{\to} b$, where each abstract rewriting step is equipped with the extra data of an invertible scalar. In that way, one can study scalar-coherence (\cref{defn:scalar-coherence}) with the tools of higher rewriting modulo, and avoid using \emph{linear} higher rewriting modulo.
At the structural level, this fact was already discussed in \cref{subsubsec:monomial_linear_nsesquipoly}.

Recall the notion of scalar 3\nbd-sesqui\-poly\-graph (\cref{subsubsec:scalar_3sesquipolygraphs}):

\begin{definition}
  A \emph{scalar higher rewriting system modulo} (\HRSM{}) $\sS=(\sR,\sE)$ is the data of two scalar 3\nbd-sesqui\-poly\-graphs $\sR$ and $\sE$ with the same underlying 2\nbd-sesqui\-poly\-graph $\sR_{2}=\sE_2$.
\end{definition}

More formally, all the notions and results of \cref{subsec:abstract_rewriting_modulo} and \cref{subsec:higher_rewriting_modulo} can be generalized to include a higher equivalence $\equiv$ on $\sP^*$, for $\sP$ a higher rewriting system; see \cite{Schelstraete_OddKhovanovHomology_2024} for a detailed account.
This higher equivalence is required to satisfy the following axiom:
\begin{quote}
  \begin{samepage}
  \emph{independence axiom:} for every pair of 1-composable 3-cells $\ssf\colon \phi\to\phi'$ and $\ssg\colon\psi\to\psi'$ in $\sP^*$, the 3-cells $(\phi'\starop_1 \ssg)\starop_2(\ssf\starop_1\psi)$ and $(\ssf\starop_1\psi')\starop_2(\phi\starop_1\ssg)$ are $\equiv$-equivalent:
  \end{samepage}
  \[\begin{tikzcd}[column sep=large]
    \satex{indep-m-l}
    &
    \satex{indep-bis}
    \\
    \satex{indep-m-r}
    &
    \satex{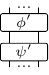}
    \tar[from=1-1,to=1-2,"\ssf\starop_1\psi"{yshift=3pt}]
    \tar[from=2-1,to=2-2,"\ssf\starop_1\psi'"'{yshift=-3pt}]
    \tar[from=1-1,to=2-1,"\phi\starop_1 \ssg"'{xshift=-1pt}]
    \tar[from=1-2,to=2-2,"\phi'\starop_1\ssg"{xshift=1pt}]
    \arrow[from=1-2,to=2-1,phantom,"\equiv"{marking,allow upside down}]
  \end{tikzcd}\]
\end{quote}
which captures the interchange of 3-cells in $\sP^*$.
The independence axiom ensures that the two sides of the canonical confluence of an independent branching are equivalent, so that \cref{lem:HRSM_independent_rewriting_context_dependent} still holds.
For our purpose, only the discrete equivalence $\equiv_{\mathrm{dis.}}$ and the scalar equivalence $\equiv_{\mathrm{scl.}}$ are relevant (see the discussion before \cref{defn:scalar-coherence}). Both equivalences satisfy the independent axiom.

\begin{definition}
  A \emph{scalar Gray rewriting system modulo} (\HRSM{}) is the data of a scalar \HRSM{} $\sS=(\sR,\sE)$ such that $\sE$ is a Gray polygraph.
\end{definition}

Again, the tools of Gray rewriting modulo generalize to include scalars.
In terms of the higher equivalence $\equiv$, what matters is that $\equiv$ encodes the naturality axioms of interchangers; this is always the case for the discrete and scalar equivalence.

\medbreak

At this point, the reader is equipped to read \cref{subsec:rewriting_foam_coherence_iso}.

\subsection{Linear rewriting modulo}
\label{subsec:linear_rewriting_modulo}

As in \cref{not:sesquicat_notation}, we fix $\Bbbk$ a commutative ring and denote $\Bbbk^\times$ its invertible elements. Given a set $\cB$, we write $\langle \cB\rangle_\Bbbk$ the free $\Bbbk$\nbd-module generated by $\cB$.

\subsubsection{Linear rewriting system modulo}
\label{subsubsec:LRSM}

Recall that linear 1\nbd-poly\-graph and linear 1-sesqui\-poly\-graph are identical notions. Unpacking the definition, a linear 1\nbd-poly\-graph is the data of a set $\cB$ of 0-cells, called \emph{monomials}, and a set $\lP$ of 1-cells, called \emph{relations}, equipped with \emph{source} and \emph{target maps} $s$ and $t$:
\[\begin{tikzcd}[cramped]
  \langle\cB\rangle_\Bbbk &  \lP
  \arrow["{t}", shift left, from=1-2, to=1-1]
  \arrow["{s}"', shift right, from=1-2, to=1-1]
\end{tikzcd}.\]
Elements in $V\coloneqq\langle\cB\rangle_\Bbbk$ are called \emph{vectors}.
The \emph{$\Bbbk$\nbd-module presented by $(\cB,\lP)$} is the module $\langle\cB\rangle_\Bbbk/\langle\lP\rangle_\Bbbk$. Note that $\langle\lP\rangle_\Bbbk$ can be viewed as a globular extension of $\langle\cB\rangle_\Bbbk$, extending $s$ and $t$ linearly.
If $\lP^=$ denotes the reflexive closure of $\lP$, we write $\lP^l\coloneqq \langle\lP^=\rangle_\Bbbk$, viewed as a globular extension:
\[\begin{tikzcd}[cramped]
  \langle\cB\rangle_\Bbbk &  \lP^l
  \arrow["{t}", shift left, from=1-2, to=1-1]
  \arrow["{s}"', shift right, from=1-2, to=1-1]
\end{tikzcd}.\]
Explicitly, relations in $\lP^l$ are of the form $\sum_i\lambda_ir_i+v$ for $\lambda_i\in\Bbbk$, $r_i\in\lP$ and $v\in V$, with source $\sum_i\lambda_is(r_i)+v$ and target $\sum_i\lambda_it(r_i)+v$.
As is explained in \cref{subsubsec:LRSM_positive_rw}, $\lP^l$ should be thought of as the set of congruences associated to $\lP$, analogous to $\aP^\top$ in the abstract case.  

We say that $\lR$ is \emph{left-monomial} if for all $r\in\lR$, we have $s(r)\in\cB$; in other words, each $r$ is of the form $b\overset{r}{\to}_\lR \sum_i\lambda_ib_i$, with $\lambda_i\in\Bbbk$ and $b,b_i\in\cB$.
We say that $\lR$ is \emph{adapted} if it is left-monomial and we have $s(r)\not\in\supp(t(r))$ for every $r\in\lR$. We refer to \cref{rem:adaptedness} for comment on the terminology and comparison with the literature.

Then:

\begin{definition}
  \label{defn:linear-rewriting-system}
  A \emph{linear rewriting system (\LRS{})} is the data $(\cB;\lP)$ of an adapted linear 1\nbd-poly\-graph $\lP$ on a set $\cB$.
\end{definition}

We sometimes leave $\cB$ implicit, and call $\lP$ a \LRS{}. 
The adaptedness condition is not an important restriction. Indeed, if $\mu b+\sum_i\mu_i b_i=0$ is some relation in a $\Bbbk$\nbd-module presentation, we can re\-write it as $b=-\mu^{-1}\sum_i\mu_i b_i$ (provided that $\mu$ is invertible). Doing so with every relation gives an adapted linear 1\nbd-poly\-graph presenting the given $\Bbbk$\nbd-module (provided we can always find such an invertible scalar $\mu$).

\medbreak

We extend to modulos the notion of \LRS{}.
Let $\lE$ be another set of linear relations on~$\cB$:
\[\begin{tikzcd}[cramped]
  \langle\cB\rangle_\Bbbk &  \lE
  \arrow["{t}", shift left, from=1-2, to=1-1]
  \arrow["{s}"', shift right, from=1-2, to=1-1]
\end{tikzcd}.\]
Denote $\Bbbk^\times\cB$ the subset of $\langle\cB\rangle_\Bbbk$ consisting of vectors of the form $\lambda b$ for $\lambda\in\Bbbk^\times$ and $b\in\cB$.
We say that $\lE$ is \emph{monomial-invertible} if it is left-monomial and for all $e\in\lE$, we have $t(e)\in\Bbbk^\times\cB$; in other words, each $e$ is of the form $b\overset{e}{\to}_\lE \lambda b'$, with $\lambda\in\Bbbk^\times$ and $b,b'\in\cB$.
This coincides with the notion of monomial-invertible linear 1\nbd-poly\-graph defined in \cref{subsubsec:monomial_linear_nsesquipoly}.
We will always assume that the modulo data is monomial-invertible.

Similarly to \ARSM{}, we write $u\sim_\lE v$ if there exists $e\in\lE^{l}$ such that $s(e)=u$ and $t(e)=v$, and in that case we say that $u$ and $v$ are \emph{$\lE$\nbd-congruent}.
In the linear context, we furthermore have a notion of \emph{projective $\lE$\nbd-congruence}, an equivalence relation on the set of monomials $\cB$, defined as $b\projrel_\lE b'$ if and only if there exists $\lambda\in\Bbbk^\times$ such that $b \sim_\lE \lambda b'$.
For $v$ a vector in $V$, we set
\[\projsupp_\lE(v)\coloneqq\left\{{b}\in\cB\mid b\projrel_\lE b'\text{ for some }b'\in\supp(v)\right\}\]
and call $\projsupp_\lE(v)$ the \emph{$\lE$\nbd-projective support of $v$}.

\begin{remark}
  \label{rem:LRSM_monomial_scalar_ARSM}
  When $\lE$ is monomial-invertible, we can alternatively treat it as a scalar 1-polygraph $\scl(\lE)$. In particular, one can use abstract rewriting rather than linear rewriting, as we explained already in \cref{subsubsec:scalar_HRSM}.
  For instance, $\lE^l$ is coherent in the sense of \cref{defn:general_coherence_modulo} if and only if $\scl(\lE)^\top$ is scalar-coherent in the sense of \cref{defn:scalar-coherence}.
\end{remark}

A set $\lR$ of linear relations on $\cB$ is said to be \emph{$\lE$\nbd-adapted} if it is left-monomial and we have $s(r)\not\in\projsupp_\lE(t(r))$ for every $r\in \lR$.

\begin{definition}
  \label{defn:linear-rewriting-system-modulo}
  A \emph{linear rewriting system modulo (\LRSM{})} is the data $\lS=(\cB;\lR,\lE)$ of two linear 1\nbd-poly\-graphs $\lR\coloneqq (\cB;\lR)$ and $\lE\coloneqq (\cB;\lE)$, such that $\lE$ is
  monomial-invertible
  and $\lR$ is $\lE$\nbd-adapted.
\end{definition}

We sometimes leave $\cB$ implicit, and call $(\lR,\lE)$ a \LRSM{}. 
Note that if $\lE=\emptyset$ (i.e.\ the set of relations in $\lE$ is empty), we recover the notion of linear rewriting system.
The \emph{module presented by $\lS$} is the module $\brak{\cB}_\Bbbk/\brak{\lR\sqcup\lE}_\Bbbk$. We write:
\[
  [-]_\lS\colon \brak{\cB}_\Bbbk\to \brak{\cB}_\Bbbk/\brak{\lR\sqcup\lE}_\Bbbk
\]
the associated quotient map. Finally, we write $\lS^l\coloneqq (\lR\sqcup\lE)^l$.

\subsubsection{Positive rewriting steps}
\label{subsubsec:LRSM_positive_rw}

Let $\lS=(\lR,\lE)$ be a \LRSM{}.
In analogy with the abstract setting, an \emph{$\lS$\nbd-rewriting step} is a composition
\[
  u\;\;\sim_\lE \;\;\lambda s(r)+v\;\;\overset{}{\dashrightarrow}_\lR \;\;\lambda t(r)+v\;\;\sim_\lE \;\;w,\quad \lambda\in\Bbbk\setminus\{0\},r\in\lR,v\in\brak{\cB}_\Bbbk.
\]
In that case, we say that the rewriting step is \emph{of type $r$}, and if $v=0$ we say that it is \emph{monomial}. We abuse notation and write $\lS$ the set of $\lS$-rewriting steps.
Contrary to the abstract case, $\lS$-rewriting steps are not suited for rewriting, as they prevent termination:

\begin{lemma}
  \label{lem:LRSM_rewriting_steps_symmetric}
  As a set of relations on $\brak{\cB}_\Bbbk$, the set of $\lS$-rewriting steps is symmetric.
\end{lemma}

\begin{proof}
  We can restrict to the non-modulo setting, i.e.\ $\lE=\emptyset$ and $\lS=\lR$, without loss of generality.
  If $\alpha = \lambda r+v$ is an $\lR$-rewriting step with the notations above, then
  \[\alpha^{-1}\coloneqq(-\lambda)r+\lambda(s(r)+t(r))+v\]
  is also an $\lR$-rewriting step with source $t(\alpha)$ and target $s(\alpha)$.
\end{proof}

To circumvent that problem, one works with so-called ``positive'' rewriting steps; see e.g.\ \cite{GHM_ConvergentPresentationsPolygraphic_2019}\footnote{This terminology appears e.g.\ in \cite{CDM_ConfluenceAlgebraicRewriting_2022} while other references \cite{GHM_ConvergentPresentationsPolygraphic_2019,Alleaume_RewritingHigherDimensional_2018,Dupont_RewritingModuloIsotopies_2022} use the terminology \emph{elementary relation} for rewriting step, and rewriting step for positive rewriting step.}.
This notion extends to rewriting modulo as follows:

\begin{definition}
  \label{defn:LRSM_positive_rewriting_steps}
  An \emph{$\lS$\nbd-rewriting step} is said to be \emph{positive} if
  \[s(r)\notin\projsupp_\lE(v).\]
  The set of positive $\lS$\nbd-rewriting steps is denoted $\lS^+$, and we sometimes abbreviate the terminology ``positive $\lS$-rewriting step'' by ``$\lS^+$-rewriting steps''.
\end{definition}

Note that the terminology ``$\lS^+$-rewriting steps'' is in line with abstract rewriting, viewing $(\brak{\cB}_\Bbbk;\lS^\st,\lE^l)$ as an abstract rewriting system.
We also say that a sequence of positive $\lS$\nbd-rewriting steps is a \emph{positive $\lS$\nbd-rewriting sequence}, or equivalently, a \emph{$\lS^+$\nbd-rewriting sequence}.

Positivity prevents the obstacle to termination described in \cref{lem:LRSM_rewriting_steps_symmetric}. In the notations of the proof, if $\alpha$ is positive then $\alpha^{-1}$ \emph{cannot} be positive; indeed, the assumptions $s(r)\notin\supp(t(r))$ (adaptedness of $\lS$) and $s(r)\notin\supp(v)$ (positiveness of $\alpha$) imply that
\[s(r)\in\supp(\lambda(s(r)+t(r))+v).\]
In diagrams, we use dashed and plain arrows to emphasize positivity:
\begin{IEEEeqnarray*}{ccc}
  \begin{tikzcd}
    s(\alpha)\arrow[dashed,r,"\lS"{subscript}] & t(\alpha)
  \end{tikzcd}
  &\mspace{80mu}&
  \begin{tikzcd}
    s(\alpha)\arrow[,r,"\lS"{subscript}] & t(\alpha)
  \end{tikzcd}
  \\*[1ex]
  \substack{
    \text{\footnotesize(not necessarily positive)}\\
    \text{\footnotesize rewriting step}
  }
  &&
  \substack{
    \text{\footnotesize positive}\\
    \text{\footnotesize rewriting step}
  }
\end{IEEEeqnarray*}
Note that a monomial rewriting step is necessarily positive.
See \cref{footnote:positivity_depends_on_modulo} for further comments on positivity and comparison with the literature.
The following lemma shows that both $\lS$\nbd-rewriting steps and positive $\lS$\nbd-rewriting steps are suited to study $\lS$\nbd-congruence:

\begin{lemma}
  \label{lem:LRSM_congruence_describe_module}
  Let $\lS=(\cB;\lR,\lE)$ be a \LRSM{}, $u,v\in\brak{\cB}_\Bbbk$ and recall the notation $[u]_\lS,[v]_\lS$ for their respective image in $\brak{\cB}_\Bbbk/\brak{\lR\sqcup\lE}_\Bbbk$. The following statements are equivalent:
  \begin{enumerate}[(i)]
    \item $[u]_\lS=[v]_\lS$ in $\brak{\cB}_\Bbbk/\brak{\lR\sqcup\lE}_\Bbbk$;
    \item $u$ and $v$ are $\lS$\nbd-congruent;
    \item $u$ and $v$ are $\lS^+$\nbd-congruent.
  \end{enumerate}
  In particular, $\lS^l=(\lS^+)^\top=\lS^\top$.
\end{lemma}

In order to prove the lemma, we need the following Church--Rosser property for positive rewriting steps, which generalizes Lemma~3.1.2 in \cite{GHM_ConvergentPresentationsPolygraphic_2019} to rewriting modulo.
The last statement is explained and used in \cref{subsubsec:LRSM_termination_order}, and can be ignored for the purpose of this subsection.

\begin{lemma}
  \label{lem:LRSM_Church-Rosser_positive_rw}
  Let $\lS=(\cB;\lR,\lE)$ be a \LRSM{}.
  If $f$ is an $\lS$\nbd-rewriting step, then there exist $g$ and $h$ either identities or positive $\lS$\nbd-rewriting steps such that $f=h^{-1}\circ g$:
  \[\begin{tikzcd}
    s(f) && t(f)\\
    & w &
    \arrow[from=1-1,to=1-3,"f",dashed,"\lS"{subscript}]
    \arrow[from=1-1,to=2-2,"g"',"\lS"{subscript}]
    \arrow[from=1-3,to=2-2,"h","\lS"'{subscript}]
  \end{tikzcd}\]
  Moreover, we have $f\relsucc w$ for any esp-order $\succ$ (see \cref{defn:LRSM_relative_relation}).
\end{lemma}

\begin{proof}
  Let $\lambda r+v$ be the $\lR$\nbd-rewriting step associated to $f$, with $r\in \lR$, $\lambda\in\Bbbk\setminus\{0\}$ and $v\in \brak{\cB}_\Bbbk$.
  Extracting $s(r)$ from the decomposition of $v$, we let $\mu\in\Bbbk$ and $v'\in \brak{\cB}_\Bbbk$ such that $v\sim_\lE\mu s(r)+v'$ and ${s(r)}\notin\projsupp_\lE(v')$.
  Set
  \[w\coloneqq (\lambda+\mu)t(r)+v'.\]
  Then $s(f)\sim_\lE(\lambda+\mu)s(r)+v'$ (resp.\ $t(f)\sim_\lE\mu s(r) + (\lambda t(r)+v')$) is either equal to $w$ if $\mu=-\lambda$ (resp.\ if $\mu=0$), or there exists a positive rewriting step $(\lambda+\mu)s(r)+v'\to_\lR w$ (resp.\ $\mu s(r) + (\lambda t(r)+v')\to_\lR w$):
  \[\begin{tikzcd}[row sep=1.5em]
    {\lambda s(r)+\mu s(r)+v'} & s(f)
    & [-12ex] & [-12ex] t(f) & {\lambda t(r)+\mu s(r)+v'}\\
    && {\lambda t(r)+\mu t(r)+v'} &&
    \arrow[from=1-1,to=1-2,snakecd,"\lE"{subscript}]
    \arrow[from=1-2,to=1-4,"f",dashed,"\lR"{subscript}]
    \arrow[from=1-4,to=1-5,snakecd,"\lE"{subscript}]
    \arrow[from=1-1,to=2-3,"="{description},"g"',"\lR"{subscript},end anchor={north west}]
    \arrow[from=1-5,to=2-3,"="{description},"h","\lR"'{subscript},end anchor={north east}]
  \end{tikzcd}\]
  Here we use the fact that ${s(r)}\notin \projsupp_\lE(t(r))$ ($\lE$\nbd-adaptedness) to ensure that ${s(r)}\notin\projsupp_\lE(\lambda t(r)+v')$.
  Finally, it follows from
  \begin{gather*}
    \projsupp_\lE\left(\lambda s(r)+\mu s(r)+v'\right)
    \cup\projsupp_\lE\left(\lambda t(r)+\mu s(r)+v'\right)
    \supset \projsupp_\lE\left((\lambda+\mu)t(r)+v'\right)
\end{gather*}
  that $f\relsucc w$ for any esp-order $\succ$.
\end{proof}

\begin{proof}[Proof of \cref{lem:LRSM_congruence_describe_module}]
  (iii) $\Leftrightarrow$ (ii) is given by \cref{lem:LRSM_Church-Rosser_positive_rw}. If $u\dashrightarrow_\lS v$, then ${[u]=[v]}$, so (ii) $\Rightarrow$ (i). To show (i) $\Rightarrow$ (ii), assume that $[u]=[v]$.
  In that case, $u\sim_\lE v+\sum_{i=0}^n\lambda_i(s(r_i)-t(r_i))$ for some scalars $\lambda_i\in\Bbbk$ and relations $r_i\in \lR$. Write $u_j=\sum_{i\leq j}\lambda_i(t(r_i)-s(r_i))$ and $\alpha_j\colon u_j\to u_{j-1}$ the obvious $\lS$\nbd-rewriting step.
  Successively applying the $\alpha_j$'s defines an $\lS$\nbd-rewriting sequence
  \[ u \sim v+u_n\overset{\alpha_n}{\dashrightarrow}_R v+u_{n-1}\overset{\alpha_n}{\dashrightarrow}_R\ldots \overset{\alpha_n}{\dashrightarrow}_Rv.\]
  This concludes.
\end{proof}

\begin{remark}
  \label{rem:LRSM_counterexample_CR_pos_rw}
  Without the adaptedness condition, \cref{lem:LRSM_Church-Rosser_positive_rw} does not hold: for instance, one can consider the \LRS{} $\lP=\{a\to 2a\}$ and the non-positive rewriting step $a+a\to 2a+a$ as a counterexample.
\end{remark}

\subsubsection{Basis from convergence}
\label{subsubsec:LRSM_basis_from_convergence_modulo}

We explain how to extract a basis from convergence modulo.

\begin{definition}
  \label{defn:LRSM_monomial_normal_form}
  Let $\lS=(\cB;\lR,\lE)$ be a \LRSM{}. A \emph{monomial $\lS$\nbd-normal form} is a monomial $b\in\cB$ which is a normal form for $\lS^+$; that is, we have $b\notprojrel_\lE s(r)$ for all $r\in\lR$. We denote $\cB\NF_{\lS}$ the set of monomial $\lS$\nbd-normal forms.
\end{definition}

\begin{remark}
\label{rem:lrs_normal_forms}
  A linear combination of $\lS$\nbd-normal forms is an $\lS$\nbd-normal form, and monomials in the support of an $\lS$\nbd-normal form are $\lS$\nbd-normal forms. The zero vector $0$ is always an $\lS$\nbd-normal form. In other words, $\NF_\lS$ is a $\Bbbk$-module and $\NF_\lS = \langle \cB\NF_\lS\rangle_\Bbbk$.
\end{remark}

\begin{proposition}
  \label{prop:LRSM_basis_from_convergence}
  Let $\lS=(\cB;\lR,\lE)$ be a \LRSM{}.
  If $\lS^+$ is convergent modulo, the canonical linear map
  \[\NF_\lS/\brak{\lE}_\Bbbk\to \brak{\cB}_\Bbbk/\brak{\lR\sqcup\lE}_\Bbbk\]
  is an isomorphism.
  In particular, if $B\subset \cB\NF_\lS$ is such that $[B]_{\lE}$ is a basis for the module $\NF_\lS/\brak{\lE}_\Bbbk$, then $[B]_{\lS}$ is a basis for $\brak{\cB}_\Bbbk/\brak{\lR\sqcup\lE}_\Bbbk$, the module presented by $\lS$.
\end{proposition}

\begin{proof}
  By definition, $\pi_0(\lS^+)$ denotes the $\Bbbk$-module of $\lS^+$-congruence classes. \Cref{lem:LRSM_congruence_describe_module} can be reformulated as stating that the canonical linear map $\pi_0(\lS^+)\to\brak{\cB}_\Bbbk/\brak{\lR\sqcup\lE}_\Bbbk$ is an isomorphism.
  \Cref{prop:pi_zero_from_normal_forms} concludes.
\end{proof}

In the situation of the above proposition, finding a basis reduces to finding a basis of the module $\NF_\lS/\brak{\lE}_\Bbbk$.
Since $\lE$ is monomial-invertible, we may instead view it as a scalar abstract rewriting system $\scl(\lE)$, and if $\scl(\lE)^\top$ is scalar-coherent, then any choice of projective $\lE$-congruence representatives on $\cB\NF_\lS$ defines a basis of $\NF_\lS/\brak{\lE}_\Bbbk$.
This leads to:

\begin{theorem}%
  [\textsc{Basis-From-Convergence Theorem}]
  \label{thm:LRSM_basis_from_convergence_theorem}
  Let $\lS=(\cB;\lR,\lE)$ be a \LRSM{}.
  If $\lS^+$ is convergent modulo and if $\scl(\lE)^\top$ is scalar-coherent on the set $\cB\NF_\lS$ of monomial $\sS^+$\nbd-normal forms,
  then the module $\brak{\cB}_\Bbbk/\brak{\lR\sqcup\lE}_\Bbbk$ presented by $\lS$ is free, and any choice of projective $\lE$\nbd-congruence representatives on $\cB\NF_\lS$ defines a basis.
  \hfill\qed
\end{theorem}

In particular, when $\lE=\emptyset$ we have:

\begin{corollary}
  \label{cor:LRSM_basis_from_convergence}
  Let $(\cB;\lP)$ be a \LRS{}. If $\lP^+$ is convergent, then $\cB \NF_\lP$ is a basis for the module $\brak{\cB}_\Bbbk/\brak{\lP}_\Bbbk$ presented by $\lP$.
  \hfill\qed
\end{corollary}

\subsubsection{Orders}
\label{subsubsec:LRSM_termination_order}

We define the linear analogues of an abstract esp-order (\cref{defn:ARS_esp_order}) and a compatible (abstract) esp-order (\cref{defn:ARSM_compatible_preorder}).

\begin{definition}
  \label{defn:LRSM_esp_order}
  Let $\lS=(\cB;\lR,\lE)$ be a \LRSM{}.
  A \emph{(linear) esp-order} is a strict partial order on $\cB$ which is $\lE$-invariant, that is, such that it is invariant with respect to projective $\lE$\nbd-congruence:
  \[(a'\projrel_\lE a\an a\succ b\an b\projrel_\lE b')\text{ implies }(a'\succ b').\]
\end{definition}

\begin{definition}
  \label{defn:LRSM_compatible_preorder}
  Let $\lS=(\cB;\lR,\lE)$ be a \LRSM{}.
  We say that an esp-order $\succ$ is \emph{compatible with $\lS$} if
  \begin{gather*}
    s(r)\succ b\qquad\text{for all } r\in\lR\an b\in\projsupp_\lE(t(r)).
  \end{gather*}
\end{definition}

Note that thanks to $\lE$-invariance, we can replace $\projsupp_\lE(t(r))$ by $\supp(t(r))$ in the above definition.

Recall from \cref{subsubsec:linear_gray_polygraph} the notion of a linear relation.
An esp-order $\succ$ on $\cB$ induces a linear relation $\succ^+$ on $\brak{\cB}_\Bbbk$, setting $u\succ^+ v$ whenever the following two conditions hold:
\begin{enumerate}[(a)]
  \item $\projsupp_\lE(u)\neq\projsupp_\lE(v)$;
  \item for every $a$ in $\projsupp_\lE(v)\setminus\projsupp_\lE(u)$, there exists ${b}\in\projsupp_\lE(u)\setminus\projsupp_\lE(v)$ such that $b\succ a$.
\end{enumerate}
This definition corresponds to setting $u\succ^+ v$ if and only if $\projsupp_\lE(u)\setsucc\projsupp_\lE(v)$ where $\setsucc$ is the multi-set relation induced by $\succ$, defined on the power set $\cP(\cB)$ as (see e.g.\ \cite[Definition~2.5.3]{BN_TermRewritingAll_1998}):
\[
  M\setsucc N\quad\Leftrightarrow\quad
  M\neq N\an\forall y\in N\setminus M, \exists x\in M\setminus N\text{ such that }x\succ y.
\]
Equivalently, $M\setsucc N$ if and only if one can go from $M$ to $N$ by a sequence of moves consisting in removing an element $b$ and adding elements $a_i$ with $b\succ a_i$. This last interpretation motivates the definition of $\succ^+$, designed precisely such that the following holds:

\begin{lemma}
  \label{lem:LRSM_linear_compatible_implies_abstract_compatible}
  Let $\lS=(\cB;\lR,\lE)$ be a \LRSM{} and $\succ$ an esp-order.
  The linear relation $\succ^+$ is an abstract esp-order on $\brak{\cB}_\Bbbk$.
  Moreover, if $\succ$ is compatible with $\lS$ in the sense of \cref{defn:LRSM_compatible_preorder}, then $\succ^+$ is compatible with $\lS^+$ in the sense of \cref{defn:ARSM_compatible_preorder}.
\end{lemma}

\begin{proof}
  It is shown in \cite[lemma~2.5.4]{BN_TermRewritingAll_1998} that if $\succ$ is a strict partial order, so is $\setsucc$; the first part of the lemma follows. 
  Consider then $r\colon s(r)\to_\lR t(r)$. Strictness (or $\lE$\nbd-adaptedness) implies that
  \[\projsupp_\lE(s(r))\setminus\projsupp_\lE(t(r))=\projsupp_\lE(s(r)),\]
  so that $s(r)\succ^+ t(r)$. The general case follows from the following lemma.
\end{proof}

\begin{lemma}
  \label{lem:LRSM_property_positive_succ}
  Let $\lS=(\cB;\lR,\lE)$ be a \LRSM{} and $\succ$ an esp-order.
  For vectors $u,v\in \langle\cB\rangle_\ring$, we have:
  \[
    u\succ^+v
    \quad\Rightarrow\quad
    \lambda u+w\succ^+\lambda v+w
  \]
  for all $\lambda\in\ring\setminus\{0\}$ and $w\in\langle\cB\rangle_\ring$ such that $\projsupp_\lE(u)\cap\projsupp_\lE(w)=\emptyset$.
\end{lemma}

\begin{proof}
  Since $\projsupp_\lE(\lambda u+w) = \projsupp_\lE(u)\sqcup\projsupp_\lE(w)$, we have:
  \begin{gather*}
    \projsupp_\lE(\lambda u+w)\setminus\projsupp_\lE(\lambda v+w)\supset\projsupp_\lE(u)\setminus\projsupp_\lE(v),\\
    \projsupp_\lE(\lambda v+w)\setminus\projsupp_\lE(\lambda u+w)\subset\projsupp_\lE(v)\setminus\projsupp_\lE(u).\qedhere
  \end{gather*}
\end{proof}

We now relate to termination:

\begin{lemma}
  \label{lem:LRSM_well-founded_implies_terminates}
  Let $\lS=(\cB;\lR,\lE)$ be a \LRSM{} and $\succ$ a compatible esp-order.
  If $\succ$ is well-founded, then $\succ^+$ is well-founded.
  In that case, $\lS^+$ terminates modulo $\lE$.
\end{lemma}

\begin{proof}
  It is shown in \cite[theorem~2.5.5]{BN_TermRewritingAll_1998} that if $\succ$ is well-founded, then $\setsucc$ is well-founded. In that case, \cref{lem:LRSM_linear_compatible_implies_abstract_compatible} implies that $\lS^+$ terminates.
\end{proof}

We conclude the subsection with the following notion, which appears in the statement of \cref{lem:LRSM_Church-Rosser_positive_rw}:

\begin{definition}
  \label{defn:LRSM_relative_relation}
  Let $\succ$ be a relation on a set $\cB$. The \emph{relative relation $\relsucc$ induced by $\succ$} is the following relation on the power set $\cP(\cB)$:
  \begin{gather*}
    M\relsucc N\quad\Leftrightarrow\quad \forall b\in\cB, (b\succ M)\Rightarrow (b\succ N).
  \end{gather*}
  If $\lS=(\cB;\lR,\lE)$ is a \LRSM{} and $\succ$ an esp-order, we set:
  \[u\relsucc v\quad\Leftrightarrow\quad \projsupp_\lE(u)\relsucc\projsupp_\lE(v).\]
\end{definition}

Note that equivalently, $M\relsucc N$ if and only if $\forall L\in\cP(\cB)$, $(L\succ M)$ implies $(L\succ N)$.

\subsubsection{Tamed linear Newman's lemma}
\label{subsubsec:LRSM_tamed_newmann_lemma}

The linear structure of a \LRSM{} $\lS=(\cB;\lR,\lE)$ induces canonical types of $\lS^+$\nbd-local triples. Let $[f,e,g]$ be a $\lS^+$\nbd-local triple such that $f$ and $g$ is of type $r_1$ and $r_2$, respectively.
Then:
\begin{itemize}
  \item if $s(r_1)\notprojrel_\lE s(r_2)$, one says that $[f,e,g]$ is \emph{additive},
  \item if $s(r_1)\projrel_\lE s(r_2)$, one says that $[f,e,g]$ is \emph{intersecting}.
\end{itemize}
Recall the notion of branchwise $\lE$-congruence (\cref{subsubsec:ARSM_E_congruence}).
A $\lS^+$\nbd-local triple $[f,e,g]$ is additive if, up to branchwise $\lE$-congruence, it has the following form:
\begin{gather*}
  [f,e,g]=[\lambda_1 r_1+\lambda_2s(r_2)+v,\id,\lambda_1s(r_1)+\lambda_2 r_2+v],
\end{gather*}
An additive branching has a canonical $\lS$\nbd-confluence:
\[\begin{tikzcd}[ampersand replacement=\&,cramped,row sep=large,column sep=small]
  \& {\lambda_1t(r_1)+\lambda_2 s(r_2)+v} \\
  {\lambda_1s(r_1)+\lambda_2 s(r_2)+v} \&\& {\lambda_1t(r_1)+\lambda_2t(s_2)+v} \\
  \& {\lambda_1 s(r_1)+\lambda_2t(r_2) + v}
  \arrow["{\lambda_1t(r_1)+\lambda_2r_2+v}",bend left=15pt, dashed, from=1-2, to=2-3]
  \arrow["{\lambda r_1+\lambda_2 s(r_2)+v}",bend left=15pt, from=2-1, to=1-2]
  \arrow["{\lambda_1 s(r_1)+\lambda_2r_2+v}"',bend right=15pt, from=2-1, to=3-2]
  \arrow["{\lambda_1r_1+\lambda_2t(r_2)+v}"',bend right=15pt, dashed, from=3-2, to=2-3]
\end{tikzcd}\]
On the other hand, $[f,e,g]$ is intersecting if, up to branchwise $\lE$-congruence, it has the following form:
\[
  [f,e,g]=[\lambda_1 r_1+w,e_r+w,\lambda_2 r_2+w],
\]
where $e_r\colon \lambda_1 s(r_1)\sim_\lE \lambda_2 s(r_2)$.
We say that $[f,e,g]$ is \emph{monomial} if it is intersecting and if, up to branchwise $\lE$-congruence, it has the form $\lambda_1=1$ and $w=0$ with the notations above.
That is, $[f,e,g]$ is monomial if, up to branchwise $\lE$-congruence, it has the following form:
\[
  [f,e,g]=[r_1,e_r,\lambda r_2].
\]
where $e_r\colon s(r_1)\sim_\lE \lambda s(r_2)$.

Note that both branches of a monomial branching are monomial rewriting steps.
By definition, every intersecting branching is of the form $\lambda_1[f,e,g]+w$ for $[f,e,g]$ a monomial branching. Still by definition, additive, intersecting and monomial branchings are positive branchings.

\medbreak

We end the subsection with a linear analogue of the tamed Newman's lemma (\cref{lem:ARSM_newmann_lemma}). Given a \LRSM{} $\lS=(\cB;\lR,\lE)$ and an esp-order $\succ$, we say ``$\succ$\nbd-tameness'' to refer to $\succ^+$\nbd-tame\-ness.

\begin{theorem}%
  [\textsc{Tamed Linear Newman's Lemma}]
  \label{thm:LRSM_tamed_newmann_lemma}
  Let $\lS=(\cB;\lR,\lE)$ be a \LRSM{} and $\succ$ a compatible esp-order.
  If $\succ$ is well-founded and every monomial $\lS$\nbd-branching is $\succ$\nbd-tamely $\lS$\nbd-congruent, then $\lS^+$ is convergent modulo $\lE$.
\end{theorem}

Our proof can be understood as a generalization of \cite[Theorem~4.2.1, part~(ii)]{GHM_ConvergentPresentationsPolygraphic_2019}, working with tamed $\lS$\nbd-congruence instead of  $\lS^+$\nbd-confluence.
The hypothesis that $\succ$ is well-founded is necessary; see \cite[Remark~4.2.4]{GHM_ConvergentPresentationsPolygraphic_2019} for a counterexample.

\begin{proof}
  Every intersecting branching is $\succ$-tamely congruent. Indeed, following the discussion above, an intersecting branching is of the form $(f,g)+w$ where $(f,g)$ is a monomial branching and $w$ is a vector.
  By hypothesis, the branching $(f,g)$ admits a $\succ$-tamed $\lS$-congruence $(f',g')$, and it follows from \cref{lem:LRSM_property_positive_succ} that the $\lS$-congruence $(f',g')+w$ is $\succ$-tamed by the source of $(f,g)+w$.

  Consider then an additive $\lS^+$\nbd-local triple as pictured above.
  We wish to show that its canonical $\lS$\nbd-confluence is tamed, that is:
  \begin{gather*}
    \lambda_1s(r_1)+\lambda_2s(r_2)+w\succ^+\lambda_1t(r_1)+\lambda_2t(s_2)+w.
  \end{gather*}
  This happens precisely if $s(r_2)\notin \projsupp_\lE(t(r_1))$ or if $s(r_1)\notin \projsupp_\lE(t(r_2))$. (These two cases correspond respectively to the upper branch or lower branch being positive.)
  Assuming otherwise, the compatibility of $\succ$ implies that $s(r_2)\succ s(r_1)$ and $s(r_1)\succ s(r_2)$. This contradicts the assumption that $\succ$ is well-founded.
  We conclude that every local $\lS^+$\nbd-branching is $\succ$\nbd-tamely $\lS$\nbd-congruent.

  Note then that:

  \begin{lemma}
    A branching is $\succ$\nbd-tamely $\lS$\nbd-congruent if and only if it is $\succ$\nbd-tamely $\lS^+$\nbd-congruent.
  \end{lemma}

  \begin{proof}
    This follows from \cref{lem:LRSM_Church-Rosser_positive_rw} and the fact that for $M,N_1,N_2,L\in\cP(\cB)$, if $M\relsucc N_1\cup N_2$ and $N_1\relsucc L$, then $M\relsucc N_1\cup N_2$.
  \end{proof}

  It follows that every local $\lS^+$\nbd-branching is $\succ$\nbd-tamely $\lS^+$\nbd-congruent. We can now apply the tamed (abstract) Newman's lemma (\cref{lem:ARSM_tamed_newmann_lemma}) and conclude.
\end{proof}

\subsection{Higher linear rewriting modulo}
\label{subsec:higher_linear_rewriting_modulo}

\subsubsection{Higher linear rewriting system modulo}
\label{subsubsec:HLRSM_definition}

\begin{definition}
  \label{defn:HLRSM}
  A \emph{higher linear rewriting system modulo} (\HLRSM{}) $\sS=(\sR,\sE)$ is the data of a left-monomial 3\nbd-sesqui\-poly\-graph $\sR$ and a
  monomial-invertible
  3\nbd-sesqui\-poly\-graph $\sE$, such that $\sR_{\leq 2}=\sE_{\leq 2}$.
\end{definition}

If we include interchangers in the modulo data, we have (recall the notations $G$ and $\mu$ from \cref{not:sesquicat_notation}):

\begin{definition}
A $(G,\mu)$\nbd-graded linear Gray rewriting system modulo is a \HLRSM{} $\sS=(\sR,\sE)$ such that $\sE$ is a $(G,\mu)$\nbd-graded linear Gray polygraph.
\end{definition}

As in the abstract higher case, a \HLRSM{} $\sS=(\sR,\sE)$ defines a category of linear 1\nbd-poly\-graphs $\lS(\square)=(\cB(\square),\lR(\square),\lE(\square))$, setting
\begin{gather*}
  \cB(\square)\coloneqq \sR_{2}^*(\square)=\sE_{2}^*(\square),
  \\
  \lR(\square)\coloneqq \Cont(\sR)(\square)
  \quad\an\quad
  \lE(\square)\coloneqq \Cont(\sE)(\square),
\end{gather*}
for each 1-sphere $\square$ in $\sR_{1}^*=\sE_{1}^*$.
Each context $\Gamma$ on $\square$ defines a morphism of linear 1\nbd-poly\-graphs:
\[
  \Gamma\colon \lS(\square)\to\lS(s_1(\Gamma),t_1(\Gamma)).
\]
Contrary to the abstract higher case however, \emph{$\sS(\square)$ needs not be a \LRSM{}}, as we did not impose the adaptedness condition appearing in \cref{defn:linear-rewriting-system-modulo}.

As for the abstract higher case, we define a notion of sub-systems, imposing now the adaptedness condition:

\begin{definition}
  \label{defn:HLRSM_subsystem}
  Let $\sS=(\sR,\sE)$ be a \HLRSM{}.
  A \emph{linear sub-system} $\lT$ is the data of a family of sub-sets $\lT(\square)\subset\sR_3^*(\square)$ for each 1-sphere $\square$ in $\sR_{1}^*$.
  Moreover, we assume the \emph{adaptedness condition}:
  \[
    s(r)\notin\projsupp_\sE(t(r)) \quad\forall r\in\lT.
  \]
  We say that $\sS$ is \emph{adapted} if $\sS$ is a linear sub-system of itself.
\end{definition}

As in the higher abstract case, we think of $\lT$ as defining a family of sub-\LRSM{} $\lT(\square)\subset\sS(\square)$ modulo $\sE(\square)$, and write $\lT\subset\sS$ to emphasize that point.
We say that $\lT$ is \emph{confluent} if every sub-\LRSM{} $\lT(\square)$, and similarly for other notions.
Note that adaptedness \emph{must} be defined via a linear sub-system, since in general, it depends on context: if $r$ is adapted, $\Gamma[r]$ needs not be.

A $\lT$\nbd-rewriting step has the following general form, where $r\in\sR$:
\[
  \alpha\colon u\;\;\sim_\lE \;\;\lambda \Gamma[s(r)]+v\;\;\overset{}{\dashrightarrow}_\lT \;\;\lambda \Gamma[t(r)]+v\;\;\sim_\lE \;\;w.
\]
We say that $\alpha$ is \emph{of type $r$}.\footnote{This is a slight abuse of notation: as belonging to a \LRSM{}, $\alpha$ is of type $\Gamma[r]$ according to \cref{subsubsec:LRSM_positive_rw}.}
As before, $\alpha$ is positive if $\Gamma[s(r)]\notin\projsupp_\sE(v)$.
Note that \emph{contextualization needs not preserve positivity}: we may have
\[\Gamma'[\Gamma[s(r)]]\in\projsupp_\sE(\Gamma'[v]),\quad\text{ even if }\quad\Gamma[s(r)]\notin\projsupp_\sE(v).\]
However, it is true that if $\alpha$ is monomial (that is, if $v=0$) then $\Gamma[\alpha]$ is monomial.

\subsubsection{Strong compatibility}
\label{subsubsec:HLRSM_strong_compatibility}

Fix $\sS=(\sR,\sE)$ a \HLRSM{} and $\lT$ a linear sub-system.
Mimicking previous definitions,
we say that an esp-order $\succ$ is \emph{compatible with $\lT$} if each \LRSM{} $\lT(\square)$ is compatible in the sense of \cref{defn:LRSM_compatible_preorder}, that is, if
\[
  s(r)\succ b,\quad\text{for all }r\in\lT\an b\in\projsupp_\sE(t(r)).
\]
As we shall see in the next subsection, compatibility is \emph{not} sufficient for higher rewriting, and a stronger condition is required:

\begin{definition}
  \label{defn:strong_compatibility}
  Let $\sS=(\sR,\sE)$ be a \HLRSM{} and $\lT$ a linear sub-system.
  An esp-order is said to be \emph{strongly compatible with $\lT$} if:
  \[
    \Gamma[s(r)]\succ b,\quad\text{for all }r\in\sR,\text{context }\Gamma\text{ such that }\Gamma[r]\in\lT,\an b\in\Gamma\big[\projsupp_\sE(t(r))\big].
  \]
\end{definition}

Strong compatibility is indeed stronger than compatibility:

\begin{lemma}
  \label{lem:strong_compatibility_implies_compatibility}
  If $\succ$ is strongly compatible with $\lT$, then $\succ$ is compatible with $\lT$.
\end{lemma}

\begin{proof}
  It suffices to choose the trivial context in the strong compatibility condition.
\end{proof}

However, the converse of \cref{lem:strong_compatibility_implies_compatibility} does not hold!
One can only say that since $\Gamma[r]\in\lT$, compatibility implies that $\Gamma[s(r)]=s(\Gamma[r])\succ b$ for all $b\in\projsupp_\sE(t(\Gamma[r]))=\projsupp_\sE(\Gamma[t(r)])$.
This is \emph{not} sufficient to conclude, as the inclusion
\[\projsupp_\sE(\Gamma[t(r)])\subset\Gamma\big[\projsupp_\sE(t(r))\big]\]
needs not be an equality!
Consider more generically a vector $w$ and the inclusion
\[\projsupp_\sE(\Gamma[w])\subset\Gamma\big[\projsupp_\sE(w)\big].\]
Recall that contextualization needs not act freely (see \cref{subsubsec:intro_independent_contextualized_branchings}): there may exist monomials $b_1,b_2$ such that $b_1\notprojrel_\lE b_2$ but $\Gamma[b_1]\projrel_\lE\Gamma[b_2]$. In particular, if $w=\lambda_1b_1+\lambda_2b_2$ and $\lambda_1,\lambda_2$ have well-chosen scalars, we have $\Gamma[w]=0$ and thus $\Gamma[b_1]\notin\projsupp_\sE(\Gamma[w])$.

\begin{remark}
  Super $\glt$-foams give an explicit example where compatibility does not imply strong compatibility.
  Consider $\Bbbk$ as in \cref{defn:ring_R_bil} with $X=Z=1$ and $Y=-1$; this is the super case.
  We let $\foamS=(\foamR,\foamE)$ and $\foamlT$ be as in \cref{subsec:convergent_pres_gltfoam}.
  Let $\widetilde{\foamlT}\subset\foamlT$ be the linear sub-system of $\foamlT$ consisting only in rewriting steps of type $\nc$; that is, $\widetilde{\foamlT}$ consists in rewriting steps of type $\nc$ for which the two strands are distinct.
  Let $\widetilde{\succ}$ be the esp-order on $\foamS_2^*$ which compares the number of \emph{open shadings}---shadings that are not enclosed by a strand.
  The esp-order is $\widetilde{\succ}$ compatible with $\widetilde{\foamlT}$.
  Compatibility is clear for $\Gamma[\nc]\in\widetilde{\foamlT}$ such that both strands are not closed, as the number of open shadings strictly decreases.
  On the other hand, if at least one of the strand is closed, then the target of $\Gamma[\nc]$ is zero:
  \begin{gather*}
    \def\scl{.6}
    \tikzpic{
      \fill[fill=colour_diag1,opacity=.25] (0,0) to[out=-90,in=0] ++(-.5,-.5) to[out=180,in=-90] ++(-.5,.5) to ++(0,1.8) to[out=90,in=180] ++(.5,.5) to[out=0,in=90] ++(.5,-.5);
      \draw[diag1=1.4pt] (0,0) to[out=-90,in=0] ++(-.5,-.5) to[out=180,in=-90] ++(-.5,.5) to ++(0,1.8) to[out=90,in=180] ++(.5,.5) to[out=0,in=90] ++(.5,-.5);
      \draw[diag1=1.4pt] (0,0) to (0,1.8);
      \fill[fill=colour_diag1,opacity=.25] (1,0) rectangle (1.5,1.8);
      \draw[diag1=1.4pt] (1,0) to (1,1.8);
    }[scale=\scl*.7]
    \qquad\to_{\widetilde{\foamlT}}\qquad
    \tikzpic{
      \begin{scope}[shift={(.5,1.7)}]
        \node[fdot1=1.7pt] at (0,0) {};
      \end{scope}
      \draw[diag1=1.4pt] (1,0) to[out=90,in=0] (.5,.75) to[out=180,in=90] (0,0) to[out=-90,in=0] ++(-.5,-.5) to[out=180,in=-90] ++(-.5,.5) to ++(0,1.8) to[out=90,in=180] ++(.5,.5) to[out=0,in=90] ++(.5,-.5) to[out=-90,in=180] (.5,.25+.8) to[out=0,in=-90] (1,1+.8);
      \fill[fill=colour_diag1,opacity=.25] (1,0) to[out=90,in=0] (.5,.75) to[out=180,in=90] (0,0) to[out=-90,in=0] ++(-.5,-.5) to[out=180,in=-90] ++(-.5,.5) to ++(0,1.8) to[out=90,in=180] ++(.5,.5) to[out=0,in=90] ++(.5,-.5) to[out=-90,in=180] (.5,.25+.8) to[out=0,in=-90] (1,1+.8);
    }[scale=\scl*.7]
    +
    \tikzpic{
      \begin{scope}[shift={(.5,.1)}]
        \node[fdot1=1.7pt] at (0,0) {};
      \end{scope}
      \draw[diag1=1.4pt] (1,0) to[out=90,in=0] (.5,.75) to[out=180,in=90] (0,0) to[out=-90,in=0] ++(-.5,-.5) to[out=180,in=-90] ++(-.5,.5) to ++(0,1.8) to[out=90,in=180] ++(.5,.5) to[out=0,in=90] ++(.5,-.5) to[out=-90,in=180] (.5,.25+.8) to[out=0,in=-90] (1,1+.8);
      \fill[fill=colour_diag1,opacity=.25] (1,0) to[out=90,in=0] (.5,.75) to[out=180,in=90] (0,0) to[out=-90,in=0] ++(-.5,-.5) to[out=180,in=-90] ++(-.5,.5) to ++(0,1.8) to[out=90,in=180] ++(.5,.5) to[out=0,in=90] ++(.5,-.5) to[out=-90,in=180] (.5,.25+.8) to[out=0,in=-90] (1,1+.8);
    }[scale=\scl*.7]
    \; \sim_\foamE \; 0.
  \end{gather*}
  Hence $\Gamma[\nc]$ is trivially compatible with $\succ$.
  However, the same situation prevents the esp-order $\widetilde{\succ}$ from being strongly compatible with $\widetilde{\foamlT}$, as the number of open shadings is the same for $\Gamma[s(\nc)]$ and $b$, for all $b\in\Gamma\big[\projsupp_\foamE(t(\nc))\big]$.
\end{remark}

\subsubsection{Contextualization}
\label{subsubsec:HLRSM_contextualization}

We now study how confluence, tamed congruence and contextualization relates.

\begin{lemma}
  \label{lem:HLRSM-contextualization-sequence}
  Let $\sS=(\sR,\sE)$ be a \HLRSM{} and let $\lT$ be a linear sub-system and $\succ$ strongly compatible esp-order.
  Consider a composition $g\starop_2 f$, where $f$ is a monomial $\sS$\nbd-rewriting step and $g$ is a positive $\sS$\nbd-rewriting sequence.
  If $\Gamma$ is a context such that both $\Gamma[f]$ and $\Gamma[g]$ are in $\lT$,
  then $\Gamma[g]$ is $\succ$\nbd-tamed by $\Gamma[s(f)]$.
\end{lemma}

\begin{proof}
  Decomposing $g$ as a composition of positive $\sS$\nbd-rewriting steps, the situation is:
  \[\begin{tikzcd}[ampersand replacement=\&,cramped]
    s(f) \& {v_0} \& {v_1} \& \ldots \& {v_n}\\
    \Gamma[s(f)] \& {\Gamma[v_0]} \& {\Gamma[v_1]} \& \ldots \& {\Gamma[v_n]}
    \arrow["{f}","\sS"{subscript}, from=1-1, to=1-2]
    \arrow["{g_0}","\sS"{subscript}, from=1-2, to=1-3]
    \arrow["{g_1}","\sS"{subscript}, from=1-3, to=1-4]
    \arrow["{g_{n-1}}","\sS"{subscript}, from=1-4, to=1-5]
    \arrow["{\Gamma[f]}","\lT"{subscript}, from=2-1, to=2-2]
    \arrow["{\Gamma[g_0]}","\lT"{subscript},dashed, from=2-2, to=2-3]
    \arrow["{\Gamma[g_1]}","\lT"{subscript},dashed, from=2-3, to=2-4]
    \arrow["{\Gamma[g_{n-1}]}","\lT"{subscript},dashed, from=2-4, to=2-5]
    \arrow["\Gamma",start anchor=east,end anchor=east,bend left=60pt,from=1-5,to=2-5]
  \end{tikzcd}\]
  We wish to show that $\Gamma[s(f)]\succ\projsupp_\sE\left(\Gamma[v_i]\right)$ for all $0\leq i\leq n$. Rather, given the inclusion $\projsupp_\sE(\Gamma[w])\subset \Gamma\left[\projsupp_\sE(w)\right]$ which holds for generic vector $w$, it suffices to show
  \begin{equation}
    \label{eq:linear_contextualization_lemma}
    \Gamma[s(f)]\succ\Gamma\left[\projsupp_\sE(v_i)\right]
  \end{equation}
  for all $0\leq i\leq n$.

  When $i=0$, we have $v_0=t(f)$. By hypothesis $\Gamma[f]\in\lT$, so that \eqref{eq:linear_contextualization_lemma} follows from strong compatibility.
  It remains to show that
  \[
    \Gamma\left[\projsupp_\sE(v_i)\right]\relsucc\Gamma\left[\projsupp_\sE(v_{i+1})\right].
  \]
  (Recall the definition of relative relation from \cref{defn:LRSM_relative_relation}.)
  Write $g_i=\lambda r+v$, its canonical decomposition as a positive $\sS$\nbd-rewriting step. On the one hand:
  \[
    \Gamma\left[\projsupp_\sE(v_{i+1})\right]
    = \Gamma\left[\projsupp_\sE(t(r)+v)\right] 
    \subset \Gamma\left[\projsupp_\sE(t(r))\right]\cup\Gamma\left[\projsupp_\sE(v)\right].
  \]
  On the other hand, thanks to the positivity of $g_i$, we have 
  $\projsupp_\sE(s(r)+v) = \{s(r)\}\sqcup\projsupp_\sE(v),$
  and hence
  \[
    \Gamma\left[\projsupp_\sE(v_i)\right]
    = \Gamma\left[\projsupp_\sE(s(r)+v)\right]
    = \{\Gamma[s(r)]\}\cup\Gamma\left[\projsupp_\sE(v)\right].
  \]
  Because $\Gamma[g_i]\in\lT$, strong compatibility implies that $\Gamma[s(r)]\succ \Gamma\left[\projsupp_\sE(t(r))\right]$.
  This concludes.
\end{proof}

\begin{corollary}%
  [\textsc{Contextualization Lemma}]
  \label{cor:HLRSM_contextualization_branching_lemma}
  Let $\sS=(\sR,\sE)$ be a \HLRSM{} and let $\lT$ be a linear sub-system and $\succ$ strongly compatible esp-order.
  Let $(f,g)$ be a monomial $\sS$\nbd-branching which admits a positive $\sS$\nbd-confluence $(f',g')$.
  If $\Gamma$ is a context such that both $\Gamma[f,g]$ and $\Gamma[f',g']$ belong to $\lT$, then $\Gamma[f',g']$ is a $\succ$\nbd-tamed $\lT$\nbd-congruence for $\Gamma[f,g]$.
  \hfill\qed
\end{corollary}

In general, $\Gamma[f',g']$ needs not be positive; see \cref{footnote:gap_critical_confluence} for an example and comparison with the literature.
In the context-agnostic setting, the \HLRSMcontextualizationlemma{} reduces to the following:

\begin{corollary} 
  \label{cor:HLRSM_contextualization_branching_lemma-agnostic}
  Let $\sS=(\sR,\sE)$ be an adapted \HLRSM{} and $\succ$ a strongly compatible esp-order.
  Let $(f,g)$ be a monomial $\sS$\nbd-branching and $\Gamma$ a context.
  If $(f,g)$ is positively $\sS$\nbd-confluent, then $\Gamma[f,g]$ is $\succ$\nbd-tamely $\sS$\nbd-congruent.
\end{corollary}

\subsubsection{Independent branching}
\label{subsubsec:HLRSM_independent_branching}

As we have seen, for \HRSM{} the statement that independent branchings are confluent is tautological (\cref{lem:HRSM_independent_branching_agnostic}).
As we now explain, the linear case is more subtle.

Fix $\sP$ a \HLRS{} and $(f,g)$ an independent $\sP$\nbd-branching. Without loss of generality, this means that $f=\ssf\starop_1 s(\ssg)$ and $g=s(\ssf)\starop_1\ssg$, for some monomial $\sP$\nbd-rewriting steps $\ssf$ and $\ssg$.
The $\sP$\nbd-confluence $(g',f')=(t(\ssf)\starop_1 \ssg,\ssf\starop_1t(\ssg))$ defines a canonical $\sP$\nbd-confluence for $(f,g)$.
Write:
\begin{gather*}
  \ssf\colon s(\ssf)\to \lambda_1x_1+\ldots+\lambda_mx_m
  \quad\an\quad
  \ssg\colon s(\ssg)\to \mu_1y_1+\ldots+\mu_ny_n,
\end{gather*}
where the $x_i$'s (resp.\ $y_j$'s) are monomials with pairwise distinct $\sE$\nbd-projective classes.
Then the $\sP$\nbd-confluence $(g',f')$ can be explicitly decomposed as a $\sP$\nbd-confluence, setting $f'=f_n'\starop_2\ldots\starop_2 f_1'$ and $g'=g_m'\starop_2\ldots\starop_2 g_1'$, where:
\begin{align*}
  f'_j&=\sum_{k<j} \mu_k(t(\ssf)\starop_1 y_k)+\mu_j (\ssf \starop_1 y_j)+\sum_{j<k} \mu_k(s(\ssf)\starop_1 y_k)
  \\
  \an\;g'_i&=\sum_{k<i} \lambda_k(x_k\starop_1t(\ssg))+\lambda_i (x_i\starop_1\ssg)+\sum_{i<k} \lambda_k(x_k\starop_1s(\ssg)).
\end{align*}

\begin{lemma}[independent branchings, context-dependent case]
  \label{lem:HLRSM_independent_branching_base_case}
  Let $\sS=(\sR,\sE)$ be a \HLRSM{} and let $\lT$ be a linear sub-system and $\succ$ strongly compatible esp-order.
  Let $(f,g)$ be an independent $\lT$\nbd-branching.
  If for each $1\leq j\leq n$ (resp.\ $1\leq i\leq m$) and with the notations above, the $\sS$\nbd-rewriting step $\ssf \starop_1 y_j$ (resp.\ $x_i\starop_1 \ssg$) is in $\lT$, then the canonical $\lT$\nbd-confluence $(g',f')$ is $\succ$\nbd-tamed by $s(f)=s(g)$.
\end{lemma}

\begin{proof}
  For each $1\leq j\leq n$, the $\sE$-projective supports of $s(f_j')$ and $t(f_j')$ are contained in the following set:
  \begin{gather*}
    M\coloneqq
    \bigcup_{1\leq j\leq n}
    \projsupp_\sE\big(t(\ssf)\starop_1 y_j\big)
    \cup\bigcup_{1\leq j\leq n}
    \projsupp\big(s(\ssf)\starop_1 y_j\big).
  \end{gather*}
  It follows from strong compatibility that $s(\ssf)\starop_1s(\ssg)\succ s(\ssf)\starop_1y_j$ (since $g=s(\ssf)\starop_1\ssg$ belongs to $\lT$)
  and that
  $s(\ssf)\starop_1y_j\succ t(\ssf)\starop_1y_j$ (since $\ssf \starop_1 y_i$ belongs to $\lT$).
  By transitivity of the esp-order $\succ$, we have $s(\ssf)\starop_1s(\ssg)\succ t(\ssf)\starop_1y_j$.
  Hence:
  \[s(\ssf)\starop_1 s(\ssg)\succ M\]
  In other words, $f'$ is $\succ$\nbd-tamed by $s(\ssf)\starop_1 s(\ssg)$.

  An analogous argument shows that $g'$ is $\succ$\nbd-tamed by $s(\ssf)\starop_1 s(\ssg)$. This concludes.
\end{proof}

Note the necessity of strong compatibility in the proof of \cref{lem:HLRSM_independent_branching_base_case}: $s(\ssf)\starop_1s(\ssg)\succ s(\ssf)\starop_1y_j$ does \emph{not} follow from compatibility alone, as while $y_j$ is in the $\sE$-projective support of $t(\ssg)$, it may be that $s(\ssf)\starop_1y_j$ is \emph{not} in the support of $s(\ssf)\starop_1t(\ssg)$.

The proof of \cref{lem:HLRSM_independent_branching_base_case} resembles the proof of Theorem 4.2.1, part (iii) in \cite{GHM_ConvergentPresentationsPolygraphic_2019}, with tameness appearing implicitly.
In our terminology, their setting (associative algebra) restrict to context-agnostic systems, and as contextualization acts freely, there is no distinction between compatibility and strong compatibility.

The following is a direct corollary of \cref{lem:HLRSM_independent_branching_base_case}:

\begin{lemma}[independent branchings, context-agnostic case]
  \label{lem:HLRSM_congruence_independent_branching}
  Let $\sS=(\sR,\sE)$ be an adapted \HLRSM{} and $\succ$ a strongly compatible esp-order.
  Every independent $\sS$\nbd-branching is $\succ$\nbd-tamely $\sS$\nbd-congruent.
  \hfill\qed
\end{lemma}

In the context-dependent case, one cannot rely on such a general statement. However:

\begin{lemma}
  \label{lem:HLRSM_independent-branching-B-confluence}
  Let $\sS=(\sR,\sE)$ be a \HLRSM{} and let $\lT$ be a linear sub-system and $\succ$ strongly compatible esp-order.
  Assume given another linear sub-system $\lB\subset\lT$, such that $\lB^+$ is convergent. If every $(\sS\setminus\lT)$\nbd-rewriting step is $\lB$\nbd-congruent, then every independent $\lT$\nbd-branching is $\succ$\nbd-tamely $\lT$\nbd-congruent.
\end{lemma}

This will be the situation of graded $\glt$\nbd-foams. 

\begin{proof}
  We adapt the proof of \cref{lem:HLRSM_independent_branching_base_case}, and borrow its notations.
  Assume, up to reordering the $y_j$'s, that $\ssf\starop_1 y_j$ is in $\lT$ for $1\leq j\leq k$ and that $\ssf\starop_1 y_j$ is \emph{not} in $\lT$ for $k<j\leq m$. By the same argument as in the proof of \cref{lem:HLRSM_independent_branching_base_case}, $f_k'\starop_2\ldots\starop_2 f_1'$ is a $\lT$\nbd-rewriting sequence $\succ$\nbd-tamed by $s(\ssf)\starop_1s(\ssg)$. On the other hand, by the hypothesis of the lemma, $f_m'\starop_2\ldots\starop_2 f_{k+1}'$ is $\lB$\nbd-congruent.
  Arguing similarly for $g'$, we get a $\lT$\nbd-congruence for $(f,g)$ given as a $\lB$\nbd-congruence ``sandwiched'' between two $\succ$\nbd-tamed $\lT$\nbd-congruence.
  Since $\lB^+$ is convergent, we can replace the $\lB$\nbd-congruence by a $\lB^+$\nbd-confluence. Compatibility of $\succ$ with $\lB$ and transitivity of $\succ^+$ conclude.
\end{proof}

\subsubsection{Independent rewriting}
\label{subsubsec:HLRSM_independent_rewriting}

We consider the linear analogue of abstract independent rewriting defined in \cref{subsubsec:HRSM_independent_rewriting}. Let us use the same notations with $\ssf$, $\ssg$ and $\ssh$ monomial
 rewriting steps, which decompose as:
\begin{align*}
  \ssf&\colon\phi\to_\sR\lambda_1x_1+\ldots+\lambda_l x_l,\\
  \ssg&\colon\phi\to_\sR\mu_1y_1+\ldots+\mu_my_m,\\
  \an\;\ssh&\colon\psi\to_\sR\nu_1z_1+\ldots+\nu_nz_n,
\end{align*}
where $x_i$'s (resp.\ $y_j$'s and $z_k's$) are monomials belonging to pairwise distinct $\sE$\nbd-projective classes.
In these decompositions, the independent rewriting of $(f,g)$ into $(f',g')$ via the triple $(h,h_1,h_2)$ is pictured as follows:
\[\begin{tikzcd}[column sep=small,row sep=small]
  &
  \sum_i\lambda_i(x_i\starop_1\psi)
  &&
  \sum_{i,k}\lambda_i\nu_k(x_i\starop_1z_k)
  \\
  \phi\starop_1\psi
  &&
  \sum_k\nu_k(\phi\starop_1z_k)
  \\
  &
  \sum_j\mu_j(y_j\starop_1\psi)
  &&
  \sum_{j,k}\mu_j\nu_k(y_j\starop_1z_k)
  \arrow[from=2-1,to=1-2,"f=\ssf\starop_1\psi"]
  \arrow[from=2-1,to=3-2,"g=\ssg\starop_1\psi"']
  \arrow[from=2-1,to=2-3,"h=\phi\starop_1\ssh"]
  \arrow[from=2-3,to=1-4,"f'=\sum_k\nu_k(\ssf\starop z_k)"{pos=.4}]
  \arrow[from=2-3,to=3-4,"g'=\sum_k\nu_k(\ssg\starop z_k)"'{pos=.4}]
  \arrow[from=1-2,to=1-4,"h_1=\sum_i\lambda_i(x_i\starop_1\ssh)"]
  \arrow[from=3-2,to=3-4,"h_2=\sum_j\lambda_j(y_j\starop_1\ssh)"']
\end{tikzcd}\]

\begin{lemma}
  [\textsc{Independent Rewriting Lemma}]
  \label{lem:HLRSM_independent_rewriting_context_dependent}
  Let $\sS=(\sR,\sE)$ be a \HLRSM{} and let $\lT$ be a linear sub-system and $\succ$ strongly compatible esp-order.
  Let $\ssf$, $\ssg$ and $\ssh$ monomial $\sR$\nbd-rewriting steps as above, such that:
  \begin{enumerate}[(i)]
    \item $f=\ssf\starop_1\psi$, $g=\ssg\starop_1\psi$ and $h=\phi\starop_1\ssh$ are in $\lT$, and for each $1\leq i\leq l$ (resp.\ $1\leq j\leq m$, resp.\ $1\leq k\leq n$), the rewriting step $x_i\starop_1\ssh$ (resp.\ $y_j\starop_1\ssh$) is in $\lT$,
    \item for each $1\leq k\leq n$, the branching $(\ssf\starop_1z_k,\ssg\starop_1z_k)$ is $\succ$\nbd-tamely $\lT$\nbd-congruent,
  \end{enumerate}
  then $(f,g)$ is $\succ$\nbd-tamely $\lT$\nbd-congruent.
\end{lemma}

\begin{proof}
  Since $h=\phi\starop_1\ssh$ is in $\lT$, strong compatibility implies that $\phi\starop_1\psi\succ\phi\starop_1z_k$ (for all suitable $k$).
  Let $h_k'$ be a $\lT$\nbd-congruence for $(\ssf\starop_1z_k,\ssg\starop_1z_k)$ $\succ$-tamed by its source $\phi\starop_1z_k$, and hence by $\phi\starop_1\psi$ by transitivity of $\succ$.
  It follows that the $\lT$\nbd-congruence $\sum_k\nu_kh_k'$ is $\succ$-tamed by $\phi\starop_1\psi$.

  Moreover, strong compatibility implies that $\phi\starop_1\psi\succ x_i\starop_1\psi$ (since $f=\ssf\starop_1\psi$ belongs to $\lT$)
  and $x_i\starop_1\phi\succ x_i\starop_1 z_k$ (since $x_i\starop\ssh$ belongs to $\lT$), for all suitable $i$ and $k$---hence $\phi\starop_1\psi\succ x_i\starop_1 z_k$ for all suitable $i$ and $k$, by transitivity of $\succ$.
  It follows that $h_1=\sum_i\lambda_i(x_i\starop_1\ssh)$ admits a decomposition as a $\lT$-congruence $\succ$-tamed by $\phi\starop_1\psi$.
  Similarly, $h_2=\sum_j\lambda_j(y_j\starop_1\ssh)$ admits a $\lT$-congruence $\succ$-tamed by $\phi\starop_1\psi$, using the hypotheses that $g=\ssg\starop_1\psi$ and $y_i\starop_1\ssh$ (for all suitable $i$) belong to $\lT$.
  
  We conclude that the triple $(h_1,\sum_k\nu_kh_k',h_2)$ defines a $\lT$\nbd-congruence for $(f,g)$, $\succ$-tamed by its source $\phi\starop_1\psi$.
\end{proof}

\begin{lemma}[independent rewriting lemma, context-agnostic case]
  \label{lem:HLRSM_independent_rewriting_context_agnostic}
  Let $\sS=(\sR,\sE)$ be an adapted \HLRSM{} and $\succ$ a strongly compatible esp-order.
  If for each $1\leq k\leq n$, the branching $(\ssf\starop_1z_k,\ssg\starop_1z_k)$ is $\succ$\nbd-tamely $\lT$\nbd-congruent,
  then $(f,g)$ is $\succ$\nbd-tamely $\lT$\nbd-congruent.
\end{lemma}

\subsection{Summary}
\label{subsec:rewriting_summary}

Given a presented linear 2\nbd-sesqui\-ca\-te\-gory, how to use rewriting theory to find a basis? This section summarizes the main tools developed throughout this section.
While we focus on the linear case (which is more involved), the same ideas apply for higher (non-linear) rewriting modulo.

\subsubsection{The setup}
\label{subsubsec:HLRSM_setup}

Let $\cC$ be a linear 2\nbd-sesqui\-ca\-te\-gory presented by a linear 3\nbd-sesqui\-poly\-graph $\sP$ in the sense of \cref{defn:presentation_graded_two_cat}. Assume the following choices of data have been made:
\begin{enumerate}[(a)]
  \item a splitting of $\sP$ as $\sP=\sR\sqcup\sE$, defining a \HLRSM{} $\sS=(\sR,\sE)$ presenting $\cC$ (see \cref{defn:HLRSM}),
  \item a linear sub-system $\lT\subset\sS$ (see \cref{defn:HLRSM_subsystem}),
  \item a preorder $\succ$ strongly compatible with $\lT$ (see \cref{defn:strong_compatibility}).
\end{enumerate}
Typically, $\cC$ would be a graded-2\nbd-cate\-gory and $\sE$ would contain the 3\nbd-sesqui\-polygraph of graded-interchangers, making the data $\sS=(\sR,\sE)$ a linear Gray rewriting system modulo.
To find a basis using the \LRSMbasisfromconvergencetheorem{}, the following needs to be checked:
\begin{enumerate}[(\greek*)]
  \item \emph{$\sS$ and $\lT$ present the same underlying (family of) module(s)}:
  it suffices to show that every $(\sS\setminus\lT)^+$\nbd-rewriting step is $\lT$\nbd-congruent, so that two vectors are $\sS$-congruent if and only if they are $\lT$-congruent (see \cref{lem:LRSM_congruence_describe_module}).

  \item \emph{$\scl(\sE)^\top$ is scalar-coherent on $\cB\NF_\lT$}: recall that $\cB\NF_\lT$ denotes the set of monomial $\lT$-normal forms (\cref{defn:LRSM_monomial_normal_form}), and that $\scl(\sE)^\top$ is scalar-coherent on $\cB\NF_\lT$ (\cref{defn:general_coherence_modulo,rem:LRSM_monomial_scalar_ARSM}) if for all $b\in\cB\NF_\lT$, the existence of an $\lE$-congruence $b\sim_\sE \lambda b$ for some scalar $\lambda\in\Bbbk$ implies that $\lambda=1$.
  
  Showing scalar-coherence of $\scl(\sE)^\top$ on $\cB\NF_\lT$ can be done either using ad-hoc arguments, or by apply higher rewriting modulo theory to $\scl(\sE)^\top$ (see \cref{subsubsec:scalar_HRSM}), together with \cref{prop:coherence_from_convergence} (coherence from convergence).

  \item \emph{$\lT^+$ is convergent:} we use the \LRSMtamednewmannlemma{}, showing on one hand the order $\succ$ is well-founded,\footnote{While this will not be the case in this paper, defining a compatible well-founded order may be hard in general. For this problem, the method of derivation of Guiraud and Malbos may be useful \cite{GM_PolygraphsFiniteDerivation_2018}.} and on the other hand that every monomial $\lT$\nbd-branching is $\succ$\nbd-tamely $\lT$\nbd-congruent.
\end{enumerate}
The analysis of monomial $\lT$\nbd-branching is the hardest task, which we now discuss.
Of course, in the context-agnostic case $\lT=\sS$, ($\alpha$) is automatic.

\begin{remark}
  In the example of superadjunction $\miniP$ given in \cref{subsec:intro_example}, we have:
  \begin{enumerate}[(\greek*)]
    \item we are in a context-agnostic setting, so ($\alpha$) is automatic;
    \item using the coherence theorem for interchangers \cite[Theorem~A.3.1]{Schelstraete_OddKhovanovHomology_2024}, one can check that $\scl(\miniE)^\top$ is scalar-coherent on $\miniP_2^*$, and a fortiori on $\cB\NF_{\miniS}$;
    \item as said in \cref{subsubsec:intro_tamed_congruence}, the order $\succ$ compares the number of generating 2-cells.
  \end{enumerate}
\end{remark}

\subsubsection{How to classify monomial branchings?}
\label{subsubsec:HLRSM_general_strategy}

Working modulo makes classifying monomial branchings difficult, as it considerably increases the number of rewriting steps.
However, in the context of diagrammatic algebra, the modulo data typically has a topological interpretation, which can be leveraged: e.g.\ rectilinear isotopies when working modulo interchangers, or planar isotopies when working modulo a pivotal structure.
We describe here a general strategy to classify local branchings, to be adapted depending on the example at hand:
\begin{enumerate}[(i)]
  \item \emph{Understand coherence of the modulo.}
  In other words, provide a topological or combinatorial description of when two monomial 2-cells are projectively $\sE$-congruent.
  this can be done either using ad-hoc arguments, or by applying higher rewriting modulo theory to $\scl(\sE)$. In practice, $(\beta)$ in the previous section comes as a byproduct.

  \item \emph{Describe naturalities of the modulo.}
  The modulo typically captures some underlying categorical structure, which should come with natural compatibilities with the rewriting steps of the rewriting system.
  For instance, a \GRSM{} comes with interchange naturalities. We call these \emph{$\sE$\nbd-naturalities}.

  \item \emph{Characterize rewriting steps modulo.}
  In other words, provide a topological or combinatorial description of when two rewriting steps are $\sE$\nbd-congruent. (One needs not show that this characterization captures \emph{all} ${\sE}$\nbd-congruences.)
  To do so, use the coherence statement to express isotopies in terms of $\sE$\nbd-naturalities.
  In fact, given the \ARSMbranchwisetamedcongruencelemma{}, one can instead characterize  when two rewriting steps are $\succ$\nbd-tamed $\lT$\nbd-congruent.

  \item \emph{Classify monomial branchings.}
  Using the characterization above, provide a list of monomial $\lT$\nbd-branch\-ings, called \emph{critical branchings}, such that every monomial $\lT$\nbd-bran\-ching rewrites into (a linear combination of) either an independent branching or a branchwise $\succ$\nbd-tamely $\lT$\nbd-congruent to a contextualization of a critical branching.
\end{enumerate}
In the last point (iv), three main tools are at play: the \HLRSMindependentrewritinglemma{}, the \ARSMbranchwisetamedcongruencelemma{} and the \HLRSMcontextualizationlemma{}. The results of \cref{subsubsec:HLRSM_independent_branching} also help to deal with independent branchings.
Hopefully, this leaves only a few critical branchings for which an explicit computation is needed, and one can conclude that every monomial $\lT$\nbd-branching is $\succ$\nbd-tamely $\lT$\nbd-congruent.

\begin{remark}
  In the example of superadjunction $\miniP$ given in \cref{subsec:intro_example}, this process was described in \cref{subsubsec:intro_classify_branching}.
  We have that
  (i) as said above, $\scl(\miniE)^\top$ is scalar-coherent,
  (ii) naturalities are interchange naturalities,
  (iii) the characterization was given in \cref{lem:intro_characterization_miniS},
  and (iv) the classification was given in \cref{lem:intro_classification_miniS}.
\end{remark}

\subsection{Notes}
\label{subsec:notes}

We gather here further comments for the experts, including some comments on the literature.

\begin{remark}
  \label{rem:adaptedness}
  The ($\lE$\nbd-)adaptedness condition (see above \cref{defn:linear-rewriting-system} and \cref{defn:linear-rewriting-system-modulo}) prevents ``obvious'' obstructions to termination. Indeed, if $s(r)\in\projsupp_\lE(t(r))$ then there exists an infinite sequence of positive $\lS$\nbd-rewriting steps, all of type $r$.
  Moreover, without this assumption the Church--Rosser property for positive rewritings steps (\cref{lem:LRSM_Church-Rosser_positive_rw}) does not hold; see \cref{rem:LRSM_counterexample_CR_pos_rw}.
  Given how fundamental this result is, we choose to enforce the adaptedness condition in the definition, in contrast with the abstract case.
  This is only a choice of presentation; indeed, in practice one eventually wishes to work with terminating rewriting systems, which implies adaptedness.
  In \cite{GHM_ConvergentPresentationsPolygraphic_2019}, ``left-monomial'' encompasses both our ``left-monomial'' and ``adapted'' notions.
  However, the adaptedness condition is dropped in \cite{Alleaume_RewritingHigherDimensional_2018} and in \cite{Dupont_RewritingModuloIsotopies_2022,Dupont_RewritingModuloIsotopies_2021}; in particular, Lemma~4.2.9 in \cite{Alleaume_RewritingHigherDimensional_2018} and Lemma~1.1.5 in \cite{Dupont_RewritingModuloIsotopies_2022} are not correct as stated.
  This motivates our change of terminology, hoping to avoid further confusion.
  Note that while the adaptedness condition is easy to check in the non-modulo setting, it can be more involved in the modulo setting, as one needs to scan through projective $\lE$\nbd-congruence classes.
\end{remark}

\begin{remark}
  \label{footnote:positivity_depends_on_modulo}
  In our approach to linear rewriting modulo, the modulo data $\lE$ must be monomial-invertible. Extending the theory to allow non-monomial modulos is an interesting open problem.
  Indeed, the positivity condition in \cref{defn:LRSM_positive_rewriting_steps} is $s(r)\notin\projsupp_\lE(v)$, and not $s(r)\notin\supp(v)$. In other words, a positive $\lS$\nbd-rewriting step is \emph{not} a positive $\lR$\nbd-rewriting step pre- and post-composed with $\lE$-congruence. The positivity condition is stronger, and depends on $\lE$. Otherwise, positive $\lS$\nbd-rewriting steps do not provide a suitable solution to termination.
  (consider the \LRSM{} $\lS=(\cB;\lR,\lE)$ with $\cB=\{a,a',b\}$, $\lR=\{a\to b\}$ and $\lE=\{a\sim a'\}$, and the $\lS$\nbd-rewriting sequence $0=a-a\sim_\lE a-a'\dashrightarrow_\lR b-a'\sim_\lE b-a\dashrightarrow_\lR b-b=0$.)
  In particular, we need $\lE$ to be monomial-invertible for the positivity condition to make sense.

  In \cite{Dupont_RewritingModuloIsotopies_2022,Dupont_RewritingModuloIsotopies_2021}, a positive $\lS$\nbd-rewriting step \emph{is} defined as a composition as above such that the middle arrow is a positive $\lR$\nbd-rewriting step, and $\lE$ is not constrained to be monomial; it is unclear to us how Dupont's approach prevents the above obstruction.
  Moreover, as far as we understand, the work \cite{DEL_SuperRewritingTheory_2024} uses Dupont's approach with a non-monomial $\lE$.
\end{remark}

\begin{remark}
    \label{footnote:same-underlying-two-cell-for-modulo}
  Recall that in each of our definitions of \HLRSM{}, we impose that $\sR$ and $\sE$ have the same underlying 2-polygraph.
  In \cite{Dupont_RewritingModuloIsotopies_2022,Dupont_RewritingModuloIsotopies_2021}, Dupont allows the more general definition $\sE_2\subset\sR_2$ in the context of linear 3-polygraph.
    We prefer to avoid this generality, as it leads to ambiguity for statements like ``$\sE$ is convergent''.
    For instance, for the linear 3-polygraphs $E$ and $R$ in \cite{Dupont_RewritingModuloIsotopies_2022,Dupont_RewritingModuloIsotopies_2021}, we have that $E$ is convergent when viewed on $E_2$, but not when viewed on $R_2$.
    Similarly, it is unclear to us why the following branching between two Reidemeister 3 moves
    \newcommand{\diagid}[2]{
      \draw (#1,#2) to (#1,#2+1);
    }
    \newcommand{\cro}[2]{
      \draw (#1,#2+0) to (#1+1,#2+1);
      \draw (#1,#2+1) to (#1+1,#2+0);
    }
    \newcommand{\cc}[2]{
      \draw (#1,#2+1) to[out=-90,in=180] (#1+.5,#2+1-.3) to[out=0,in=-90] (#1+1,#2+1);
      \draw (#1,#2) to[out=90,in=180] (#1+.5,#2+.3) to[out=0,in=90] (#1+1,#2);
    }
    \begin{gather*}
      \tikzpic{
        \cro{0}{4}\diagid{2}{4}\diagid{3}{4}
        \diagid{0}{3}\cro{1}{3}\diagid{3}{3}
        \cro{0}{2}\cc{2}{2}
        \diagid{0}{1}\cro{1}{1}\diagid{3}{1}
        \cro{0}{0}\diagid{2}{0}\diagid{3}{0}
      }[scale=.5,yscale=.7]
    \end{gather*}
    does not seem to be considered in \cite{Dupont_RewritingModuloIsotopies_2021}. On a related note, it is unclear to us how the same branching does not contradict Lemma~5.7 in \cite{DEL_SuperRewritingTheory_2024}.
\end{remark}

\begin{remark}
  \label{rem:quasi-normal-forms}
  In \cite{Alleaume_RewritingHigherDimensional_2018,Alleaume_HigherdimensionalLinearRewriting_2018}, \cite{Dupont_RewritingModuloIsotopies_2022,Dupont_RewritingModuloIsotopies_2021} and \cite{DEL_SuperRewritingTheory_2024}, the authors used the notion of quasi-termi\-nation to deal with context-dependent termination issues, and suggested a basis theorem based on quasi-terminating normal forms.
  Context-dependent rewriting gives a different way to deal with the same issue, avoiding the use of quasi-terminating normal forms altogether.
  It is also more general, as it applies to settings that are not quasi-terminating, such as the setting of graded $\glt$-foams studied in \cref{sec:rewriting_foam}.

  Quasi-terminating normal forms are typically badly behaved.
  For instance, a monomial in the support of a quasi-normal form needs not be a quasi-normal form, and a linear combination of quasi-normal forms needs not be a quasi-normal form; compare with \cref{rem:lrs_normal_forms}.
  It is unclear how the authors of \cite{Alleaume_RewritingHigherDimensional_2018,Alleaume_HigherdimensionalLinearRewriting_2018,Dupont_RewritingModuloIsotopies_2022,Dupont_RewritingModuloIsotopies_2021,DEL_SuperRewritingTheory_2024} extract a basis from a quasi-terminating system, and in particular, how \cite[Theorem~2.5.6]{Dupont_RewritingModuloIsotopies_2022} and \cite[Theorem~2.32]{DEL_SuperRewritingTheory_2024} should be interpreted.
  For instance, one can consider the quasi-terminating linear rewriting system (non-modulo) $(\cB=\{a,b,c\},\lR=\{a\to b+c,b\to a-c\})$: the set of monomial quasi-normal forms is $\cB\NF_\lR=\cB=\{a,b,c\}$, which is not a basis of the module presented by $\lR$.
\end{remark}

\begin{remark}
  \label{footnote:gap_critical_confluence}
  We give an example that in the setting of \cref{cor:HLRSM_contextualization_branching_lemma}, the confluence $\Gamma[f',g']$ need not be positive;
  it is taken from \cref{lem:foam_spatial_B-confluence}.
  Recall that plain arrows denote positive rewriting steps, while dashed arrows denote not-necessarily-positive rewriting steps.
  In the example below, the graded interchange law is part of the modulo.
  On the left-hand side, the branches of the confluence (labelled ``$\dd$'') evaluate
  $\tikzpic{
    \node[fdot1=1.5pt] at (0,0){};
    \node[fdot1=1.5pt] at (0,.5){};
  }[scale=.5]$ to zero, so that only the other diagram remains.
  The two diagrams on the top (resp.\ on the bottom) are \emph{not} projectively congruent modulo, so that both branches of the confluence are positive.
  The same confluent branching is pictured on the right-hand side, contextualized with a cap and cup. Because the graded interchange law is part of the modulo, dots can now freely move between the top and the bottom. In particular, the two diagrams on the top (resp.\ on the bottom) are now projectively congruent modulo, and the branches of the confluence are not positive anymore.
  \begin{figure}
    \begin{gather*}
    \label{eq:example_contextualization_positivity}
    \def\scl{.4}
    \def\yscl{.8}
    \newcommand{\tempdiag}{
      \funzip[-1][3]\frst[1][3]
      \fzip[-1][-2]\frst[1][-2]
      \begin{scope}[yscale=4]
        \frst[-1][-1/4]
      \end{scope}
    }
    \begin{tikzcd}[ampersand replacement=\&,row sep=-1em,column sep=tiny]
      \&
      \tikzpic{
        \flst[0][2]\fdot[.5][2.5][1]\frst[1][2]
        \fdot[.5][2-.2][1]
        \fcap\fcup[0][1]
        \flst[0][-1]\frst[1][-1]
      }[scale=\scl,yscale=\yscl]
      +
      \tikzpic{
        \flst[0][2]\fdot[.5][2.5][1]\frst[1][2]
        \fcap\fcup[0][1]
        \fdot[.5][.2][1]
        \flst[0][-1]\frst[1][-1]
      }[scale=\scl,yscale=\yscl]
      \&
      \\
      \tikzpic{
        \begin{scope}[yscale=4]
          \flst[0][-1/4]\frst[1][-1/4]
        \end{scope}
        \fdot[.5][1][1]
      }[scale=\scl,yscale=\yscl]
      \&\&
      \tikzpic{
        \flst[0][2]\frst[1][2]
        \fdot[.5][2-.2][1]
        \fcap\fcup[0][1]
        \fdot[.5][.2][1]
        \flst[0][-1]\frst[1][-1]
      }[scale=\scl,yscale=\yscl]
      \\
      \&
      \tikzpic{
        \flst[0][2]\frst[1][2]
        \fdot[.5][2-.2][1]
        \fcap\fcup[0][1]
        \flst[0][-1]\fdot[.5][-.5][1]\frst[1][-1]
      }[scale=\scl,yscale=\yscl]
      +
      \tikzpic{
        \flst[0][2]\frst[1][2]
        \fcap\fcup[0][1]
        \fdot[.5][.2][1]
        \flst[0][-1]\fdot[.5][-.5][1]\frst[1][-1]
      }[scale=\scl,yscale=\yscl]
      \&
      \arrow[from=2-1,to=1-2,"\foamlT"{subscript,pos=.85},"\nc",bend left=20]
      \arrow[from=2-1,to=3-2,"\foamlT"{subscript,pos=.95},"\nc"',bend right=20]
      \arrow[from=1-2,to=2-3,"\foamlT"{subscript},"\dd",bend left=20]
      \arrow[from=3-2,to=2-3,"\foamlT"{subscript},"\dd"',bend right=20]
    \end{tikzcd}
    \qquad
    \begin{tikzcd}[ampersand replacement=\&,row sep=-1em,column sep=tiny]
      \&
      \tikzpic{
        \flst[0][2]\fdot[.5][2.5][1]\frst[1][2]
        \fdot[.5][2-.2][1]
        \fcap\fcup[0][1]
        \flst[0][-1]\frst[1][-1]
        \tempdiag
      }[scale=\scl,yscale=\yscl]
      +
      \tikzpic{
        \flst[0][2]\fdot[.5][2.5][1]\frst[1][2]
        \fcap\fcup[0][1]
        \fdot[.5][.2][1]
        \flst[0][-1]\frst[1][-1]
        \tempdiag
      }[scale=\scl,yscale=\yscl]
      \&
      \\
      \tikzpic{
        \begin{scope}[yscale=4]
          \flst[0][-1/4]\frst[1][-1/4]
        \end{scope}
        \fdot[.5][1][1]
        \tempdiag
      }[scale=\scl,yscale=\yscl]
      \&\&
      \tikzpic{
        \flst[0][2]\frst[1][2]
        \fdot[.5][2-.2][1]
        \fcap\fcup[0][1]
        \fdot[.5][.2][1]
        \flst[0][-1]\frst[1][-1]
        \tempdiag
      }[scale=\scl,yscale=\yscl]
      \\
      \&
      \tikzpic{
        \flst[0][2]\frst[1][2]
        \fdot[.5][2-.2][1]
        \fcap\fcup[0][1]
        \flst[0][-1]\fdot[.5][-.5][1]\frst[1][-1]
        \tempdiag
      }[scale=\scl,yscale=\yscl]
      +
      \tikzpic{
        \flst[0][2]\frst[1][2]
        \fcap\fcup[0][1]
        \fdot[.5][.2][1]
        \flst[0][-1]\fdot[.5][-.5][1]\frst[1][-1]
        \tempdiag
      }[scale=\scl,yscale=\yscl]
      \&
      \arrow[from=2-1,to=1-2,"\foamlT"{subscript,pos=.85},"\nc",bend left=20]
      \arrow[from=2-1,to=3-2,"\foamlT"{subscript,pos=.95},"\nc"',bend right=20]
      \arrow[from=1-2,to=2-3,"\foamlT"{subscript},"\dd",bend left=20,dashed]
      \arrow[from=3-2,to=2-3,"\foamlT"{subscript},"\dd"',bend right=20,dashed]
    \end{tikzcd}
    \end{gather*}
    \caption{In the setting of \cref{cor:HLRSM_contextualization_branching_lemma}, the confluence $\Gamma[f',g']$ need not be positive.}
    \label{fig:example_positivity_and_confluence}
  \end{figure}
  In particular, there are related gaps in the proofs of \cite[lemma~4.2.12]{Alleaume_RewritingHigherDimensional_2018} and \cite[theorem 2.2.9]{Dupont_RewritingModuloIsotopies_2022}.
\end{remark}

\begin{remark}
  \label{rem:mistake_article_dupont}
  \newcommand{\diagid}[2]{
    \draw (#1,#2) to (#1,#2+1);
  }
  \newcommand{\cro}[2]{
    \draw (#1,#2+0) to (#1+1,#2+1);
    \draw (#1,#2+1) to (#1+1,#2+0);
  }
  \newcommand{\diagcap}[2]{
    \draw (#1,#2) to[out=90,in=180] (#1+.5,#2+.3) to[out=0,in=90] (#1+1,#2);
  }
  \newcommand{\diagdot}[2]{
    \node[fill=black,circle,inner sep=1.5pt] at (#1,#2) {};
  }
  In \cite{Dupont_RewritingModuloIsotopies_2021}, Dupont developed a rewriting approach to categorified quantum groups.
  In our understanding, the rewriting rules (\cite[Definition~2.3.1 and Definition~3.2.1, iv), 1)]{Dupont_RewritingModuloIsotopies_2021}) and the modulo (\cite[Section~3.2.5]{Dupont_RewritingModuloIsotopies_2021}) respectively contain, in particular, the following 3-cells:
  \begin{gather*}
    \tikzpic{\cro{0}{0}\diagdot{.25}{.75}}[scale=.4]
    \to
    \tikzpic{\cro{0}{0}\diagdot{.75}{.25}}[scale=.4]
    +
    \tikzpic{\diagid{0}{0}\diagid{1}{0}}[scale=.4]
    \qquad
    \tikzpic{\cro{0}{0}\diagdot{.75}{.75}}[scale=.4]
    \to
    \tikzpic{\cro{0}{0}\diagdot{.25}{.25}}[scale=.4]
    -
    \tikzpic{\diagid{0}{0}\diagid{1}{0}}[scale=.4]
    \qquad\an\qquad
    \tikzpic{\diagcap{0}{0}\diagdot{.2-.05}{.2}}[xscale=.4]\sim \tikzpic{\diagcap{0}{0}\diagdot{.8+.05}{.2}}[xscale=.4]
  \end{gather*}
  Given this definition, it is unclear to us how to define a confluence for the following branching:
  \begin{gather*}
    \ldots\leftarrow
    \tikzpic{\cro{0}{0}\diagcap{1}{1}\cro{2}{0}\diagdot{1}{1}}[scale=.4]
    \sim
    \tikzpic{\cro{0}{0}\diagcap{1}{1}\cro{2}{0}\diagdot{2}{1}}[scale=.4]
    \to\ldots
  \end{gather*}
\end{remark}

\section{A basis for graded \texorpdfstring{$\glt$}{gl2}-foams}
\label{sec:rewriting_foam}

In this section, we apply linear Gray rewriting modulo as developed in \cref{sec:foundation_rewriting} to a certain graded-2-category $\gfoam_d$, the \emph{graded-2-category of $\glt$-foams}, and show that it has the expected basis.
The strategy follows the blueprint given in \cref{subsubsec:HLRSM_general_strategy}.
\Cref{subsec:review_foams} reviews $\gfoam_d$ (\cref{defn:graded2cat_foams}) and states the basis theorem (\cref{thm:foam_basis_theorem}).
\Cref{subsec:convergent_pres_gltfoam} then defines the working data for the rewriting theory.
The core of the proof is given in \cref{subsec:rewriting_foam_coherence_iso} and \cref{subsec:rewriting_foam_confluence_modulo_iso}, respectively dealing with coherence of the modulo data and confluence of monomial branchings.
We conclude with an addendum (\cref{subsec:addendum_deformation_foam}), showing that $\gfoam_d$ admits a variant $\gfoam_d'$ satisfying the same basis theorem.

\subsection{Graded \texorpdfstring{$\glt$}{gl2}-foams}
\label{subsec:review_foams}

\subsubsection{Presentation by a linear Gray polygraph}

Fix a positive integer $d\in\bN$.
We define a certain linear Gray polygraph $\sgfoam_d$.
The objects of $\sgfoam_d$ are\footnote{In \cite{SV_OddKhovanovHomology_2023}, $(\sgfoam_d)_0$ is denoted $\underline{\Lambda}_d$.}
\begin{gather*}
  (\sgfoam_d)_0\coloneqq
  \bigsqcup_{k\in\bN}\{\lambda\in\{1,2\}^k\mid \lambda_1+\ldots+\lambda_k=d\}.
\end{gather*}
For each $\lambda\in (\sgfoam_d)_0$ with $k$ coordinates, we label the $r$th coordinate by $l_\lambda(r)=\sum_{s<r}\lambda_s+1$.
For instance:
\[(\underset{1}{1},\underset{2}{2},\underset{4}{1},\underset{5}{1})\]
where the value of $l_\lambda$ is written under each coordinate.
In other words, the label $l_\lambda(r)$ can be understood as a weighted coordinate, where coordinates with value $2$ count double.
We write $r=l_\lambda^{-1}(i)$ whenever $i=l_\lambda(r)$, keeping in mind that $l_\lambda^{-1}$ is not well-defined for all $1\leq i\leq d$.

The set of 1-cells $(\sgfoam_d)_1$ is given by the following, for all $1\leq i<d$:
\begin{gather*}
  (\ldots,\underset{i}{1},\underset{i+1}{1},\ldots)
  \mspace{10mu}
  \tikzpic{\flst\node[below=-2pt] at (0,0) {\small $i$};}[yscale=1.5][(0,.75)]
  \mspace{10mu}
  (\ldots,\underset{i}{2},\ldots)
  \;\quad\an\quad\;
  (\ldots,\underset{i}{2},\ldots)
  \mspace{10mu}
  \tikzpic{\frst\node[below=-2pt] at (0,0) {\small $i$};}[yscale=1.5][(0,.75)]
  \mspace{10mu}
  (\ldots,\underset{i}{1},\underset{i+1}{1},\ldots)
\end{gather*}
We read 1\nbd-morphism from right to left: for instance if $d=2$, then the first 1\nbd-cell is a 1\nbd-morphism from $(2)$ to $(1,1)$.
The label $i$, called the \emph{colour} of the strand, has value $1\leq i\leq d-1$.
One side of the strand is shaded: this encodes the fact that $\lambda_{l_\lambda^{-1}(i)}=2$, where $\lambda$ labels the shaded region.
The label of one region determines the label on the other region.

\begin{terminology}
  \label{term:i_strand}
  A strand labelled by a colour $i$ is called an \emph{$i$-strand} and a region shaded with a colour $i$ is called an \emph{$i$-shading}.
\end{terminology}

We shall omit to label regions with objects, and instead record shadings. For instance, the two diagrams
\begin{gather*}
\tikzpic{
  \draw[diag1] (0,0) node[below]{\tiny $1$} to (0,1.5);
  \draw[diag3] (1,0) node[below]{\tiny $3$} to (1,1.5);
  \draw[diag1] (2,0) node[below]{\tiny $1$} to (2,1.5);
  \draw[diag2] (3,0) node[below]{\tiny $2$} to (3,1.5);
  \node at (4.2,.75) {\small $(1,2,1,1)$};
}[][(0,.75)]
\qquad\an\qquad
\tikzpic{
  \draw[diag1] (0,0) node[below]{\tiny $1$} to (0,1.5);
  \draw[diag3] (1,0) node[below]{\tiny $3$} to (1,1.5);
  \draw[diag1] (2,0) node[below]{\tiny $1$} to (2,1.5);
  \draw[diag2] (3,0) node[below]{\tiny $2$} to (3,1.5);
  \fill[color=colour_diag3,opacity=.2] (-.8,0) rectangle (1,1.5);
  \fill[color=colour_diag1,opacity=.2] (0,0) rectangle (2,1.5);
  \fill[color=colour_diag2,opacity=.2] (3,0) rectangle (3.8,1.5);
}[][(0,.75)]
\end{gather*}
encodes the same data.
Note that an $i$-shading and a $j$-shading can only overlap when $\abs{i-j}>1$.

The set of 2-cells $(\sgfoam_d)_{2}$ is given in \cref{fig:foam-2-cells}. Following the above terminology, we call an dot labelled $i$ an \emph{$i$-dot}, and similarly for other generators.
In addition to the shadings pictured in the diagrams, the diagram may carry further $k$-shadings, according to the following conditions:\label{page:label_condition}
\begin{enumerate}[(i)]
  \item $\abs{i-k}>1$ for cups, caps, zips and unzips;
  \item $\abs{i-k}>1$ and $\abs{j-k}>1$ for crossings;
  \item $i\neq k,k-1$ for dots.
\end{enumerate}
Conditions (i) and (ii) state that an $i$-shading and a $j$-shading can only overlap when $\abs{i-j}>1$, as already required for identity diagrams.
In terms of coordinates, condition (iii) translates to the condition that $\lambda_{l_\lambda^{-1}(i)}=1$, where $\lambda$ labels the region where the $i$-dot lies.

\begin{figure}
  \begin{gather*}
    \begin{tabular}{*{4}{c@{\hskip 6ex}}c}
      $\tikzpic{\fdot\node at (.2,-.2) {\scriptsize $i$};}$
      &
      $\tikzpic{\fcup\node at (1,.4) {\scriptsize $i$};}$
      &
      $\tikzpic{\fcap\node at (1,.4) {\scriptsize $i$};}$
      &
      $\tikzpic{\fzip\node at (1,.6) {\scriptsize $i$};}$
      &
      $\tikzpic{\funzip\node at (1,.6) {\scriptsize $i$};}$
      \\*[2ex]
      $\degd$ & $\degrcu$ & $\deglca$ & $\deglcu$ & $\degrca$
      \\*[1ex]
      {\small dot} & {\small cup} & {\small cap} & {\small zip} & {\small unzip}
    \end{tabular}
    \\[2ex]
    \begin{tabular}{*{4}{c@{\hskip 6ex}}c}
      $\tikzpic{
        \fbcro
        \node[left=-3pt] at (0,0) {\scriptsize $i$};
        \node[right=-3pt] at (1,0){\scriptsize $j$};
      }$
      &
      $\tikzpic{
        \flcro
        \node[left=-3pt] at (0,0) {\scriptsize $i$};
        \node[right=-3pt] at (1,0){\scriptsize $j$};
      }$
      &
      $\tikzpic{
        \ftcro
        \node[left=-3pt] at (0,0) {\scriptsize $i$};
        \node[right=-3pt] at (1,0){\scriptsize $j$};
      }$
      &
      $\tikzpic{
        \frcro
        \node[left=-3pt] at (0,0) {\scriptsize $i$};
        \node[right=-3pt] at (1,0){\scriptsize $j$};
      }$
      &
      $\mspace{40mu}\text{if }\abs{i-j}>1$
      \\*[3ex]
      $(0,0)$ & $(0,0)$ & $(0,0)$ & $(0,0)$ &
    \end{tabular}
    \\
    \text{\small crossings}
  \end{gather*}
  \caption{The set of 2-cells $(\sgfoam_d)_2$. In each case, the diagram may be shaded by further colours as described in the main text: here we only give the shadings induced by the strands in the diagram.
  Furthermore, it is always assumed that the region where the $i$-dot lies does not carry an $i$- or $i+1$-shading. Each 2-cell is equipped with a $\bZ\times\bZ$-degree, denoted below its diagram.}
  \label{fig:foam-2-cells}
\end{figure}

The following gives the structural data that defines $(G,\bilfoam)$-graded-interchangers:

\begin{definition}
  \label{defn:ring_R_bil}
  Let $\ringfoam$ be a commutative ring together with three invertible elements $X$, $Y$ and $Z$ in $\ringfoam^\times$ such that $X^2=Y^2=1$.
  Given this data, let $\bilfoam$ be the following bilinear form for the abelian group $G\coloneqq\bZ^2$:
  \begin{align*}
    \bilfoam\colon\bZ^2\times \bZ^2&\to \ringfoam^\times,\\*
    ((a,b),(c,d)) &\mapsto X^{ac}Y^{bd}Z^{ad-bc}.
  \end{align*}
\end{definition}

Finally:

\begin{definition}
  Fix $d\in\bN$.
  The linear Gray polygraph $\sgfoam_d$ has 0-cells $(\sgfoam_d)_{0}$ and 1-cells $(\sgfoam_d)_{1}$ as given above, 2-cells given in \cref{fig:foam-2-cells} and 3-cells
  \[(\sgfoam_d)_3=\foamE_3\sqcup\foamR_3,\]
  where the sets $\foamE_3$ and $\foamR_3$ are given in \cref{fig:foam_E} and \cref{fig:foam_R}, respectively.
\end{definition}

The diagrams encode certain singular surfaces called \emph{$\glt$\nbd-foams}.
A foam is made of \emph{1\nbd-facets} and \emph{2\nbd-facets}; two 1\nbd-facets and one 2-facet can join at a singular line called a \emph{seam}.
The string diagrammatics encodes $\glt$-foams via their 2-facets, as follows:
\begin{IEEEeqnarray*}{CcC}
  \figfoam[scale=.15]{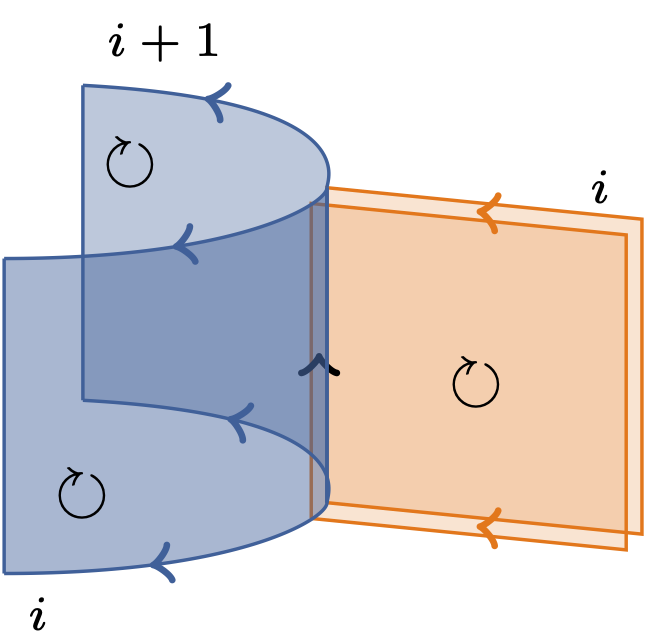}
  &
  \mspace{60mu}
  \leftrightarrow
  \mspace{30mu}
  &
  \tikzpic{\frst}[scale=2]
  \\*
  \text{\small $\glt$-foams}
  &&
  \text{\small shading diagrammatics}
\end{IEEEeqnarray*}
On $\glt$-foams, the fact that $i$-dots cannot sit on $i$- or $i+1$-shadings translates to the fact that the dots cannot sit on 2-facets.

\begin{figure}
  \centering
  \captionsetup[subfigure]{justification=centering}

  \begin{subfigure}{\textwidth}
    \begin{gather*}
      \tikzpic{
        \draw (0,2) to
          node[black_dot,pos=.3] {}
          node[right,pos=.3] {\scriptsize $\alpha$}
          (0,0);
        \draw (1,2) to
          node[black_dot,pos=.7] {}
          node[right,pos=.7] {\scriptsize $\beta$}
          (1,0);
      }
      \;\overset{\bilfoam(\deg\alpha,\deg\beta)}{\Rrightarrow}\;
      \tikzpic{
        \draw (0,2) to
          node[black_dot,pos=.7] {}
          node[right,pos=.7] {\scriptsize $\alpha$}
          (0,0);
        \draw (1,2)  to
          node[black_dot,pos=.3] {}
          node[right,pos=.3] {\scriptsize $\beta$}
          (1,0);
      }
      \\
      \text{\scriptsize for all $\alpha,\beta\in(\sgfoam_d)_2^*$}
    \end{gather*}
    \caption{The set of 3-cells $\foamGray_3$ of $(\bZ^2,\bilfoam)$-graded interchangers.}
    \label{subfig:foam_FGray}
  \end{subfigure}

  \begin{subfigure}{\textwidth}
    \def\scl{.8}
    \begin{gather*}
      \begingroup
      \tikzpic{
        \draw[foamdraw1]  +(0,-.75) node[below] {\textcolor{black}{\scriptsize $i$}}
          .. controls (0,-.375) and (1,-.375) .. (1,0)
          .. controls (1,.375) and (0, .375) .. (0,.75);
        \draw[foamdraw2]  +(1,-.75) node[below] {\textcolor{black}{\scriptsize $j$}}
          .. controls (1,-.375) and (0,-.375) .. (0,0)
          .. controls (0,.375) and (1, .375) .. (1,.75);
      }[scale=.8*\scl][(0,0)]
      \Rrightarrow
      \tikzpic{
        \draw[foamdraw1] (0,-.75) node[below] {\textcolor{black}{\scriptsize $i$}} to (0,.75);
        \draw[foamdraw2] (1,-.75) node[below] {\textcolor{black}{\scriptsize $j$}} to (1,.75);
      }[scale=.8*\scl][(0,0)]
      \mspace{60mu}
      %
      \tikzpic{
        \draw[foamdraw1]  +(0,0)node[below] {\textcolor{black}{\scriptsize $i$}}
          .. controls (0,0.5) and (2, 1) ..  +(2,2);
        \draw[foamdraw3]  +(2,0)node[below] {\textcolor{black}{\scriptsize $k$}}
          .. controls (2,1) and (0, 1.5) ..  +(0,2);
        \draw[foamdraw2]  (1,0)node[below] {\textcolor{black}{\scriptsize $j$}}
          .. controls (1,0.5) and (0, 0.5) ..  (0,1)
          .. controls (0,1.5) and (1, 1.5) ..  (1,2);
      }[scale=.7*\scl][(0,1*.7*\scl)]
      \Rrightarrow
      \tikzpic{
        \draw[foamdraw1]  +(0,0)node[below] {\textcolor{black}{\scriptsize $i$}}
          .. controls (0,1) and (2, 1.5) ..  +(2,2);
        \draw[foamdraw3]  +(2,0)node[below] {\textcolor{black}{\scriptsize $k$}}
          .. controls (2,.5) and (0, 1) ..  +(0,2);
        \draw[foamdraw2]  (1,0)node[below]{\textcolor{black} {\scriptsize $j$}}
          .. controls (1,0.5) and (2, 0.5) ..  (2,1)
          .. controls (2,1.5) and (1, 1.5) ..  (1,2);
      }[scale=.7*\scl][(0,1*.7*\scl)]
      \endgroup
      \\[2ex]
      \begingroup
        \tikzpic{
          \draw[foamdraw1](-.5,.4) to (0,-.3);
          \draw[foamdraw2] (0.3,-0.3)
            to[out=90, in=0] (0,0.2)
            to[out = -180, in = 40] (-0.5,-0.3);
          \node at (-0.5,-.5) {\scriptsize $j$};
          \node at (0,-.5) {\scriptsize $i$};
        }[scale=1.2*\scl][(0,0)]
        \Rrightarrow
        \tikzpic{
          \draw[foamdraw1](.6,.4) to (.1,-.3);
          \draw[foamdraw2] (0.6,-0.3)
            to[out=140, in=0] (0.1,0.2)
            to[out = -180, in = 90] (-0.2,-0.3);
          \node at (-0.2,-.5) {\scriptsize $j$};
          \node at (0.1,-.5) {\scriptsize $i$};
        }[scale=1.2*\scl][(0,0)]
        \mspace{40mu}
        \tikzpic{
          \draw[foamdraw1](-.5,-.3) to (0,.4);
          \draw[foamdraw2] (0.3,0.4)
            to[out=-90, in=0] (0,-0.1)
            to[out = 180, in = -40] (-0.5,0.4);
          \node at (-0.5,.6) {\scriptsize $j$};
          \node at (-0.05,.6) {\scriptsize $i$};
        }[scale=1.2*\scl][(0,0)]
        \Rrightarrow
        \tikzpic{
          \draw[foamdraw1](.6,-.3) to (.1,.4);
          \draw[foamdraw2] (0.6,0.4)
            to[out=-140, in=0] (0.1,-0.1)
            to[out = 180, in = -90] (-0.2,0.4);
          \node at (-0.25,.6)  {\scriptsize $j$};
          \node at (0.15,.6) {\scriptsize $i$};
        }[scale=1.2*\scl][(0,0)]
      \endgroup
      \mspace{40mu}
      \begin{array}{c}
        \tikzpic{
          \fdot[-1][1][2]\node at (-.7,.8) {\scriptsize $j$};
          \frst[0][1]\frst\node[left=-3pt] at (0,2) {\scriptsize $i$};
        }[scale=.5][(0,1*.5)]
        \Rrightarrow
        \tikzpic{
          \fdot[1][1][2]\node at (-.7+2,.8) {\scriptsize $j$};
          \frst[0][1]\frst\node[left=-3pt] at (0,2) {\scriptsize $i$};
        }[scale=.5][(0,1*.5)]
        \\
        \text{\footnotesize if $j\neq i,i+1$}
      \end{array}
    \end{gather*}
    \caption{The set of 3-cells $\foamX_3$, called \emph{braid moves}, \emph{pitchforks} and \emph{dot slides}.}
    \label{subfig:foam_X}
  \end{subfigure}

  \begin{subfigure}{\textwidth}
    \begin{gather*}
      \def\scl{.7}
      \tikzpic{
        \fcap[0][1]\flst[2][1]\flst\fzip[1][0]
        \node[left=-3pt] at (2,2) {\scriptsize $i$};
      }[scale=.5][(0,1*.5)]
      \overset{}{\Rrightarrow}
      \tikzpic{
        \flst[0][1]\flst\node[left=-3pt] at (0,2) {\scriptsize $i$};
      }[scale=.5][(0,1*.5)]
      \mspace{20mu}
      \tikzpic{
        \frst[0][1]\fcap[1][1]\fzip[0][0]\frst[2][0]
        \node[left=-3pt] at (0,2) {\scriptsize $i$};
      }[scale=.5][(0,1*.5)]
      \overset{X}{\Rrightarrow}
      \tikzpic{
        \frst[0][1]\frst\node[left=-3pt] at (0,2) {\scriptsize $i$};
      }[scale=.5][(0,1*.5)]
      \mspace{20mu}
      \tikzpic{
        \funzip[0][1]\frst[2][1]\frst\fcup[1][0]
        \node[left=-3pt] at (2,2) {\scriptsize $i$};
      }[scale=.5][(0,1*.5)]
      \overset{Z^2}{\Rrightarrow}
      \tikzpic{
        \frst[0][1]\frst\node[left=-3pt] at (0,2) {\scriptsize $i$};
      }[scale=.5][(0,1*.5)]
      \mspace{20mu}
      \tikzpic{
        \flst[0][1]\funzip[1][1]\fcup[0][0]\flst[2][0]
        \node[left=-3pt] at (0,2) {\scriptsize $i$};
      }[scale=.5][(0,1*.5)]
      \overset{YZ^2}{\Rrightarrow}
      \tikzpic{
        \flst[0][1]\flst\node[left=-3pt] at (0,2) {\scriptsize $i$};
      }[scale=.5][(0,1*.5)]
    \end{gather*}
    \caption{The set of 3-cells $\foamZ_3$, called \emph{zigzags}.}
    \label{subfig:foam_Z}
  \end{subfigure}

  \caption{The set of 3-cells $\foamE_3=\foamGray_3\sqcup\foamX_3\sqcup\foamZ_3$. We use the convention for scalar 3\nbd-sesqui\-poly\-graph already used in \cref{ex:polygraphC}, writing the associated scalar on top of the arrow. For the 3-cells in $\foamX_3$ involving crossings, we have a 3-cell for every choice of shadings, and colours should be so that crossings exist.}
  \label{fig:foam_E}

  \def\vs{1ex}
  \begin{IEEEeqnarray*}{CcC}
    \left(\tikzpic{
      \fdot\node at (.2,-.2) {\scriptsize $i$};
    }[][(0,0)]\right)^2
    \overset{}{\Rrightarrow}
    0
    &\mspace{50mu}&
    \tikzpic{
      \fdot[-1][1]\node at (-.7,.8) {\scriptsize $i$};
      \frst[0][1]\frst\node[left=-3pt] at (0,2) {\scriptsize $i$};
    }[scale=.5][(0,1*.5)]
    \overset{}{\Rrightarrow}\norma{}{-\;\;}{}
    \tikzpic{
      \fdot[-1][1][2]\node at (-.7,.5) {\scriptsize $i+1$};
      \frst[0][1]\frst\node[left=-3pt] at (0,2) {\scriptsize $i$};
    }[scale=.5][(0,1*.5)]
    \\[1ex]
    \text{\small dot annihilation ($\dd$)}
    &&
    \text{\small dot migration ($\dm$)}
    \\[2ex]
    \IEEEeqnarraymulticol{3}{c}{
      {}
      \tikzpic{
        \fdot[.5][1][2]\node[above right=-2pt] at (.5,1) {\scriptsize $i+1$};
        \fcap[0][1]\fcup\node at (.8+.5,-.3+1) {\scriptsize $i$};
      }[scale=.7]
      \overset{}{\Rrightarrow}
      1
      \mspace{50mu}
      \tikzpic{
        \fcap[0][1]\fcup\node at (.8+.5,-.3+1) {\scriptsize $i$};
      }[scale=.7]
      \overset{}{\Rrightarrow}
      0
      \mspace{50mu}
      \tikzpic{
        \funzip[0][1]\fzip\node at (.8+.5,-.3+1) {\scriptsize $i$};
      }[scale=.7]
      \mspace{5mu}\overset{}{\Rrightarrow}\mspace{5mu}
      \norma{}{}{Z}
      \mspace{10mu}
      \tikzpic{
        \fdot\node at (.2,-.2) {\scriptsize $i$};
      }[][(0,0)]
      \mspace{10mu}
      \norma{+XY}{-XY}{+XYZ}
      \mspace{10mu}
      \tikzpic{
        \fdot[0][0][2]\node at (.5,-.2) {\scriptsize $i+1$};
      }[][(0,0)]
    }
    \\[1ex]
    \IEEEeqnarraymulticol{3}{c}{
      \text{\small bubble evaluations ($\bb$), either with external shading ($\bbev$) or internal shading ($\bbodd$)}
    }
    \\[2ex]
    \begingroup
      \tikzpic{
        \flst[0][1]\flst\node[left=-3pt] at (0,2) {\scriptsize $i$};
        \frst[1][1]\frst[1][0];
      }[scale=.7][(0,1*.7)]
      \overset{}{\Rrightarrow}
      \tikzpic{
        \fdot[.5][1.8]\fcup[0][1]\fcap
        \node[left=-3pt] at (0,2) {\scriptsize $i$};
      }[scale=.7][(0,1*.7)]
      \;+\;
      \tikzpic{
        \fdot[.5][.2]\fcup[0][1]\fcap
        \node[left=-3pt] at (0,2) {\scriptsize $i$};
      }[scale=.7][(0,1*.7)]
    \endgroup
    & \mspace{120mu} &
    \begingroup
      \tikzpic{
        \flst[0][0]\frst[.5][0][2]\flst[1.5][0][2]\frst[2][0]
        \node[below=-2pt] at (0,0) {\scriptsize $i$};
        \node[below=-2pt] at (1.5,0) {\scriptsize $i+1$};
      }[yscale=4,scale=.35][(0,2*.35)]
      \Rrightarrow Z^{-1}\;
      \tikzpic{
        \clip (-1,0) rectangle (2,4);
        \fzip[0][3][2]
        \funzip[0][0][2]
        \begin{scope}[scale=2]
          \fcup[-.25][1]
          \fcap[-.25][0]
        \end{scope}
        \node[below=-2pt] at (-.5,0) {\scriptsize $i$};
        \node[below=-2pt] at (1,0) {\scriptsize $i+1$};
      }[scale=.35][(0,2*.35)]
    \endgroup
    \\[1ex]
    \text{\small neck-cutting relation ($\nc$)}
    &&
    \text{\small squeezing relation ($\sq$)}
  \end{IEEEeqnarray*}
  
  \caption{%
    The set of 3-cells $\foamR_3=\{\dd,\dm,\bb,\nc,\sq\}$.
  }
  \label{fig:foam_R}
\end{figure}

\begin{definition}
  \label{defn:graded2cat_foams}
  The graded-2-category $\gfoam_d$ is the graded-2-category presented by $\sgfoam_d$:
\[\gfoam_d \coloneqq [\sgfoam_d]_\sim.\]
\end{definition}

\subsubsection{Foams isotopies}
\label{subsubsec:foam_isotopies}

The 3-cells in $\foamE_3$ encode a braid-like pivotal structure in $\gfoam_d$.
This will be our modulo data (\cref{subsec:rewriting_foam_coherence_iso}).
Topologically, these relations encode certain isotopies of singular surfaces called \emph{foam isotopies}.

\begin{terminology}
  \label{term:foam_isotopy}
  A 3-cell in $\foamE^\top$ is called a \emph{foam isotopy}, or simply an \emph{isotopy}; accordingly, we shall say that two diagrams are \emph{isotopic} if they are the source and target of a 3-cell in $\foamE^\top$.
\end{terminology}

The data of dots and strands in a given diagram is preserved under foam isotopies; that is, if $e\colon D_0\to D_1$ is an isotopy, then there is a canonical bijection between the set of dots and strands of $D_0$, with the set of dots and strands of $D_1$.
The following statement is shown in \cref{subsec:rewriting_foam_coherence_iso}:

\begin{proposition}[coherence of foam isotopies]
  \label{prop:coherence_foam_isotopies}
  If two parallel morphisms in $\foamE^\top$ define the same bijection on dots and strands, then they have the same associated scalar.
\end{proposition}

\begin{terminology}
  \label{term:dot_isotopy}
  Following on \cref{term:foam_isotopy}, we say that two dots in a diagram are \emph{isotopic} if there exists an isotopic diagram in which the two dots are next to each other.
  Similarly, a dot and a strand, or two strands, can be \emph{isotopic}.
\end{terminology}

\subsubsection{The basis theorem}

We review the notion of a reduced family (\cite{SV_OddKhovanovHomology_2023}).
As it will not be essential for us, we only give a quick definition, and refer to \cite[Section~2.5]{SV_OddKhovanovHomology_2023} for details.

Recall that our diagrams can be understood as foams. Given a foam $F\colon W\to W'$, write $\undcomp(F)$ the surface obtained by removing its 2\nbd-facets\footnote{More precisely, by removing the interior of the 2-facets, leaving 1-facets glued along seams. Since precisely two 1-facets meet at each seam, this gives a topological surface with boundary.}.
Similarly, given two parallel 1-morphisms $W$ and $W'$ in $\gfoam_d$, one can associate a certain union of circles $S^1$ denoted $\undcomp(W\sqcup_\partial W')$.
Write $\pi_0(\undcomp(W\sqcup_\partial W'))$ the set of connected components in $\undcomp(W\sqcup_\partial W')$.

A foam $F\colon W\to W'$ is \emph{reduced} if $\undcomp(F)$ is a union of disks, with at most one dot on each disk.
In terms of shading diagrammatics, ``reduced'' corresponds to the following two properties:
\begin{itemize}
  \item any path between two distinct $i$-strands intersect with an $i\pm 1$-strand;
  \item no two dots are related by a sequence of dot slides and dot migrations.
\end{itemize}
Up to foam isotopies, a reduced foam is characterized by a subset $\delta\subset \pi_0(\partial \undcomp(W\sqcup_\partial W'))$: a dot lies on a disk $D$ if and only if $\partial D$ belongs to $\delta$.
In that case, we say that $F$ is \emph{$\delta$-dotted}.

\begin{definition}
  \label{defn:reduced_family}
  Let $W$ and $W'$ be two parallel 1-morphisms in $\gfoam_d$.
  A \emph{reduced family} is a family $(F_\delta)_{\delta\subset \pi_0(\undcomp(W\sqcup_\partial W'))}$, where $F_\delta\colon W\to W'$ is a $\delta$-dotted reduced foam.
\end{definition}

The following is the main result of this section:

\begin{theorem}
  \label{thm:foam_basis_theorem}
  Let $W$ and $W'$ be two parallel 1-morphisms in $\gfoam_d$.
  Any reduced family defines a basis of the $\Bbbk$-module $\Hom_{\gfoam_d}(W,W')$.
\end{theorem}

\subsection{A convergent presentation of graded \texorpdfstring{$\glt$}{gl2}-foams}
\label{subsec:convergent_pres_gltfoam}

This section sets up the ``working data'' for the rewriting theory of graded $\glt$\nbd-foams. This is the data (a), (b) and (c) as described in \cref{subsubsec:HLRSM_setup}.

\subsubsection{A \texorpdfstring{\LGRSM{}}{linear Gray RSM} for graded \texorpdfstring{$\glt$}{gl2}-foams}
\label{subsubsec:GRSM_graded_foam}

Our rewriting approach to \cref{thm:foam_basis_theorem} is based on the following \LGRSM{}:
\[\foamS\coloneqq(\foamR,\foamE).\]
We further define $\foamB$, a linear sub-3\nbd-sesqui\-poly\-graph of $\foamR$, with $\foamB_{\leq 2}\coloneqq\foamR_{\leq 2}$ and
\begin{gather*}
  \foamB_3\coloneqq\{\dd,\dm,\bb\},
\end{gather*}
where $\bb=\{\bbev,\bbodd\}$. We write $\foamlR\coloneqq\Cont(\foamR)$, $\foamlE\coloneqq\Cont(\foamE)$ and $\foamlB\coloneqq\Cont(\foamB)$ the associated (family of) \LRS{}, and $\foamlS\coloneqq(\foamlR,\foamlE)$ the associated (family of) \LRSM{}.

Note that:

\begin{lemma}
  \label{lem:foam_other_rewriting_steps}
  The following are $\foamS$-rewriting steps:
  \begin{gather*}
    \tikzpic{
      \flst\node[below=-2pt] at (0,0) {\scriptsize $i$};
      \fdot[.5][.5]\node[below right] at (.5,.5) {\scriptsize $i$};
    }[xscale=.5][(0,.5)]
    \;\overset{}{\Rrightarrow}\;\norma{}{-\;\;}{}
    \tikzpic{
      \flst\node[below=-2pt] at (0,0) {\scriptsize $i$};
      \fdot[.5][.5][2]\node[below right] at (.5,.5) {\scriptsize $i+1$};
    }[xscale=.5][(0,.5)]
    \mspace{50mu}
    \begingroup
    \tikzpic{
      \clip (-1,0) rectangle (2,4);
      \fzip[0][3][1]
      \funzip[0][0][1]
      \begin{scope}[scale=2]
        \fcup[-.25][1][2]
        \fcap[-.25][0][2]
      \end{scope}
      \node[below=-2pt] at (-.5,0) {\scriptsize $i+1$};
      \node[below=-2pt] at (1,0) {\scriptsize $i$};
    }[scale=.35][(0,2*.35)]
    \Rrightarrow XYZ
    \tikzpic{
      \flst[0][0][2]\frst[.5][0][1]\flst[1.5][0][1]\frst[2][0][2]
      \node[below=-2pt] at (0,0) {\scriptsize $i+1$};
      \node[below=-2pt] at (1.5,0) {\scriptsize $i$};
    }[yscale=4,scale=.35][(0,2*.35)]
    \endgroup
    \\[1ex]
    \tikzpic{
      \funzip\fzip[0][1]\node[below=-2pt] at (0,0) {\scriptsize $i$};
    }[scale=.5][(0,.5)]
    \;\overset{}{\Rrightarrow}\;
    XZ^2\;
    \tikzpic{
      \frst\flst[1][0]\node[below=-2pt] at (0,0) {\scriptsize $i$};
      \fdot[1.5][.5]\node[below right=-2pt] at (1.5,.5) {\scriptsize $i$};
    }[yscale=2,scale=.5][(0,.5)]
    \;+\norma{YZ}{YZ}{YZ^2}\;
    \tikzpic{
      \frst\flst[1][0]\node[below=-2pt] at (0,0) {\scriptsize $i$};
      \fdot[-.7][.5]\node[below right=-2pt] at (-.7,.5) {\scriptsize $i$};
    }[yscale=2,scale=.5][(0,.5)]
  \end{gather*}
  which we denote $\ov{\dm}$, $\ov{sq}$ and $\ov{\nc}$, respectively.
\end{lemma}

\subsubsection{A context-dependent linear sub-system}
\label{subsubsec:rw_foam_context_sub_syst}

The \LGRSM{} $\foamS$ is not terminating, and hence not suited for a reduction algorithm. For instance, we have the following infinite sequence:
\begin{gather*}
  \tikzpic{
    \fcap[0][1]
    \flst\frst[1][0]
  }[scale=.8]
  \;\overset{\nc}{\Rrightarrow}\;
  \tikzpic{
    \fcap[0][2]\fcup[0][1]\fcap
    \fdot[.5][2-.2]
  }[scale=.5]
  \;+\;
  \tikzpic{
    \fcap[0][2]\fcup[0][1]\fcap
    \fdot[.5][.2]
  }[scale=.5]
  \;\overset{\nc}{\Rrightarrow}\;
  \tikzpic{
    \fcap[0][4]\fcup[0][3]\fdot[.5][4-.2]
    \fcap[0][2]\fcup[0][1]\fcap
    \fdot[.5][2-.2]
  }[scale=.5,yscale=.7]
  \;+\;
  \tikzpic{
    \fcap[0][4]\fcup[0][3]\fdot[.5][2+.2]
    \fcap[0][2]\fcup[0][1]\fcap
    \fdot[.5][.2]
  }[scale=.5,yscale=.7]
  \;+\;
  \tikzpic{
    \fcap[0][2]\fcup[0][1]\fcap
    \fdot[.5][.2]
  }[scale=.5]
  \;\overset{\nc}{\Rrightarrow}\;
  \ldots
\end{gather*}
To avoid this obstruction to termination, we define a linear sub-system of $\foamS$ that prevents one from applying a neck-cutting relation ($\nc$) on two pieces of the same strand.
Taking into account the analogue obstruction for the squeezing relation ($\sq$) leads to the following definition:

\begin{definition}
  The (family of) \LRSM{}(s) $\foamlT$ is the linear sub-system $\foamlT\subset\foamS$ such that for each $\ssr\in\foamR_3$, $\Gamma[\ssr]\in\foamlT$ if and only if either $\ssr\in\foamB$, or $\ssr=\nc$ (resp.\ $\sq$) of colour $i$ (resp.\ $(i,i+1)$), and the two pieces of $i$-strands belong to distinct strands in $s(\Gamma[\ssr])$.
\end{definition}

It is not hard to check that $\foamlT$ is adapted, and hence is indeed a linear sub-system in the sense of \cref{defn:HLRSM_subsystem}; in fact, $\foamS$ is already adapted. Note however that $\foamlT$ is \emph{not} of the form $\Cont(\foamT)$ for some sub-\GRSM{} $\foamT\subset\foamS$.

The following lemma shows that $\foamlT$ is a suitable choice to prove \cref{thm:foam_basis_theorem}:

\begin{lemma}
  \label{lem:normal_form_equal_reduced}
  A foam is a $\foamlT$-normal form if and only if it is reduced.
  In particular, any choice of foam isotopy representatives for $\cB\NF_{\foamlT}$ is a reduced family.
\end{lemma}

The proof of \cref{lem:normal_form_equal_reduced} is given in \cite[Proposition~1.6.5]{Schelstraete_OddKhovanovHomology_2024}; in order to keep the focus on the rewriting theory, we omit it here.

\medbreak

Given a foam $F$, we associate to $F$ the following data (recall \cref{term:i_strand}):
\begin{align*}
  \#\mathrm{sh}_i(F) &\coloneqq \text{number of $i$-shadings in $F$},\\
  \#\mathrm{cl}_i(F) &\coloneqq \text{number of closed $i$-strands in $F$},\\
  \#\mathrm{d}_i(F) &\coloneqq \text{number of $i$-dots in $F$.}
\end{align*}
Each of these data defines an esp-order on foams. Let $\foamsucc$ be the lexicographic esp-order induced by these esp-orders:
\[{\foamsucc}\coloneqq(\#\mathrm{sh}_1,\ldots,\#\mathrm{sh}_{d-1},\#\mathrm{cl}_1,\ldots,\#\mathrm{cl}_{d-1},\#\mathrm{d}_1,\ldots,\#\mathrm{d}_d).\]
Since $\#\mathrm{sh}_i(F)$, $\#\mathrm{cl}_i(F)$ and $\#\mathrm{d}_i(F)$ are valued in the natural numbers, $\foamsucc$ is well-founded.

\begin{lemma}
  \label{prop:foam_T_terminates}
  The esp-order $\foamsucc$ is strongly compatible (\cref{defn:strong_compatibility}) with $\foamlT$.
\end{lemma}

\begin{proof}
  For rewriting steps of type $\foamB=\{\dd,\dm,\bb\}$, this is fairly straightforward.
  On the other hand, a rewriting step $\Gamma[\ssr]$ of type $\nc$ or $\sq$ does not change the number of $k$-shadings for $k<i$ and reduces the number of $i$-shadings if and only if the two pieces of $i$-strand do not belong to the same $i$-strand in $s(\Gamma[\ssr])$; this is precisely the definition of $\foamlT$.
\end{proof}

\subsubsection{A basis for graded \texorpdfstring{$\glt$}{gl2}-foams}
\label{subsubsec:rw_foam_basis_statement}

The following is a corollary of the coherence of foam isotopies (\cref{prop:coherence_foam_isotopies}):

\begin{corollary}
  \label{cor:Morse-regular_iso_coherent_on_reduced_foams}
  $\foamE^\top$ is scalar-coherent on reduced foams.
\end{corollary}

\begin{proof}
  A reduced foam, expressed in diagrammatics, has no closed strands (see \cite[Corollary~1.6.8]{Schelstraete_OddKhovanovHomology_2024}), nor isotopic dots with the same colour. Hence, foam isotopies on reduced foams all induce the same bijection on dots and strands.
  The statement then follows from \cref{prop:coherence_foam_isotopies}.
\end{proof}

To show \cref{thm:foam_basis_theorem} using the \LRSMbasisfromconvergencetheorem{}, it remains to show the coherence of foam isotopies (\cref{prop:coherence_foam_isotopies})---this is done in \cref{subsec:rewriting_foam_coherence_iso}---and to check that, on one hand:

\begin{lemma}
  \label{lem:EquivTandSfoam}
  The family of \LRSM{} $\foamlT$ and $\foamlS$ present the same underlying family of modules.
\end{lemma}

and on the other hand, thanks to the \LRSMtamednewmannlemma{}:

\begin{proposition}
  \label{prop:foamTmonomialtamecongruence}
  Every monomial $\foamlT$-branching is $\foamsucc$-tamely $\foamlT^\st$\nbd-con\-gruent.
\end{proposition}

Both statements are shown in \cref{subsec:rewriting_foam_confluence_modulo_iso}.

\subsection{Coherence of foam isotopies}
\label{subsec:rewriting_foam_coherence_iso}

Recall the terminologies introduced in \cref{term:i_strand,term:foam_isotopy,term:dot_isotopy}. Furthermore:

\begin{terminology}
  \label{term:distant_colour}
  Two colours $i$ and $j$ are said to be \emph{distant} if $\abs{i-j}>1$.
  Given a colour $i$ and a diagram $\psi$, we say that \emph{$i$ is distant from $\psi$} if for each $j$-strand (resp.\ $j$-dot) in $\psi$, we have $\abs{i-j}>1$ (resp.\ $j\neq i,i+1$).
\end{terminology}

Recall the scalar 3\nbd-sesqui\-po\-ly\-graphs defined in \cref{subsubsec:foam_isotopies}.
Let
\[\foamoE=\foamGray\sqcup\foamX.\]
In this subsection, we prove the coherence of foam isotopies (\cref{prop:coherence_foam_isotopies}) using the higher Newman's lemma (\cref{lem:HRSM_higher_newmann_lemma}) for the \GRSM{} $(\foamZ,\foamoE)$. The strategy of proof is close to the one described for superadjunction in \cref{subsec:intro_example}.

\medbreak

We start by providing the analogous coherence result for the modulo data:

\begin{lemma}[coherence of foam isotopies, except zigzag relations]
  \label{lem:coherence_foam_iso_without_zigzag}
  If two parallel morphisms in $\foamoE^\top$ define the same bijection on dots and strands, then they have the same associated scalar.
\end{lemma}

\begin{proof}
   Since $\bilfoam$ is symmetric, the scalar of a relation in $\foamGray^\top$ only depends on how generating 2-cells vertically permute (see \cite[Section~2.1]{SV_OddKhovanovHomology_2023}, and also the coherence of interchangers in Gray categories \cite[Appendix~A]{Schelstraete_OddKhovanovHomology_2024}). In fact, it only depends on how generating 2-cells \emph{with a non-trivial grading} vertically permute.
   3-cells in $\foamX^\top$ have trivial associated scalar. While they do not preserve the data of the set of generating 2-cells, they do preserve the data of the set of \emph{non-trivially graded} generating 2-cells. That is, if $e\colon D_0\to D_1$ is a 3-cell in $\foamX$, there is a canonical way to identify the non-trivially graded generating 2-cells of $D_0$ with those of $D_1$.
   In other words, the permutation of non-trivially graded generating 2-cells is a well-defined data associated to any 3-cell in $\foamoE^\top$, and this data determines the associated scalar.
   Finally, we check that this permutation data is itself determined by the bijection on dots and strands.
\end{proof}

The modulo $\foamoE^\top$ encompasses a braided-like structure.
Indeed, $\foamoE^\top$ is equivalently generated by the braided-like relations and the following 3\nbd-morphisms:
\begin{IEEEeqnarray*}{CcCcC}
  \tikzpic{
    \pic[transform shape] at (1,0) {cap=diag2};
    \draw[diag1] (0,0) to[out=90,in=-90] (2,3);
    \node[below=-1pt] at (0,0) {\scriptsize $i$};
    \node[below=-1pt] at (2,0) {\scriptsize $j$};
  }[xscale=-.5,yscale=.6,scale=.7]
  \Rrightarrow
  \tikzpic{
    \pic[transform shape] at (0,0) {crossing=diag1/diag2};
    \pic[transform shape] at (1,1) {crossing=diag1/diag2};
    \pic[transform shape] at (0,2) {cap=diag2};
    \draw[diag2] (0,1) to (0,2);
    \draw[diag2] (2,0) to (2,1);
    \draw[diag1] (2,2) to (2,3);
    \node[below=-1pt] at (0,0) {\scriptsize $i$};
    \node[below=-1pt] at (2,0) {\scriptsize $j$};
  }[xscale=-.5,yscale=.6,scale=.7]
  &\mspace{70mu}&
  \tikzpic{
    \pic[transform shape] at (0,2) {crossing=diag2/diag1};
    \pic[transform shape] at (1,1) {crossing=diag2/diag1};
    \pic[transform shape] at (0,0) {cup=diag2};
    \draw[diag2] (0,1) to (0,2);
    \draw[diag1] (2,0) to (2,1);
    \draw[diag2] (2,2) to (2,3);
    \node[above=-4pt] at (0,0) {\scriptsize $j$};
    \node[below=-1pt] at (2,0) {\scriptsize $i$};
  }[xscale=.5,yscale=.6,scale=.7]
  \Rrightarrow
  \tikzpic{
    \pic[transform shape] at (1,2) {cup=diag2};
    \draw[diag1] (2,0) to[out=90,in=-90] (0,3);
    \node[below=-1pt] at (2,2.5) {\scriptsize $j$};
    \node[below=-1pt] at (2,0) {\scriptsize $i$};
  }[xscale=.5,yscale=.6,scale=.7]
  &\mspace{70mu}&
  \tikzpic{
    \node[fdot2] at (.5,.5) {};
    \draw[diag1] (2,0) to[out=90,in=-90] (0,3);
    \node[above=-2pt] at (0,0) {\scriptsize $j$};
    \node[below=-1pt] at (2,0) {\scriptsize $i$};
  }[xscale=.5,yscale=.6,scale=.7]
  \Rrightarrow
  \tikzpic{
    \node[fdot2] at (1.5,2.5) {};
    \draw[diag1] (2,0) to[out=90,in=-90] (0,3);
    \node[above=-2pt] at (2,2) {\scriptsize $j$};
    \node[below=-1pt] at (2,0) {\scriptsize $i$};
  }[xscale=.5,yscale=.6,scale=.7]
  \\*
  \text{if $\abs{i-j}>1$}
  &&
  \text{if $\abs{i-j}>1$}
  &&
  \text{if $j\neq i,i+1$}
\end{IEEEeqnarray*}
which capture the naturality of the braiding with respect to the cap, the cup and the dot.
We say ``braided-like'' to emphasize the restriction on the labels.

Recall that any \GRSM{} has canonical interchange naturalities (\cref{subsubsec:scalar_HRSM}). 
For $(\foamZ,\foamoE)$, we further have the following braided-like naturalities:

\begin{lemma}[braided-like $\foamoE$-naturalities]
  \label{lem:regular_foam_braid_naturalities}
  For any relation $f\colon\psi_0\to\psi_1$ in $\foamZ$ and $i$ any colour distant from $\psi_0$ (see \cref{term:distant_colour}), we have the following $\foamoE$-congruences:
  \begin{center}
    \hfill
    \begin{tikzcd}[ampersand replacement=\&,anchor=south]
      \satex{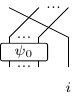}
      \tar[r,"{\satex[scale=.7]{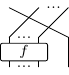}}"{yshift=1pt}]\tar[d,snakecd]
      \&
      \satex{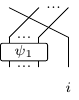}
      \tar[d,snakecd]
      \\
      \satex{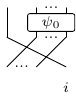}
      \tar[r,"{\satex[scale=.7]{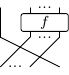}}"{yshift=1pt}]
      \&
      \satex{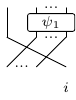}
    \end{tikzcd}
    \hspace{1cm}
    \begin{tikzcd}[ampersand replacement=\&,anchor=south]
      \satex{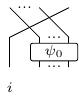}
      \tar[r,"{\satex[scale=.7]{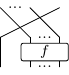}}"{yshift=1pt}]\tar[d,snakecd]
      \&
      \satex{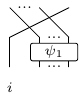}
      \tar[d,snakecd]
      \\
      \satex{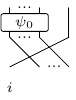}
      \tar[r,"{\satex[scale=.7]{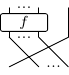}}"{yshift=1pt}]
      \&
      \satex{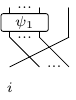}
    \end{tikzcd}
    \hfill\qed
  \end{center}
\end{lemma}

Our next step is to upgrade our knowledge of $\foamoE$-naturalities to a general characterization of $\foamoE$-congruence classes, leveraging our understanding of coherence in $\foamoE$:

\begin{lemma}[characterization of $\foamoE$-congruence classes for $\foamZ$]
  \label{lem:characterization_zigzag}
  If $(f,g)$ is a local $\foamZ$-branching such that $f$ and $g$ are of the same type with identical cup and cap, then $(f,g)$ is $\foamoE$-scalar-congruent.
\end{lemma}

\begin{proof}
  Let $[f,e,g]$ be a $\foamZ$-local triple with $f$ and $g$ as in the statement of the lemma. Thanks to coherence of the modulo (\cref{lem:coherence_foam_iso_without_zigzag}), we can choose the isotopy $e\colon s(f)\to s(g)$ such that each step either does not overlap $s(f)$, or consists of a braided-like $\foamoE$-naturality (\cref{lem:regular_foam_braid_naturalities}).
  This concludes.
\end{proof}

With this characterization at hand, we can classify local $\foamZ$-branchings and reduce scalar-confluence to four explicit branchings, the critical branchings:

\begin{lemma}
  \label{lem:confluence_zigzag}
  Every local $\foamZ$-branching is, up to branchwise $\foamoE$-congruence, either an independent branching or a contextualization of one of the following scalar-confluent branchings:
  \begin{gather*}
    \begin{tikzpicture}
      \node (A) at (0,0) {
        \tikzpic{
          \flst[0][2]\funzip[1][2]\frst[3][2]
          \fcup[0][1]\flst[2][1]\frst[3][1]
          \fcup[2][0]
          \fill[foamshade1] (-.5,0) rectangle (1.5,1);
        }[scale=.5,xscale=.8,yscale=.7]
      };
      \node (B) at (3,0) {
        \tikzpic{
          \fcup
        }[scale=.5,xscale=.8]
      };
      \draw[->] (A) to[out=20,in=160] node[above]{$YZ^2$} (B);
      \draw[->] (A) to[out=-20,in=-160] node[below]{$Y\cdot Z^2$} (B);
    \end{tikzpicture}
    \mspace{50mu}
    \begin{tikzpicture}
      \node (A) at (0,0) {
        \tikzpic{
          \funzip[2][2]
          \funzip[0][1]\frst[2][1]\flst[3][1]
          \frst[0][0]\fcup[1][0]\flst[3][0]
        }[scale=.5,xscale=.8,yscale=.7]
      };
      \node (B) at (3,0) {
        \tikzpic{
          \funzip
        }[scale=.5,xscale=.8]
      };
      \draw[->] (A) to[out=20,in=160] node[above]{$Z^2$} (B);
      \draw[->] (A) to[out=-20,in=-160] node[below]{$Y\cdot YZ^2$} (B);
    \end{tikzpicture}
    \\[1ex]
    \begin{tikzpicture}
      \node (A) at (0,0) {
        \tikzpic{
          \frst[0][2]\fcap[1][2]\flst[3][2]
          \fzip[0][1]\frst[2][1]\flst[3][1]
          \fzip[2][0]
        }[scale=.5,xscale=.8,yscale=.7]
      };
      \node (B) at (3,0) {
        \tikzpic{
          \fzip
        }[scale=.5,xscale=.8]
      };
      \draw[->] (A) to[out=20,in=160] node[above]{$X$} (B);
      \draw[->] (A) to[out=-20,in=-160] node[below]{$X\cdot 1$} (B);
    \end{tikzpicture}
    \mspace{50mu}
    \begin{tikzpicture}
      \node (A) at (0,0) {
        \tikzpic{
          \fill[foamshade1] (-.5,2) rectangle (1.5,3);\fcap[2][2]
          \fcap[0][1]\flst[2][1]\frst[3][1]
          \flst[0][0]\fzip[1][0]\frst[3][0]
        }[scale=.5,xscale=.8,yscale=.7]
      };
      \node (B) at (3,0) {
        \tikzpic{
          \fcap
        }[scale=.5,xscale=.8]
      };
      \draw[->] (A) to[out=20,in=160] node[above]{$1$} (B);
      \draw[->] (A) to[out=-20,in=-160] node[below]{$X\cdot X$} (B);
    \end{tikzpicture}
  \end{gather*}
\end{lemma}

\begin{proof}
  Let $[f,e,g]=[\Gamma_f[\ssr_f],e,\Gamma_f[\ssr_g]]$ be a $\foamZ$-local triple, with $\ssr_f,\ssr_g\in\foamZ_3$. We wish to simplify $[f,e,g]$ by replacing it with a branchwise $\foamoE$-congruent branching using the characterization given in \cref{lem:characterization_zigzag}.
  Given a cup $\cup$ and a cap $\cap$ belonging to the same strand $S$, we refer to a ``strand between the $\cup$ and $\cap$'' to mean a strand crossing $S$ on the piece of strand connecting $\cup$ to $\cap$.

  Consider the cup and cap of $s(\ssr_g)$, respectively denoted $\cup_g$ and $\cap_g$, as sitting in $s(f)$. They necessarily belong to the same strand $S$, possibly with strands between $\cup_g$ and $\cap_g$. Using isotopies in $\foamoE$, we can slide these strands away, so that no strand lies between $\cup_g$ and $\cap_g$.
  Moreover, we can do so without adding strands between $\cup_f$ and $\cap_f$, the cup and cap of $s(\ssr_f)$ (viewed as sitting in $s(f)$). Thanks to the characterization given in \cref{lem:characterization_zigzag}, this does not change the branchwise $\foamoE$-congruence class of $[f,e,g]$.

  If the cups and caps of $s(\ssr_f)$ and $s(\ssr_g)$ are disjoint, we can further use interchanges to move the cup and cap of $s(\ssr_g)$ one below the other, again without affecting $s(\ssr_f)$.
  This shows that $[f,e,g]$ is branchwise $\foamoE$-congruent to an independent $\foamZ$-branching.
  If $s(\ssr_f)$ and $s(\ssr_g)$ have precisely one cap or one cup in common, we find the four critical branchings given in the lemma.
\end{proof}

Moreover:

\begin{lemma}
  \label{lem:termination_zigzag}
  $\foamZ$ terminates modulo $\ov{\foamE}$.
\end{lemma}

\begin{proof}
  $\foamZ$ strictly diminishes the number of caps and cups, which is kept constant by $\foamoE$.
\end{proof}

We now have all the ingredients to show \cref{prop:coherence_foam_isotopies}:

\begin{proof}[Proof of \cref{prop:coherence_foam_isotopies}]
  Using the higher Newman's lemma (\cref{lem:HRSM_higher_newmann_lemma}), the analysis of local branchings (\cref{lem:confluence_zigzag}) and termination (\cref{lem:termination_zigzag}), we conclude that $\foamZ$ is convergent modulo $\foamoE$.
  Hence, $\foamZ^\top$ is scalar-coherent modulo $\foamoE^\top$ (coherence modulo from convergence modulo; \cref{prop:coherence_from_convergence}).
  
  To show the proposition, it suffices to show that if a loop $e$ in $\foamE^\top$ defines the identity bijection of dots and strands, then its associated scalar is one.
  Since $\foamZ^\top$ is scalar-coherent modulo $\foamoE^\top$, we can write $e$ as $e=g\circ h\circ g^{-1}$ with $h$ in $\foamoE^\top$. We conclude using coherence of $\foamoE^\top$ (\cref{lem:coherence_foam_iso_without_zigzag}).
\end{proof}

\subsection{Analysis of monomial branchings}
\label{subsec:rewriting_foam_confluence_modulo_iso}

\begin{notation}
  Let $\foamlA\coloneqq\foamlT\setminus\foamlB$; that is, $\foamlA$ consists of the rewriting step of type $\nc$ and $\sq$ which are in $\foamlT$.
\end{notation}

We study monomial $\lT$-branchings.
It is divided in three parts: both branches are of type $\foamlB$ (\cref{prop:confluence_modulo_B_foam}), one branch is of type $\foamlB$ and the other is of type $\foamlA$ (\cref{lem:foam_monomial_confluence_A_B}), and both branches are of type $\foamlA$ (\cref{lem:foam_monomial_confluence_A_A}).
The general strategy (with some variations) is the one described in \cref{subsubsec:HLRSM_general_strategy}: give naturalities of the modulo, characterize rewriting steps modulo, and finally enumerate monomial branchings.

Throughout we use the terminology of \cref{term:i_strand,term:foam_isotopy,term:dot_isotopy,term:distant_colour}.

\subsubsection{Foam isotopy naturalities}

\begin{lemma}[braided-like $\foamE$-naturalities]
  \label{lem:full_foam_braided_naturalities}
  The relations in $\foamR_3\sqcup\{\ov{\dm},\ov{\nc},\ov{\sq}\}$ satisfy the same braided-like naturalities as pictured in \cref{lem:regular_foam_braid_naturalities}.
\end{lemma}

\begin{proof}
  This follows from the fact that if $i$ is distant from $\psi_0$, then $i$ is distant from $\psi_1$ (using the same notations as in \cref{lem:regular_foam_braid_naturalities}). This is straightforward in most cases; we only detail type $\dm$. Let $j$ be the colour of the migrating dot in $\psi_0$. By assumption, the associated $j$-strand crosses an $i$-strand, so that we have $\abs{i-j}>1$. In particular $j+1\neq i,i+1$ and the $(j+1)$-dot can still slide through the $i$-strand once it has migrated.
\end{proof}

For $f\colon\psi_0\to\psi_1$ a 3-cell in $\foamR_3\sqcup\{\ov{\dm},\ov{\nc},\ov{\sq}\}$, we write $\rotatebox[origin=c]{180}{$f$}\colon\rotatebox[origin=c]{180}{$\psi_0$}\to\rotatebox[origin=c]{180}{$\psi_1$}$ the 3-cell obtained by rotating each diagram by a half-turn. Most 3-cells rotate to themselves, except for $\dm$ (resp.\ $\ov{\dm}$), which rotate to $\ov{\dm}$ (resp.\ $\dm$).

\begin{lemma}[pivotal $\foamE$-naturalities]
  \label{lem:foam_pivotal_naturalities}
  Let $f\colon\psi_0\to\psi_1$ be a 3-cell in $\foamR_3\sqcup\{\ov{\dm},\ov{\nc},\ov{\sq}\}$.
  Then the following squares commute:
  \begin{gather*}
    \begin{tikzcd}[ampersand replacement=\&]
      \satex{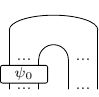}
      \tar[r,"{\satex[scale=.7]{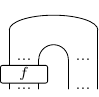}}"{yshift=1pt}]\tar[d,snakecd]
      \&
      \satex{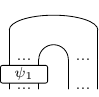}
      \tar[d,snakecd]
      \\
      \satex{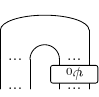}
      \tar[r,"{\satex[scale=.7]{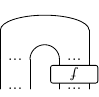}}"{yshift=1pt}]
      \&
      \satex{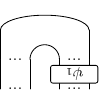}
    \end{tikzcd}
    \mspace{80mu}
    \begin{tikzcd}[ampersand replacement=\&]
      \satex{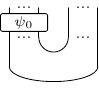}
      \tar[r,"{\satex[scale=.7]{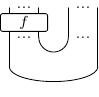}}"{yshift=1pt}]\tar[d,snakecd]
      \&
      \satex{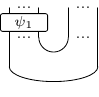}
      \tar[d,snakecd]
      \\
      \satex{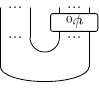}
      \tar[r,"{\satex[scale=.7]{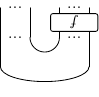}}"{yshift=1pt}]
      \&
      \satex{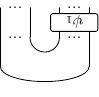}
    \end{tikzcd}
  \end{gather*}
\end{lemma}

\begin{proof}
  The statement is trivial for types $\dd$ and $\bb$, and follows from graded interchange for types $\dm$.
  Note that:
  \renewcommand{\scl}{.7}
  \def\tempdiag{
    \clip (-.5,0) rectangle (3.5,4);
    \fcap\fcup[0][1]
    \begin{scope}[yscale=2]
      \flst[2][0]\frst[3][0]
    \end{scope}
    \funzip[1][2]
    \begin{scope}[xscale=3,yscale=2]
      \fcap[0/3][2/2]
    \end{scope}
  }
  \def\tempdiagbis{
    \begin{scope}
      \clip (-1,0) rectangle (2,4);
      \fzip[0][3][2]
      \funzip[0][0][2]
      \begin{scope}[scale=2]
        \fcup[-.25][1]
        \fcap[-.25][0]
      \end{scope}
    \end{scope}
    \begin{scope}[shift={(2.5,0)},yscale=4]
      \flst
      \frst[1][0][2]
      \flst[2][0][2]
      \frst[3][0]
    \end{scope}
    \funzip[1.5][4]
    \draw[foamdraw2] (1,4) to[out=90,in=180] (2.25,5.5) to[out=0,in=90] (3.5,4);
    \draw[foamdraw2] (0,4) to[out=90,in=180] (2.25,6.5) to[out=0,in=90] (4.5,4);
    \fill[foamshade2] (1,4) to[out=90,in=180] (2.25,5.5) to[out=0,in=90] (3.5,4) to (4.5,4) to[out=90,in=0] (2.25,6.5) to[out=180,in=90] (0,4);
    \begin{scope}
      \clip (-1,0) rectangle (6,8);
      \begin{scope}[xscale=6,yscale=4.5]
        \fcap[-.5/6][4/4.5]
      \end{scope}
    \end{scope}
  }
  \begin{IEEEeqnarray*}{rClcrCl}
    \tikzpic{
      \tempdiag
    }[xscale=.6,yscale=.7,scale=.7*\scl]
    &\sim_\foamE&XY\;
    \tikzpic{
      \tempdiag
    }[xscale=-.6,yscale=.7,scale=.7*\scl]
    &\mspace{50mu}\an\mspace{50mu}&
    \tikzpic{
      \tempdiag
    }[xscale=.6,yscale=.7,scale=-.7*\scl]
    &\sim_\foamE&XY\;
    \tikzpic{
      \tempdiag
    }[xscale=-.6,yscale=.7,scale=-.7*\scl]\;.
  \end{IEEEeqnarray*}
  It follows that 
  \begin{gather*}
    \tikzpic{
      \tempdiag
      \fdot[.5][.2][1][1.5]
    }[xscale=.6,yscale=.7,scale=.7*\scl]
    \;+\;
    \tikzpic{
      \tempdiag
      \fdot[.5][2-.2][1][1.5]
    }[xscale=.6,yscale=.7,scale=.7*\scl]
    \;\sim_\foamE\;
    \tikzpic{
      \tempdiag
      \fdot[.5][2-.2][1][1.5]
    }[xscale=-.6,yscale=.7,scale=.7*\scl]
    \;+\;
    \tikzpic{
      \tempdiag
      \fdot[.5][.2][1][1.5]
    }[xscale=-.6,yscale=.7,scale=.7*\scl]
  \end{gather*}
  and similarly for the other case; this shows the lemma for type $\nc$. On the other hand, we have
  \renewcommand{\scl}{.5}
  \begin{IEEEeqnarray*}{rClcrCl}
    \tikzpic{
      \tempdiagbis
    }[xscale=.6,yscale=.7,scale=.6*\scl]
    &\;\sim_\foamE&
    \tikzpic{
      \tempdiagbis
    }[xscale=-.6,yscale=.7,scale=.6*\scl]
    &\mspace{50mu}\an\mspace{50mu}&
    \tikzpic{
      \tempdiagbis
    }[xscale=.6,yscale=.7,scale=-.6*\scl]
    &\sim_\foamE\;&
    \tikzpic{
      \tempdiagbis
    }[xscale=-.6,yscale=.7,scale=-.6*\scl]\;,
  \end{IEEEeqnarray*}
  which shows the lemma for type $\sq$.
  This also implies the lemma for types $\ov{\dm}$, $\ov{\nc}$ and $\ov{\sq}$, using coherence of $\foamE$.
\end{proof}

\subsubsection{Confluence of monomial \texorpdfstring{$\foamlB$}{B}-branchings}

We begin our study of confluence with positive branchings in the context-agnostic linear sub-system $\foamlB\subset\foamlT$, which derives from a \LGRSM{}: $\foamlB=\Cont(\foamB)$.
In fact, we show that in that case, confluences can also be taken in $\foamB$:

\begin{proposition}
  \label{prop:confluence_modulo_B_foam}
  Every monomial $\foamB$-branching is $\foamsucc$-tamely $\foamB$\nbd-cong\-ruent.
\end{proposition}

Since $\succ$ is strongly compatible with $\foamlT$ (\cref{prop:foam_T_terminates}), it is also strongly compatible with $\foamB$.
As a preliminary step, we give a topological characterization of $\foamE$-congruence classes:

\begin{lemma}[characterization of $\foamE$-congruence classes for $\foamB$]
  \label{lem:characterization_foam_B}
  Let $(f,g)$ be a monomial $\foamB$-branching, with $f$ and $g$ of the same type. The following holds:
  \begin{enumerate}[(a)]
    \item \emph{type $\dd$:} $(f,g)$ is $\foamE$\nbd-congruent;
    \item \emph{type $\dm$:} If $f$ and $g$ have isotopic $i$-dot and $i$-strand, then $(f,g)$ is $\foamE$\nbd-congruent;
    \item \emph{type $\bb$:} If $f$ and $g$ have isotopic $i$-strand, then $(f,g)$ is $\foamE$\nbd-congruent.
  \end{enumerate}
  We call the data associated to each type its \emph{combinatorial data}, and say that this combinatorial data \emph{characterizes} the $\foamE$-congruence class of the type.
\end{lemma}

\begin{proof}
  Denote by $\psi$ the local picture of the rewriting step $f$; that is, $\psi=s(\ssr)$ for $\ssr\in \foamB$ and $f=\Gamma[\ssr]$ for some context $\Gamma$.
  In each type, we use coherence of $\foamE$ (\cref{prop:coherence_foam_isotopies}) to present the isotopy $e$ as a composition of $\foamE$\nbd-natu\-ra\-lities, as described in \cref{lem:full_foam_braided_naturalities} and \cref{lem:foam_pivotal_naturalities}:
  \begin{enumerate}[(a)]
    \item Trivial, since both $f$ and $g$ rewrite to zero.
    \item Through the isotopy $e$, the $i$-dot starts and ends next to the $i$-strand. Hence, we can choose $e$ such that the $i$-dot always remains next to the $i$-strands. In that case, the only isotopies that overlap $\psi$ are braided-like naturalities.
    \item Given that the bubble $\psi$ starts and ends without any strand crossing it, we can choose the isotopy $e$ such that $D$ always crosses strands ``at once''. In that case, the only isotopies that overlap $\psi$ are braided-like isotopies.
  \end{enumerate}
  This concludes.
\end{proof}

As $\foamB$ is context-agnostic, the case of independent $\foamB$-branchings comes for free (\cref{lem:HLRSM_congruence_independent_branching}):

\begin{lemma}
  \label{lem:foam_independent_B_branching_tamed}
  Every independent $\foamB$-branchings is $\succ$-tamely $\foamB$\nbd-cong\-ruent.
  \hfill\qed
\end{lemma}

We can now prove the proposition:

\begin{proof}[Proof of \cref{prop:confluence_modulo_B_foam}]
  We use the \ARSMbranchwisetamedcongruencelemma{} and the fact that independent $\foamB$-branchings are $\succ$-tamely $\foamB$\nbd-cong\-ruent (\cref{lem:foam_independent_B_branching_tamed}) without further mention.

  Let $[f,e,g]$ be a monomial $\foamB$-local triple.
  If $f$ and $g$ are of type $\dd$ or $\bb$ (with possibly distinct types for $f$ and $g$), then either they have the same combinatorial data and hence are $\foamE$\nbd-congruent, or their combinatorial data are disjoint. We can isotope them away from each other, and $[f,e,g]$ is branchwise $\foamE$\nbd-congruent to an independent branching. A similar argument applies if $f$ is of type $\dm$ and $g$ is of type $\bbev$.

  Assume $f$ is of type $\dm$ and $g$ is of type $\dd$. If their combinatorial data are disjoint, $[f,e,g]$ is branchwise $\foamE$\nbd-congruent to an independent branching. Otherwise, they share the data of an $i$-dot.
  In that case, $[f,e,g]$ is $\foamE$\nbd-congruent to a contextualization of the following $\foamB$-branching, shown to be positively $\foamB$\nbd-conf\-luent:
  \begin{gather*}
    \def\scl{.6}
    \begin{tikzcd}[ampersand replacement=\&,row sep=0em,column sep=6em]
      \tikzpic{
        \begin{scope}[yscale=2]
          \frst
        \end{scope}
        \fdot[-.5][.5]
        \fdot[-.5][1.5]
      }[scale=.4]
      \&
      \tikzpic{
        \begin{scope}[yscale=2]
          \frst
        \end{scope}
        \fdot[-.5][.5]
        \fdot[-.5][1.5][2]
      }[scale=.4]
      \&
      \tikzpic{
        \begin{scope}[yscale=2]
          \frst
        \end{scope}
        \fdot[-.5][.5][2]
        \fdot[-.5][1.5][2]
      }[scale=.4]
      \&
      0
      \arrow[from=1-1,to=1-2,"\dm"]
      \arrow[from=1-2,to=1-3,"\dm"]
      \arrow[from=1-3,to=1-4,"\dd"]
      \arrow[from=1-1,to=1-4,"\dd",bend right=20]
    \end{tikzcd}
  \end{gather*}
  It follows from the \HLRSMcontextualizationlemma{} that $[f,e,g]$ is $\foamsucc$-tamely $\foamB$-confluent.

  Assume finally that $f$ is of type $\dm$ and $g$ is of type $\bbodd$.
  A similar reasoning reduces the statement to the following positive $\foamB$-confluence:
  \begin{gather*}
    \def\scl{.5}
    \begin{tikzcd}[ampersand replacement=\&,row sep=2em,column sep=1em]
      \&[-2.2em]
      \tikzpic{
        \fzip\funzip[0][1]
        \fdot[-.5][1]
      }[scale=.7*\scl]
      \& 
      \tikzpic{
        \fzip\funzip[0][1]
        \fdot[-.5][1][2]
      }[scale=.7*\scl]
      \&\&[-10em]
      YZ^{-1}\left(
      Z\;
      \tikzpic{
        \fdot[-.5][0][2]
        \fdot[.5][1][1]
      }[scale=.7*\scl]
      \;+\;
      XYZ\;
      \tikzpic{
        \fdot[-.5][0][2]
        \fdot[.5][1][2]
      }[scale=.7*\scl]
      \right)
      \&
      Y
      \&[-1.8em]
      \tikzpic{
        \fdot[-.5][0][2]
        \fdot[.5][1][1]
      }[scale=.7*\scl]
      \\
      YZ^{-1}\&[-3em]
      \tikzpic{
        \fzip\funzip[0][1]
        \fdot[-.5][0]
      }[scale=.7*\scl]
      \&\&
      YZ^{-1}\left(
      Z\;
      \tikzpic{
        \fdot[-.5][0][1]
        \fdot[.5][1][1]
      }[scale=.7*\scl]
      \;+\;
      XYZ\;
      \tikzpic{
        \fdot[-.5][0][1]
        \fdot[.5][1][2]
      }[scale=.7*\scl]
      \right)
      \&\&
      X
      \&[-1.8em]
      \tikzpic{
        \fdot[-.5][0][1]
        \fdot[.5][1][2]
      }[scale=.7*\scl]
      \arrow[from=1-2,to=1-3,"\dm"]
      \arrow[from=1-3,to=1-5,"\bbodd"{pos=.7},decorate,decoration={snake,amplitude=1.5pt,segment length=8pt,post length=1.5cm}]
      \arrow[from=1-5,to=1-6,"\dd"]
      \arrow[from=2-2,to=2-4,"\bbodd"]
      \arrow[from=2-4,to=2-6,"\dd"]
      \arrow[from=1-2,to=2-2,snakecd]
      \arrow[from=1-7,to=2-7,snakecd]
    \end{tikzcd}
  \end{gather*}
  This concludes.
\end{proof}

Finally, it follows from the \LRSMtamednewmannlemma{} that:

\begin{corollary}
  \label{cor:foamB_convergent}
  $\foamB^+$ is convergent.
  \hfill\qed
\end{corollary}

A \emph{bubble} is an endomorphism of an identity 1\nbd-morphism, possibly viewed inside a bigger diagram.
One can check (see \cite[Lemma~1.6.7]{Schelstraete_OddKhovanovHomology_2024}) that any bubble can be ``evaluated'' using $\foamB^+$\nbd-rewri\-ting steps, in the sense that it rewrites into a sum of diagrams, each consisting only of dots.
\Cref{prop:confluence_modulo_B_foam} shows that this evaluation is uniquely defined up to $\foamE$\nbd-congruence, so that we can speak of \emph{the} $\foamB^+$\nbd-rewri\-ting sequence evaluating a bubble:

\begin{definition}
  \label{defn:bubble_evalutation}
  For each bubble $\phi$, its \emph{bubble evaluation} is the $\foamB^+$\nbd-rewri\-ting sequence, well-defined up to $\foamE$-congruence, which rewrites $\phi$ into a sum of dots:
  \begin{gather*}
    \bb^*\colon\quad
    \tikzpic{
      \node[draw,circle,inner sep=1pt] at (0.5,2.2) {$\scriptstyle \phi$};
    }
    \overset{*}{\longrightarrow}_\foamB \sum \text{dots}.
  \end{gather*}
  A $\foamB$-rewriting sequence defined as a contextualized $\bb^*$ is said to be \emph{of type $\bb^*$}.
\end{definition}

\subsubsection{Characterization of neck-cutting and squeezing relation}

In order to study rewriting steps of type $\nc$ and $\sq$, we would like to give a topological characterization of their $\foamE$-congruence classes akin to the one given for $\foamB$ in \cref{lem:characterization_foam_B}.
In fact, it is easier to characterize their $\foamB$-congruence classes.

\begin{proposition}
  \label{prop:charactization_nc_sq}
  Let $(f,g)$ be a monomial
  $\foamS$-branching with $f$ and $g$ of type $\nc$ and label $i$ (resp. of type $\sq$ and label $(i,i+1)$). If $f$ and $g$ apply to the same $i$-strand(s), then $(f,g)$ is $\foamB$\nbd-congruent.
\end{proposition}

Note that since $\foamB^+$ is convergent (\cref{cor:foamB_convergent}), a branching is $\foamB$-congruent if and only if it is $\foamB^+$\nbd-confluent.
We call the data of the $i$-strand(s) the \emph{combinatorial data} of type $\nc$ (resp.\ $\sq$), and say that it \emph{characterizes} their $\foamB$-congruence class.
Before proving \cref{prop:charactization_nc_sq}, we show the following elementary $\foamB$-congruences:

\begin{lemma}[spatial-like $\foamB$-congruence]
  \label{lem:foam_spatial_B-confluence}
  The following branchings are $\foamB$\nbd-congruent, where $\phi$ is an arbitrary bubble and the dotted wiggly line denotes either a neck-cutting or a squeezing relation:
  \begin{IEEEeqnarray*}{CcC}
    \begin{tikzcd}[ampersand replacement=\&,column sep=2em]
      \tikzpic{
        \draw[,thin,snakecd,dottedcd] (0,1) to (1,1);
        \begin{scope}[yscale=2]
          \flst\frst[1][0]
        \end{scope}
        \node[draw,circle,inner sep=.5pt,anchor=center] at (0.5,1.5) {\scriptsize $\phi$};
      }[scale=.7]
      \&
      \tikzpic{
        \draw[,thin,snakecd,dottedcd] (0,1) to (1,1);
        \begin{scope}[yscale=2]
          \flst\frst[1][0]
        \end{scope}
        \node[draw,circle,inner sep=.5pt,anchor=center] at (0.5,.5) {\scriptsize $\phi$};
      }[scale=.7]
      \tar[from=1-1,to=1-2,snakecd]
    \end{tikzcd}
    &\qquad&
    \begin{tikzcd}[ampersand replacement=\&,column sep=2em]
      \tikzpic{
        \draw[,thin,snakecd,dottedcd] (0,1) to (3,1);
        \begin{scope}[yscale=2]
          \flst\frst[1][0][2]\flst[2][0][2]\frst[3][0]
        \end{scope}
        \node[draw,circle,inner sep=.5pt,anchor=center] at (.5,1.5) {\tiny $\phi$};
      }[xscale=.4,yscale=.6]
      \&
      \tikzpic{
        \draw[,thin,snakecd,dottedcd] (0,1) to (3,1);
        \begin{scope}[yscale=2]
          \flst\frst[1][0][2]\flst[2][0][2]\frst[3][0]
        \end{scope}
        \node[draw,circle,inner sep=.5pt,anchor=center] at (.5,.5) {\tiny $\phi$};
      }[xscale=.4,yscale=.6]
      \tar[from=1-1,to=1-2,snakecd]
    \end{tikzcd}
    \\
    \begin{tikzcd}[ampersand replacement=\&,column sep=2em]
      \tikzpic{
        \draw[,thin,snakecd,dottedcd] (0,1) to (3,1);
        \begin{scope}[yscale=2]
          \flst\frst[1][0][2]\flst[2][0][2]\frst[3][0]
        \end{scope}
        \node[draw,circle,inner sep=.5pt,anchor=center] at (1.5,1.5) {\tiny $\phi$};
      }[xscale=.4,yscale=.6]
      \&
      \tikzpic{
        \draw[,thin,snakecd,dottedcd] (0,1) to (3,1);
        \begin{scope}[yscale=2]
          \flst\frst[1][0][2]\flst[2][0][2]\frst[3][0]
        \end{scope}
        \node[draw,circle,inner sep=.5pt,anchor=center] at (1.5,.5) {\tiny $\phi$};
      }[xscale=.4,yscale=.6]
      \tar[from=1-1,to=1-2,snakecd]
    \end{tikzcd}
    &&
    \begin{tikzcd}[ampersand replacement=\&,column sep=2em]
      \tikzpic{
        \draw[,thin,snakecd,dottedcd] (0,1) to (3,1);
        \begin{scope}[yscale=2]
          \flst\frst[1][0][2]\flst[2][0][2]\frst[3][0]
        \end{scope}
        \node[draw,circle,inner sep=.5pt,anchor=center] at (2.5,1.5) {\tiny $\phi$};
      }[xscale=.4,yscale=.6]
      \&
      \tikzpic{
        \draw[,thin,snakecd,dottedcd] (0,1) to (3,1);
        \begin{scope}[yscale=2]
          \flst\frst[1][0][2]\flst[2][0][2]\frst[3][0]
        \end{scope}
        \node[draw,circle,inner sep=.5pt,anchor=center] at (2.5,.5) {\tiny $\phi$};
      }[xscale=.4,yscale=.6]
      \tar[from=1-1,to=1-2,snakecd]
    \end{tikzcd}
  \end{IEEEeqnarray*}
\end{lemma}

\begin{proof}
  We shall see that in each case, it suffices to evaluate the bubble $\phi$ (see \cref{defn:bubble_evalutation}) and apply some additional dot migrations to achieve $\foamB^+$-confluence.
  Denote
  \begin{gather*}
    \bb^*\colon\quad
    \tikzpic{
      \node[draw,circle,inner sep=1pt] at (0.5,2.2) {$\scriptstyle \phi$};
    }
    \overset{*}{\longrightarrow}_\foamB \sum_{\delta}\;
    \tikzpic{
      \node[draw,circle,inner sep=1pt] at (0.5,2.2) {$\scriptstyle \delta$};
    }
  \end{gather*}
  the bubble evaluation of $\phi$. Denote by $\delta$ a generic union of dots appearing in this bubble evaluation.
  
  Consider the first branching. We compare
  \begin{gather*}
    \tikzpic{
      \fcap\fcup[0][1]
      \fdot[.5][2-.2]
      \node[draw,circle,inner sep=1pt] at (0.5,2+.3) {$\scriptstyle \delta$};
    }[scale=.6][(0,1*.6)]
    \;+\;
    \tikzpic{
      \fcap\fcup[0][1]
      \fdot[.5][.2]
      \node[draw,circle,inner sep=1pt] at (0.5,2+.3) {$\scriptstyle \delta$};
    }[scale=.6][(0,1*.6)]
    \qquad\text{with}\qquad
    \tikzpic{
      \fcap\fcup[0][1]
      \fdot[.5][2-.2]
      \node[draw,circle,inner sep=1pt] at (0.5,-.3) {$\scriptstyle \delta$};
    }[scale=.6][(0,1*.6)]
    \;+\;
    \tikzpic{
      \fcap\fcup[0][1]
      \fdot[.5][.2]
      \node[draw,circle,inner sep=1pt] at (0.5,-.3) {$\scriptstyle \delta$};
    }[scale=.6][(0,1*.6)]
    \;.
  \end{gather*}
  If $\delta$ only consists of $j$-dots with $j\neq i,i+1$, then $\delta$ can slide from top to bottom, without additional scalar.
  If $\delta$ contains at least two $j$-dots with $j=i,i+1$, then both sides rewrite to zero, possibly migrating an $i$-dot into a $(i+1)$-dot first.
  Finally, if $\delta$ contains exactly one $j$-dot with $j=i,i+1$, then both sides rewrite to a single diagram consisting of a $(i+1)$-dot on top and a $(i+1)$-dot on the bottom (again using dot migrations).

  The other branchings are treated similarly. For the second and fourth branchings, one can use dot migrations to rewrite $i$- and $(i+1)$-dots into $(i+2)$-dots, allowing them to cross the $i$-strands. For the third branching, the fact that a dot cannot sit on a shaded region (for certain labels, see condition (iii) on page \pageref{page:label_condition}) prevents $(i+1)$- and $(i+2)$-dots, and $i$-dots can be treated as before, sliding them first across the $(i+1)$-strand.
\end{proof}

We can now prove the proposition:

\begin{proof}[Proof of \cref{prop:charactization_nc_sq}]
  Let $[f,e,g]$ be a local triple as in the proposition and denote $D_f$ and $D_g$ the source of $f$ and $g$, respectively.
  As in \cref{lem:characterization_foam_B}, denote $\psi_f$ the local picture of the rewriting step $f$; e.g.\ for type $\nc$, the diagram $\psi_f$ consists of two vertical $i$-strands.
  Denote $\psi_g$ similarly.
  The main idea is to treat $\psi_f$ and $\psi_g$ as formal generators, thanks to the following fact:
  \begin{quote}
    \textsc{Fact}. If $D_f$ and $D_g$ are isotopic, there exists an isotopy between $D_f$ and $D_g$ which isotope $\psi_f$ to $\psi_g$ and which is a composition of isotopies independent of $\psi_f$, braid-like isotopies with $\psi_f$, spatial-like isotopies with $\psi_f$ and possibly one pivotal isotopy with $\psi_f$.
  \end{quote}
  We leave the details to the reader; it is a variation of the standard decomposition of an isotopy of curves into elementary moves.
  Here is a simple example:
  \begin{gather*}
    \def\scl{.7}
    \def\tempdiag{
      \draw[thin,snakecd] (1,1) to (2,1);
      \frst[0][1]\flst[1][1]\funzip[2][1]
      \fzip[0][0]\frst[2][0]\flst[3][0]
    }
    \tikzpic{
      \tempdiag
    }[scale=.7*\scl]
    \;\sim_\foamE\;
    \tikzpic{
      \draw[thin,snakecd,] (1,1) to (2,1);
      \funzip[2][1]
      \fzip[0][0]\frst[2][0]\flst[3][0]\fzip[4][0]
      \begin{scope}
        \clip (.5,1) rectangle (4.5,3);
        \begin{scope}[xscale=3,yscale=2]
          \fcap[1/3][1/2]
        \end{scope}
      \end{scope}
      \begin{scope}[yscale=2]
        \frst[0][1/2]\flst[5][1/2]
      \end{scope}
    }[scale=.5*\scl]
    \;\sim_\foamE\;
    \tikzpic{
      \draw[thin,snakecd,] (3,1) to (4,1);
      \funzip[2][1]
      \fzip[0][0]\frst[2][0]\flst[3][0]\fzip[4][0]
      \begin{scope}
        \clip (.5,1) rectangle (4.5,3);
        \begin{scope}[xscale=3,yscale=2]
          \fcap[1/3][1/2]
        \end{scope}
      \end{scope}
      \begin{scope}[yscale=2]
        \frst[0][1/2]\flst[5][1/2]
      \end{scope}
    }[scale=.5*\scl]
    \;\sim_\foamE\;
    \tikzpic{
      \tempdiag
    }[scale=.7*\scl,yscale=-1]
  \end{gather*}
  We conclude using braid-like naturalities (\cref{lem:full_foam_braided_naturalities}), pivotal naturalities (\cref{lem:foam_pivotal_naturalities}) and spatial-like $\foamB$-congruences (\cref{lem:foam_spatial_B-confluence}).
\end{proof}

\begin{lemma}
  \label{lem:foam_T_congruent_iff_S_congruent}
  Every monomial $(\foamlS\setminus\foamlT)$-rewriting step is $\foamE$\nbd-congruent to a $\foamB$\nbd-con\-gruence.
\end{lemma}

\begin{proof}
  By \cref{prop:charactization_nc_sq}, a monomial $(\foamlS\setminus\foamlT)$-rewriting step $f$ of type $\nc$ is $\foamB$\nbd-congruent to a rewriting step which, up to contextualization, has the form on the left:
  \begin{gather*}
    \tikzpic{
      \fcap[0][1]
      \flst\frst[1][0]
    }[scale=.8]
    \;\overset{\nc}{\Rrightarrow}\;
    \tikzpic{
      \fcap[0][2]\fcup[0][1]\fcap
      \fdot[.5][2-.2]
    }[scale=.5]
    \;+\;
    \tikzpic{
      \fcap[0][2]\fcup[0][1]\fcap
      \fdot[.5][.2]
    }[scale=.5]
    \;\overset{\bb^*}{\Rrightarrow}\;
    \tikzpic{
      \fcap[0][1]
      \flst\frst[1][0]
    }[scale=.8]
    \mspace{80mu}
    \tikzpic{
      \clip (-.5,0) rectangle (2.5,6);
      \begin{scope}[yscale=4]
        \flst[0][0]\frst[.5][0][2]\flst[1.5][0][2]\frst[2][0]
      \end{scope}
      \funzip[.5][4][2]
      \begin{scope}[xscale=2,yscale=2]
        \fcap[0/2][4/2][1]
      \end{scope}
    }[scale=.35,yscale=.7]
    \Rrightarrow Z^{-1}\;
    \tikzpic{
      \clip (-1,0) rectangle (2,6);
      \fzip[0][3][2]
      \funzip[0][0][2]
      \begin{scope}[scale=2]
        \fcup[-.25][1]
        \fcap[-.25][0]
      \end{scope}
      \funzip[0][4][2]
      \begin{scope}[xscale=2,yscale=2]
        \fcap[-.5/2][4/2][1]
      \end{scope}
    }[scale=.35,yscale=.7]
  \end{gather*}
  It is readily seen to be $\foamB$-congruent, as pictured.
  Similarly, for type $\sq$ we can reduce the argument to the rewriting step pictured above on the right, evaluating the bubble to find a $\foamB$-congruence (see \cite[Lemma~2.24]{SV_OddKhovanovHomology_2023} for the relation in $\gfoam_d$).
\end{proof}

As a direct corollary, we have that $\foamlT$ and $\foamlS$ present the same underlying module (\cref{lem:EquivTandSfoam}). Moreover:

\begin{lemma}
  \label{lem:foam_independent_branching_confluates}
  Every independent $\foamlT$-branching is $\foamsucc$-tamely $\foamlT$\nbd-cong\-ruent.
\end{lemma}

\begin{proof}
  This follows from \cref{lem:HLRSM_independent-branching-B-confluence}.
\end{proof}

\subsubsection{Confluence of monomial \texorpdfstring{$(\foamlA,\foamlB)$}{(A,B)}-branchings}

Recall the notation $\lA=\lT\setminus\lB$.
We use the terminology \emph{$(\foamlA,\foamlB)$-branching} to refer to a branching which has one branch in $\foamlA$ and the other in $\foamlB$.

\begin{lemma}
  \label{lem:foam_monomial_confluence_A_B}
  Every monomial $(\foamlA,\foamlB)$-branching is $\foamsucc$-tamely $\foamlT$\nbd-con\-gruent.
\end{lemma}

Here and in what follows, we shall use that independent $\foamlT$-branchings are $\foamsucc$-tamely $\foamlT$\nbd-cong\-ruent (\cref{lem:foam_independent_branching_confluates}) without further mention. Similarly, we also stop explicitly mentioning the \ARSMbranchwisetamedcongruencelemma{}.

\begin{proof}
  Let $[f,e,g]$ be a monomial $(\foamlA,\foamlB)$-local triple.
  As a result of the characterization of $\foamB$ (\cref{lem:characterization_foam_B}) and the characterization of $\{\nc,\sq\}$ (\cref{prop:charactization_nc_sq}), we can freely choose the combinatorial representatives of $f$ and $g$. For instance, if $f$ and $g$ have distinct combinatorial data then we can choose the combinatorial representatives of $f$ and $g$ so that $[f,e,g]$ is an independent branching and hence $\foamsucc$-tamely $\foamlT$\nbd-cong\-ruent.
  In particular, if $g$ is of type $\dd$ then $[f,e,g]$ is automatically $\foamsucc$-tamely $\foamlT$\nbd-cong\-ruent. A similar reasoning applies if $g$ is of type $\dm$, choosing to apply the dot migration away from where the neck-cutting or the squeezing relation happens.

  Assume that $f$ is of type $\nc$ and $g$ is of type $\bbev$, such that the strand in the combinatorial data of $g$ (call it $s_1$) is one of the two $i$-strands in the combinatorial data of $f$ (call them $s_1$ and $s_2$).
  The fact that $s_1$ is closed forces $s\coloneqq s_1=s_2$, so that $f$ was in fact not a $\foamlT$-rewriting step to begin with.
  If instead $g$ is of type $\bbodd$ sharing part of its combinatorial data with $f$,
  then $[f,e,g]$ is branchwise $\foamsucc$-tamely $\foamlT$\nbd-cong\-ruent to a contextualization of the $\foamlT^+$-confluent branching pictured in \cref{fig:critical_branching_B_nc}. Since the given $\foamlT^+$-confluence belongs to $\foamB$, the \HLRSMcontextualizationlemma{} applies.
  \begin{figure}[h]
    \begin{gather*}
      \def\scl{.8}
      \begin{tikzcd}[ampersand replacement=\&,row sep=small] 
        \tikzpic{
          \funzip[0][1]\frst[2][1]
          \fzip\frst[2][0]
        }[scale=.4]
        \& 
        Z\;
        \tikzpic{
          \fdot[1][1]\frst[2][1]
          \frst[2][0]
        }[scale=.4]
        \;+\;
        XYZ\;
        \tikzpic{
          \fdot[1][1][2]\frst[2][1]
          \frst[2][0]
        }[scale=.4]
        \&
        (1+XY)Z\;
        \tikzpic{
          \fdot[1][1][2]\frst[2][1]
          \frst[2][0]
        }[scale=.4]
        \\
        \tikzpic{
          \funzip[0][1]\frst[2][1]
          \fzip\frst[2][0]
        }[scale=.4]
        \&
        \tikzpic{
          \funzip[0][3]\frst[2][3]
          \frst[0][2]\fcup[1][2]
          \frst[0][1]\fcap[1][1]
          \fzip\frst[2][0]
          \fdot[1.5][3]
        }[scale=.4]
        \;+\;
        \tikzpic{
          \funzip[0][3]\frst[2][3]
          \frst[0][2]\fcup[1][2]
          \frst[0][1]\fcap[1][1]
          \fzip\frst[2][0]
          \fdot[1.5][1]
        }[scale=.4]
        \&
        (1+XY)Z\;
        \tikzpic{
          \fdot[1][1][2]\frst[2][1]
          \frst[2][0]
        }[scale=.4]
        \arrow[from=1-1,to=1-2,"\bbodd"]
        \arrow[from=1-2,to=1-3,"\dm"]
        \arrow[from=2-1,to=2-2,"\nc"]
        \arrow[from=2-2,to=2-3,"\dm"]
        \arrow[from=1-1,to=2-1,equals]
        \arrow[from=1-3,to=2-3,equals]
      \end{tikzcd}
    \end{gather*}
    
    \caption{Critical branching between types $\foamB$ and type $\nc$.}
    \label{fig:critical_branching_B_nc}
  \end{figure}

  Similar arguments can be given if $f$ is of type $\sq$ and $g$ is of type $\bbev$ (resp.\ $\bbodd$), showing that $[f,e,g]$ is branchwise $\foamsucc$-tamely $\foamlT$\nbd-con\-gruent to a contextualization of the $\foamlT^+$\nbd-con\-fluent branching pictured in \cref{fig:critical_branching_B_sq}.
  In that case however, the $\foamlT^+$\nbd-con\-fluences are not in $\foamB$, and more care is needed to apply the \HLRSMcontextualizationlemma{}; in fact, it \emph{does not} apply to the first of the two critical branchings in \cref{fig:critical_branching_B_sq}.
  One can give yet another critical branching to deal with this case. Instead, we describe another argument that does not require further computation and works for all critical branchings. This argument is essentially the same argument as used in \cref{lem:HLRSM_independent-branching-B-confluence}.

  Note that each confluence in \cref{fig:critical_branching_B_nc} and \cref{fig:critical_branching_B_sq} contains at most one $\foamlA^+$\nbd-rewri\-ting step; call it $r$. Assume that for some context $\Gamma$, one of these confluences fails to verify the hypothesis of the \HLRSMcontextualizationlemma{}.
  That means that $\Gamma[r]$ is not in $\foamlA$, and hence by \cref{lem:foam_T_congruent_iff_S_congruent}, $[f,e,g]$ is $\foamB$\nbd-cong\-ruent. We conclude that $[f,e,g]$ is $\foamB^+$\nbd-conf\-luent, thanks to the convergence of $\foamB$ (\cref{cor:foamB_convergent}).
\end{proof}

\begin{figure}
  \def\scl{.4}
  \begin{gather*}
    \begin{tikzcd}[ampersand replacement=\&,row sep=3em]
      Z\;
      \tikzpic{
        \begin{scope}[yscale=2]
          \flst\frst[1][0]
        \end{scope}
        \fdot[.5][1][2]
      }[scale=\scl]
      +XYZ\;
      \tikzpic{
        \begin{scope}[yscale=2]
          \flst\frst[1][0]
        \end{scope}
        \fdot[.5][1][3]
      }[scale=\scl]
      \&
      Z\;
      \tikzpic{
        \fdot[.5][2-.2][1]
        \fcap\fcup[0][1]
        \flst[0][-1]\fdot[.5][-.5][2]\frst[1][-1]
      }[scale=\scl]
      +Z\;
      \tikzpic{
        \fcap\fcup[0][1]
        \fdot[.5][.2][1]
        \flst[0][-1]\fdot[.5][-.5][2]\frst[1][-1]
      }[scale=\scl]
      +XYZ\;
      \tikzpic{
        \fdot[.5][2-.2][1]
        \fcap\fcup[0][1]
        \flst[0][-1]\fdot[.5][-.5][3]\frst[1][-1]
      }[scale=\scl]
      +XYZ\;
      \tikzpic{
        \fcap\fcup[0][1]
        \fdot[.5][.2][1]
        \flst[0][-1]\fdot[.5][-.5][3]\frst[1][-1]
      }[scale=\scl]
      \\
      \tikzpic{
        \begin{scope}[yscale=2]
          \flst\frst[2][0]
        \end{scope}
        \fzip[.5][0][2]\funzip[.5][1][2]
      }[scale=\scl]
      \&
      Z\;
      \tikzpic{
        \fdot[.5][2-.2][2]
        \fcap\fcup[0][1]
        \fdot[.5][.2][2]
      }[scale=\scl]
      +XYZ\;
      \tikzpic{
        \fdot[.5][2-.2][2]
        \fcap\fcup[0][1]
        \fdot[.5][.2][3]
      }[scale=\scl]
      +XYZ\;
      \tikzpic{
        \fdot[.5][2-.2][3]
        \fcap\fcup[0][1]
        \fdot[.5][.2][2]
      }[scale=\scl]
      \\
      Z^{-1}\;
      \tikzpic{
        \clip (-1,-1) rectangle (2,5);
        \flst[-.5][4]\funzip[0][4][2]\frst[1.5][4]
        \fzip[0][3][2]
        \funzip[0][0][2]
        \begin{scope}[scale=2]
          \fcup[-.25][1]
          \fcap[-.25][0]
        \end{scope}
        \flst[-.5][-1]\fzip[0][-1][2]\frst[1.5][-1]
      }[scale=\scl*.7][(0,2*.35)]
      \&
      Z\;
      \tikzpic{
        \fdot[.5][2-.2][2]
        \fcap\fcup[0][1]
        \fdot[.5][.2][2]
      }[scale=\scl]
      +XYZ\;
      \tikzpic{
        \fdot[.5][2-.2][3]
        \fcap\fcup[0][1]
        \fdot[.5][.2][2]
      }[scale=\scl]
      +XYZ\;
      \tikzpic{
        \fdot[.5][2-.2][2]
        \fcap\fcup[0][1]
        \fdot[.5][.2][3]
      }[scale=\scl]
      +Z\;
      \tikzpic{
        \fdot[.5][2-.2][3]
        \fcap\fcup[0][1]
        \fdot[.5][.2][3]
      }[scale=\scl]
      \arrow[from=2-1,to=1-1,"\bbodd"]
      \arrow[from=1-1,to=1-2,"\nc"]
      \arrow[from=1-2,to=2-2,"{\{\dm,\dd\}}"{pos=.7},decorate,
      decoration={snake,amplitude=1.5pt,segment length=8pt,post length=.6cm}]
      \arrow[from=2-1,to=3-1,"\sq"']
      \arrow[from=3-1,to=3-2,"\bbodd^*"]
      \arrow[from=3-2,to=2-2,"\dd"'{pos=.7},decorate,
      decoration={snake,amplitude=1.5pt,segment length=8pt,post length=1cm}]
    \end{tikzcd}
    \\[3ex]
    \begin{tikzcd}[ampersand replacement=\&,row sep=3em,column sep=2em]
      \&[-6em]\&[-6em]
      Z\;
      \tikzpic{
        \begin{scope}[yscale=2]
          \frst[0][0][2]\flst[1][0][2]\frst[2][0]
        \end{scope}
        \fdot[-.5][1]
      }[scale=\scl]
      \;+\;
      XYZ\;
      \tikzpic{
        \begin{scope}[yscale=2]
          \frst[0][0][2]\flst[1][0][2]\frst[2][0]
        \end{scope}
        \fdot[-.5][1][2]
      }[scale=\scl]
      \&[-6em]\&[-7em]
      \\
      \tikzpic{
        \begin{scope}[yscale=2]
          \frst[0][0][2]\flst[1][0][2]\frst[2][0]
        \end{scope}
        \fzip[-1.5][0]\funzip[-1.5][1]
      }[scale=\scl]
      \&\&\&\&
      Z\;
      \tikzpic{
        \begin{scope}[yscale=2]
          \frst[0][0][2]\flst[1][0][2]\frst[2][0]
        \end{scope}
        \fdot[1.5][1][2]
      }[scale=\scl]
      \;+\;
      XYZ\;
      \tikzpic{
        \begin{scope}[yscale=2]
          \frst[0][0][2]\flst[1][0][2]\frst[2][0]
        \end{scope}
        \fdot[-.5][1][2]
      }[scale=\scl]
      \\
      \&
      Z^{-1}\;
      \tikzpic{
        \clip (-2,-1) rectangle (2,5);
        \funzip[-1.5][4][1]\frst[0][4][2]\flst[1][4][2]\frst[1.5][4]
        \fzip[0][3][2]
        \funzip[0][0][2]
        \begin{scope}[scale=2]
          \fcup[-.25][1]
          \fcap[-.25][0]
        \end{scope}
        \fzip[-1.5][-1][1]\frst[0][-1][2]\flst[1][-1][2]\frst[1.5][-1]
        \draw[foamdraw1] (-1.5,0) to (-1.5,4);
      }[scale=\scl*.7]
      \&\&
      XZ^{-1}\;
      \tikzpic{
        \begin{scope}[yscale=2]
          \frst[2][0]
        \end{scope}
        \funzip[0][0][2]\fzip[0][1][2]
      }[scale=\scl]
      \&
      %
      \arrow[from=2-1,to=1-3,"\bbodd"]
      \arrow[from=1-3,to=2-5,"\dm"]
      \arrow[from=2-1,to=3-2,"\sq"']
      \arrow[from=3-2,to=3-4,snakecd]
      \arrow[from=3-4,to=2-5,"\nc"']
    \end{tikzcd}
  \end{gather*}

  \caption{Critical branchings between types $\foamB$ and type $\sq$. Colours are $i=\textdot[1]$, $i+1=\textdot[2]$ and $i+2=\textdot[3]$.}
  \label{fig:critical_branching_B_sq}
\end{figure}

\subsubsection{Confluence of monomial \texorpdfstring{$\foamlA$}{A}-branchings}

\begin{lemma}
  \label{lem:foam_monomial_confluence_A_A}
  Every monomial $\foamlA$-branching is $\foamsucc$-tamely $\foamlT$\nbd-cong\-ruent.
\end{lemma}

\begin{proof}
  Let $[f,e,g]$ be a monomial $\foamlA$-local triple.
  As in the proof of \cref{lem:foam_monomial_confluence_A_B}, the characterization of $\{\nc,\sq\}$ \cref{prop:charactization_nc_sq} implies that we can freely choose the combinatorial representatives of $f$ and $g$.
  Contrary to the proof of \cref{lem:foam_monomial_confluence_A_B} however, even if $f$ and $g$ have distinct combinatorial data, there may not exist combinatorial representatives for which $f$ and $g$ are independent.
  On the other hand, we shall see that in most cases, we can choose the combinatorial representatives of $f$ and $g$ such that $[f,e,g]$ \emph{rewrites} into an independent branching, in the sense of \cref{subsubsec:HLRSM_independent_rewriting}.
  Only two cases will not follow this scheme: they are the critical branchings given in \cref{fig:critical_branching_nc_sq}.

  \medbreak

  Consider first the case where both $f$ and $g$ are of type $\nc$. If their respective colours are $i$ and $j$, we can assume that either $j=i$ or $j=i+1$ (otherwise, we can choose combinatorial representatives such that $[f,e,g]$ is an independent branching).
  In these cases, we can choose combinatorial representatives such that $[f,e,g]=\Gamma[f',e',g']$ for some context $\Gamma$ and $[f',e',g']$ is encoded in the diagram on the left-hand side below (the two wiggly lines encode $f$ and $g$):
  \begin{gather*}
    \def\scl{.5}
    \tikzpic{
      \draw[snakecd,thin] (0,1) to (3,1);
      \draw[snakecd,thin] (1.5,1-.3) to (1.5,1+.3);
      \begin{scope}[yscale=2]
        \flst[0][0][2]\frst[3][0][2]
      \end{scope}
      \funzip[1][0]\fzip[1][1]
      \node[below=-1pt] at (0,0) {\scriptsize $i$};
      \node[below=-1pt] at (1,0) {\scriptsize $j$};
    }[scale=\scl,yscale=.8*6/2][(0,.8*6/2*\scl)]
    \mspace{130mu}
    \tikzpic{
      \draw[snakecd,thin] (0,1) to (3,1);
      \draw[snakecd,thin] (2,.3) to (3,.3);
      \draw[snakecd,thin] (0,2-.3) to (1,2-.3);
      \draw[snakecd,thin] (1.5,1-.3) to (1.5,1+.3);
      \begin{scope}[yscale=2]
        \flst\frst[3][0]
      \end{scope}
      \funzip[1][0]\fzip[1][1]
      \node[below=-1pt] at (0,0) {\scriptsize $i$};
    }[scale=\scl,yscale=.8*6/2][(0,.8*6/2*\scl)]
    \;\overset{\nc^*}{\longrightarrow}\; 
    \sum_{\substack{\text{suitable}\\\text{dots}}}\;
    \tikzpic{
      \draw[snakecd,thin] (0,3) to (3,3);
      \draw[snakecd,thin] (1.5,3-.4) to (1.5,3+.4);
      \fcup[0][5]\flst[2][5]\frst[3][5]
      \fcap[0][4]\flst[2][4]\frst[3][4]
      \flst[0][3]\fzip[1][3]\frst[3][3]
      \flst[0][2]\funzip[1][2]\frst[3][2]
      \flst[0][1]\frst[1][1]\fcup[2][1]
      \flst[0][0]\frst[1][0]\fcap[2][0]
      \node[below=-1pt] at (0,0) {\scriptsize $i$};
    }[scale=\scl,yscale=.8][(0,3*.8*\scl)]
  \end{gather*}
  If $j=i$, then $\Gamma[f',e',g']$ rewrites via two neck-cuttings into a linear combination of branchings $\Gamma[f'',e'',g'']$ with $f''$ and $g''$ associated to the same combinatorial data, as pictured in the schematic on the right-hand side above (we only picture the $i$-strands, leaving the dots implicit).
  If $\Gamma$ does not connect any of the $i$-strands involved, then all monomial rewriting steps involved are in $\foamlT$, and we can apply the \HLRSMindependentrewritinglemma{}.
  Otherwise, it means that only one, or none, neck-cutting is necessary to get to the same situation, and the \HLRSMindependentrewritinglemma{} is still applicable. In any case, $[f,e,g]$ is $\foamsucc$-tamely $\foamlT$\nbd-cong\-ruent.

  If $j=i+1$, the $\foamlT^+$-confluence of $[f',e',g']$ is the first critical branching pictured in \cref{fig:critical_branching_nc_sq}. As this confluence only uses one rewriting step of type $\foamlA$, we conclude as in \cref{lem:foam_monomial_confluence_A_B} that $[f,e,g]=\Gamma[f',e',g']$ is $\foamsucc$-tamely $\foamlT$\nbd-cong\-ruent.

  \medbreak

  Consider now the case where $f$ is of type $\sq$ and $g$ is of type $\nc$. Denote $(i,i+1)$ and $j$ their respective colour. As before, we can assume that $j=i-1$, $i$, $i+1$ or $i+2$. Then, up to choice of combinatorial representatives and vertical symmetry, $[f,e,g]$ is equal to $\Gamma[f',e',g']$ for some context $\Gamma$ and $[f',e',g']$ as pictured below:
  \begin{gather*}
    \tikzpic{
      \draw[snakecd,thin] (0,1) to (5,1);
      \draw[snakecd,thin] (1.5,1-.4) to (1.5,1+.4);
      \begin{scope}[yscale=2]
        \flst\frst[3][0][2]\flst[4][0][2]\frst[5][0]
      \end{scope}
      \fzip[1][1][3]\funzip[1][0][3]
      \node[below=-1pt] at (0,0) {\scriptsize $i$};
      \node[below=-1pt] at (2,0) {\scriptsize $j$};
      \node[below=-1pt] at (4,0) {\scriptsize $i+1$};
    }[scale=.5]
  \end{gather*}
  If $j=i+2$, we can further isotope the branching to get an independent branching. If $j=i$ or if $j=i+1$, we can rewrite it using neck-cuttings and conclude similarly as above.

  The last case $j=i-1$ leads to the second critical branching pictured in \cref{fig:critical_branching_nc_sq}. Note that flipping everything vertically leads to the same critical branching because $j$ and $i+1$ are distant colours.
  The confluence uses two rewriting steps of type $\foamlA$, one of type $\sq$ and the other of type $\nc$; we refer to them simply as ``$\sq$'' and ``$\nc$''. If $\Gamma[\text{``$\sq$''}]$ does not belong to $\foamlT$, then $f$ does not belong to $\foamlT$ either; hence if $[f,e,g]=\Gamma[f',e',g']$ is a $\foamlT^+$-branching, $\Gamma[\text{``$\sq$''}]$ necessarily belongs to $\foamlT$.
  On the other hand, if $\Gamma[\text{``$\nc$''}]$ does not belong to $\foamlT$, we can replace it with a $\foamB^+$-confluence.
  We conclude that $[f,e,g]$ is $\foamsucc$-tamely $\foamlT$\nbd-cong\-ruent.
  
  \medbreak

  \begin{figure}
    \def\tempdiagA{
      \begin{scope}[yscale=4]
        \flst[0][0]\frst[.5][0][2]\flst[1.5][0][2]\frst[2][0]
      \end{scope}
    }
    \def\tempdiagB{
      \funzip[.5][0][2]\fzip[.5][3][2]
      \begin{scope}[yscale=4]
        \flst[0][0]\frst[2][0]
      \end{scope}
    }
    \def\tempdiagC{
      \clip (-1,0) rectangle (2,4);
      \fzip[0][3][2]
      \funzip[0][0][2]
      \begin{scope}[scale=2]
        \fcup[-.25][1]
        \fcap[-.25][0]
      \end{scope}
    }
    \def\tempdiagD{
      \begin{scope}[yscale=4]
        \flst[0][0][2]\frst[1][0][1]\flst[2][0][1]\frst[3][0][3]\flst[4][0][3]\frst[5][0][2]
      \end{scope}
    }
    \begin{gather*}
      \def\scl{.35}
      \begin{tikzcd}[ampersand replacement=\&,row sep=3em,column sep=2em]
        \tikzpic{
          \tempdiagB
        }[scale=\scl]
        \&
        \tikzpic{
          \tempdiagC
          \fdot[.5][3][1][1.5]
        }[scale=\scl]
        \;+\;
        \tikzpic{
          \tempdiagC
          \fdot[.5][1][1][1.5]
        }[scale=\scl]
        \&
        (X+Y)Z\;
        \tikzpic{
          \tempdiagC
          \fdot[1.2][3.5][3][1.5]
        }[scale=\scl]
        \\
        \tikzpic{
          \tempdiagB
        }[scale=\scl]
        \&
        XZ^2
        \tikzpic{
          \tempdiagA
          \fdot[1.75][2][2][1.5]
        }[scale=\scl]
        \;+\;
        YZ^2
        \tikzpic{
          \tempdiagA
          \fdot[.25][2][2][1.5]
        }[scale=\scl]
        \&
        (X+Y)Z\;
        \tikzpic{
          \tempdiagC
          \fdot[1.2][3.5][3][1.5]
        }[scale=\scl]
        \arrow[from=1-1,to=1-2,"\nc"]
        \arrow[from=1-2,to=1-3,"\dm"{pos=.5},
        ]
        \arrow[from=2-1,to=2-2,"\nc"]
        \arrow[from=2-2,to=2-3,"\sq\cdot\dm"]
        \arrow[from=1-1,to=2-1,equals]
        \arrow[from=1-3,to=2-3,equals]
      \end{tikzcd}
      \\[8ex]
      \def\scl{.35}\def\xscl{.8}
      \begin{tikzcd}[ampersand replacement=\&,column sep=0em,row sep=3em]
        \&[-3em]\&[-8em]
        XZ^2\;
        \tikzpic{
          \tempdiagD
          \fdot[2.5][2][1][1.5]
        }[scale=\scl,xscale=\xscl]
        +
        YZ^2\;
        \tikzpic{
          \tempdiagD
          \fdot[.5][2][1][1.5]
        }[scale=\scl,xscale=\xscl]
        \&[-6em]\&[-9em]
        \\
        \tikzpic{
          \funzip[1][0][1]\fzip[1][3][1]
          \begin{scope}[yscale=4]
            \flst[0][0][2]\frst[3][0][3]\flst[4][0][3]\frst[5][0][2]
          \end{scope}
        }[scale=\scl,xscale=\xscl]
        \&\&\&\&
        XZ^2\;
        \tikzpic{
          \tempdiagD
          \fdot[4.5][2][3][1.5]
        }[scale=\scl,xscale=\xscl]
        +
        YZ^2\;
        \tikzpic{
          \tempdiagD
          \fdot[2.5][2][3][1.5]
        }[scale=\scl,xscale=\xscl]
        \\
        \&
        Z^{-1}\;
        {}\tikzpic{
          \clip (-.5,0) rectangle (5.5,7);
          \fzip[1][5][1]\frst[3][5][3]\flst[4][5][3]
          \fzip[3][4][3]
          \begin{scope}[xscale=5,yscale=3]
            \fcup[0][3/3][2]
            \fcap[0/5][0/3][2]
          \end{scope}
          \funzip[3][1][3]
          \funzip[1][0][1]\frst[3][0][3]\flst[4][0][3]
        }[scale=\scl,xscale=\xscl,yscale=4/6][(0,3*\scl*4/6)]
        \&\&
        \tikzpic{
          \funzip[3][0][3]\fzip[3][3][3]
          \begin{scope}[yscale=4]
            \flst[0][0][2]\frst[1][0][1]\flst[2][0][1]\frst[5][0][2]
          \end{scope}
        }[scale=\scl,xscale=\xscl]
        \&
        \arrow[
          from=2-1,to=1-3,"\nc",
          end anchor={[xshift=3em]south west}
        ]
        \arrow[
          from=1-3,to=2-5,"\dm"{pos=.2},"\dm"{pos=.75},
          start anchor={[xshift=-4em]south},
          end anchor={[xshift=-4em]north},
          decorate,
          decoration={snake,amplitude=1.5pt,segment length=8pt,post length=1cm,pre length=1cm}
        ]
        \arrow[
          from=1-3,to=2-5,"\dm^*"{pos=.3},
          start anchor={[xshift=6em]south},
          end anchor={[xshift=6em]north},
          decorate,
          decoration={snake,amplitude=1.5pt,segment length=8pt,post length=.1cm,pre length=2cm}
        ]
        \arrow[from=2-1,to=3-2,"\sq"]
        \arrow[
          from=3-2,to=3-4,"\sq"{pos=.75},
          decorate,
          decoration={snake,amplitude=1.5pt,segment length=8pt,post length=.5cm}
        ]
        \arrow[from=3-4,to=2-5,"\nc"]
      \end{tikzcd}
    \end{gather*}
  
    \caption{Critical branchings in types $\{\nc,\sq\}$. Colours are $i=\textdot[1]$, $i+1=\textdot[2]$ and $i+2=\textdot[3]$.}
    \label{fig:critical_branching_nc_sq}
  \end{figure}

  Finally, consider the case where both $f$ and $g$ are of type $\sq$, with respective colours $(i,i+1)$ and $(j,j+1)$. Without loss of generality we can assume that $j\geq i$, and furthermore that either $j=i$, $j=i+1$ or $j=i+2$.
  Each case will allow different choices of combinatorial representatives. To help the exposition, we fix the positions of the $i$ and $(i+1)$-strands associated to $f$, simply called the \emph{$i$\nbd- and $(i+1)$-strands} below, and discuss how the $j$- and $(j+1)$-strands associated to $g$, simply called the \emph{$j$- and $(j+1)$-strands} below, can be isotoped with regard to the $i$- and $(i+1)$-strands.

  If $j=i+2$, then $j$ is adjacent to $i+1$. In particular, the two $j$-strands cannot be isotoped through the $(i+1)$-strands. As there exists an isotopy joining the $j$ and $(j+1)$-strands, the four strands must lie on one side of the $(i+1)$-strands (at least partially for the $(j+1)$-strands), as pictured in the following schematic (replacing wiggly lines with straight lines to avoid clutter):
  \tikzset{rdirected/.style={}}
  \tikzset{directed/.style={}}
  \[\begin{tikzpicture}[scale=.6]
    \node[above left=-4pt] at (140:2) {\scriptsize $i$};
    \node[below left=-4pt] at (-110:2) {\scriptsize $i+1$};
    \node[above=-2pt] at (80:2) {\scriptsize $j$};
    \node[above left=-4pt] at (180:2) {\scriptsize $j+1$};
    \draw (0,0) circle (2cm);
    \clip (0,0) circle (2cm);
    \fill[colour_diag4,opacity=.2] (-2,-.2) rectangle (2,.2);
    \draw[diag4=1,directed=.6] (-2,.2) to (2,.2);
    \draw[diag4=1,rdirected=.6] (-2,-.2) to (2,-.2);
    \fill[colour_diag3,opacity=.2] (2,.3) to (.3,2) to (.3+1,2+1) to (2+1,.3+1) to (2,.3);
    \fill[colour_diag3,opacity=.2] (2,-.3) to (.3,-2) to (.3+1,-2-1) to (2+1,-.3-1) to (2,-.3);
    \draw[diag3=1,directed] (2,.3) to (.3,2);
    \draw[diag3=1,rdirected] (2,-.3) to (.3,-2);
    \fill[colour_diag1,opacity=.2] (1,2) rectangle (2,-2);
    \fill[colour_diag1,opacity=.2] (-1.5,2) rectangle (-2,-2);
    \draw[diag1=1,rdirected=.7] (1,2) to (1,-2);
    \draw[diag1=1,directed=.7] (-1.5,2) to (-1.5,-2);
    \fill[colour_diag2,opacity=.2] (-.75,2) rectangle (-.25,-2);
    \draw[diag2=1,rdirected=.7] (-.75,2) to (-.75,-2);
    \draw[diag2=1,directed=.7] (-.25,2) to (-.25,-2);
    \draw (-1.5,.8) to (1,.8);
    \draw (2-.4,.3+.4) to (2-.4,-.3-.4);
  \end{tikzpicture}\]
  On the side where the $j$-strands are, there is one $i$-strand. Given that $i$ is distant from $j$ and $j+1$, it can be isotoped through the $j$ and $(j+1)$-strands. In that way, $[f,e,g]$ is branchwise $\foamE$\nbd-congruent to an independent branching.

  If $j=i+1$, then the $(i+1)$-strands and $(j+1)$-strands cannot intersect. This leads to three possible schematics, up to symmetries:
  \begin{gather*}
    \def\scl{.8}
    \tikzpic{
      \draw (0,0) circle (2cm);
      \clip (0,0) circle (2cm);
      \draw[diag1=1,directed=.5] (-1.5-.2,2) to (-1.5-.2,-2);
      \fill[fill=colour_diag1,opacity=.2] (-2,-2) rectangle (-1.7,2);
      \draw[diag2=1,rdirected=.5] (1.5+.2,2) to (1.5+.2,-2);
      \fill[fill=colour_diag2,opacity=.2] (2,-2) rectangle (1.7,2);
      \draw[diag2=1,rdirected=.7] (-1-.2,2) to (-1-.2,-2);
      \draw[diag2=1,directed=.4] (-.7-.2,2) to (-.7-.2,-2);
      \fill[fill=colour_diag2,opacity=.2] (-.9,-2) rectangle (-1.2,2);
      \draw[diag3=1,rdirected=.7] (1,2) to (1,-2);
      \draw[diag3=1,directed=.5] (.7,2) to (.7,-2);
      \fill[fill=colour_diag3,opacity=.2] (1,2) rectangle (.7,-2);
      \draw[diag1=1,rdirected=.3] 
        (.4,2) .. controls ++(0,-2) and ++(-.2,-1.5) .. (1.5,2);
      \fill[fill=colour_diag1,opacity=.1] 
        (.4,2) .. controls ++(0,-2) and ++(-.2,-1.5) .. (1.5,2);
      \draw[] (-1.7,-.2) to[out=0,in=-150] (.6,.8);
      \draw[] (-.9,-.3) to (1.7,-.3);
    }[scale=.6*\scl]
    \mspace{50mu}
    \tikzpic{
      \draw (0,0) circle (2cm);
      \clip (0,0) circle (2cm);
      \draw[diag1=1,directed=.5] (-1.5-.2,2) to (-1.5-.2,-2);
      \fill[fill=colour_diag1,opacity=.2] (-2,-2) rectangle (-1.7,2);
      \draw[diag2=1,rdirected=.5] (1.5+.2,2) to (1.5+.2,-2);
      \fill[fill=colour_diag2,opacity=.2] (2,-2) rectangle (1.7,2);
      \draw[diag2=1,rdirected=.7] (-1-.2,2) to (-1-.2,-2);
      \draw[diag2=1,directed=.4] (-.7-.2,2) to (-.7-.2,-2);
      \fill[fill=colour_diag2,opacity=.2] (-.9,-2) rectangle (-1.2,2);
      \draw[diag3=1,rdirected=.7] (1,2) to (1,-2);
      \draw[diag3=1,directed=.5] (.7,2) to (.7,-2);
      \fill[fill=colour_diag3,opacity=.2] (1,2) rectangle (.7,-2);
      \draw[diag2=1,directed=.4] 
        (-.5,-2) .. controls ++(.2,1.5) and ++(-.2,1.5) .. (.2,-2);
      \fill[fill=colour_diag2,opacity=.2] 
        (-.5,-2) .. controls ++(.2,1.5) and ++(-.2,1.5) .. (.2,-2);
      \draw[diag1=1,rdirected=.3] 
        (.4,2) .. controls ++(0,-2) and ++(-.2,-1.5) .. (1.5,2);
      \fill[fill=colour_diag1,opacity=.1] 
        (.4,2) .. controls ++(0,-2) and ++(-.2,-1.5) .. (1.5,2);
        \draw[] (-1.7,-.2) to[out=0,in=-150] (.6,.8);
        \draw[] (0,-1.1) to[out=60,in=180] (1.7,-.3);
    }[scale=.6*\scl]
    \mspace{50mu}
    \tikzpic{
      \draw (0,0) circle (2cm);
      \clip (0,0) circle (2cm);
      \draw[diag1=1,directed=.5] (-1.5-.2,2) to (-1.5-.2,-2);
      \fill[fill=colour_diag1,opacity=.2] (-2,-2) rectangle (-1.7,2);
      \draw[diag2=1,rdirected=.5] (1.5+.2,2) to (1.5+.2,-2);
      \fill[fill=colour_diag2,opacity=.2] (2,-2) rectangle (1.7,2);
      \draw[diag2=1,rdirected=.7] (-1-.2,2) to (-1-.2,-2);
      \draw[diag2=1,directed=.4] (-.7-.2,2) to (-.7-.2,-2);
      \fill[fill=colour_diag2,opacity=.2] (-.9,-2) rectangle (-1.2,2);
      \draw[diag3=1,rdirected=.7] (1,2) to (1,-2);
      \draw[diag3=1,directed=.5] (.7,2) to (.7,-2);
      \fill[fill=colour_diag3,opacity=.2] (1,2) rectangle (.7,-2);
      \draw[diag2=1,rdirected=.4] 
        (-.5,2) .. controls ++(.2,-1.5) and ++(-.2,-1.5) .. (.2,2);
      \fill[fill=colour_diag2,opacity=.2] 
        (-.5,2) .. controls ++(.2,-1.5) and ++(-.2,-1.5) .. (.2,2);
      \draw[diag1=1,rdirected=.3] 
        (.4,2) .. controls ++(0,-2) and ++(-.2,-1.5) .. (1.5,2);
      \fill[fill=colour_diag1,opacity=.2] 
        (.4,2) .. controls ++(0,-2) and ++(-.2,-1.5) .. (1.5,2);
      \draw[] (-1.7,.2) to[out=0,in=-150] (.6,.8);
      \draw[] (0,1.1) to[out=-60,in=180] (1.7,.3);
    }[scale=.6*\scl]
    \overset{\nc}{\longrightarrow}
    \tikzpic{
      \draw (0,0) circle (2cm);
      \clip (0,0) circle (2cm);
      \draw[diag1=1,directed=.5] (-1.5-.2,2) to (-1.5-.2,-2);
      \fill[fill=colour_diag1,opacity=.2] (-2,-2) rectangle (-1.7,2);
      \draw[diag2=1,rdirected=.5] (1.5+.2,2) to (1.5+.2,-2);
      \fill[fill=colour_diag2,opacity=.2] (2,-2) rectangle (1.7,2);
      \draw[diag3=1,rdirected=.7] (1,2) to (1,-2);
      \draw[diag3=1,directed=.5] (.7,2) to (.7,-2);
      \fill[fill=colour_diag3,opacity=.2] (1,2) rectangle (.7,-2);
      \draw[diag2=1,rdirected=.5] (-1.2,2) to (-1.2,-2);
      \draw[diag2=1,directed=.5] (.2,2) to (-.9,-2);
      \draw[diag2=1,rdirected=.4] 
        (-.5,2) .. controls ++(-.2,-1.5) and ++(.2,-1.5) .. (-1,2);
      \fill[fill=colour_diag2,opacity=.2] 
      (-.5,2) .. controls ++(-.2,-1.5) and ++(.2,-1.5) .. (-1,2)
        to (-1.2,2) to (-1.2,-2) to (-.9,-2) to (.2,2) to (-1,2);
      \draw[diag1=1,rdirected=.3] 
        (.4,2) .. controls ++(0,-2) and ++(-.2,-1.5) .. (1.5,2);
      \fill[fill=colour_diag1,opacity=.2] 
        (.4,2) .. controls ++(0,-2) and ++(-.2,-1.5) .. (1.5,2);
        \draw[] (-1.7,.2) to[out=0,in=-150] (.6,.8);
        \draw[] (-.65,-1.1) to[out=0,in=-130] (1.7,.3);
    }[scale=.6*\scl]
  \end{gather*}
  In the first schematic, one of the $(i+1)$-strands coincides with one of the $j$-strands. The first two schematics are independent branchings on the nose. The last one rewrites into an independent branching, as pictured.

  Finally, if $j=i$, either the $(i+1)$-strands and $(j+1)$-strands do not intersect, or they coincide.
  Below we only picture the two schematics for which rewriting the branching is necessary:
  \begin{gather*}
    \def\scl{.8}
    \tikzpic{
      \draw (0,0) circle (2cm);
      \clip (0,0) circle (2cm);
      \draw[diag2=1,rdirected] (-2,-2) .. controls ++(1,1) and ++(0,-1) .. (-.2,2);
      \draw[diag2=1,rdirected=.3] (-.2-.4,2) .. controls ++(0,-1) and ++(1,1) .. (-2-.4,-2);
      \fill[fill=colour_diag2,opacity=.2] 
        (-.2-.4,2) .. controls ++(0,-1) and ++(1,1) .. (-2-.4,-2)
        to (-2,-2) .. controls ++(1,1) and ++(0,-1) .. (-.2,2)
        to (-.2-.4,2);
      \draw[diag1=1,rdirected] (-2,-.5) .. controls ++(.5,.5) and ++(.5,-.5) .. (-1.5,2);
      \fill[fill=colour_diag1,opacity=.2] (-2,-.5) .. controls ++(.5,.5) and ++(.5,-.5) .. (-1.5,2) to (-2,2) to (-2,-.5);
      \draw[diag2=1,directed=.7] (2,-2) .. controls ++(-1,1) and ++(0,-1) .. (.2,2);
      \draw[diag2=1,directed] (.2+.4,2) .. controls ++(0,-1) and ++(-1,1) .. (2+.4,-2);
      \fill[fill=colour_diag2,opacity=.2] 
        (.2+.4,2) .. controls ++(0,-1) and ++(-1,1) .. (2+.4,-2)
        to (2,-2) .. controls ++(-1,1) and ++(0,-1) .. (.2,2)
        to (.2+.4,2);
      \draw[diag1=1,directed] (2,-.5) .. controls ++(-.5,.5) and ++(-.5,-.5) .. (1.5,2);
      \fill[fill=colour_diag1,opacity=.2] (2,-.5) .. controls ++(-.5,.5) and ++(-.5,-.5) .. (1.5,2) to (2,2) to (2,-.5);
      \draw[diag1=1,rdirected=.7] (-.2,-2) .. controls ++(0,1.5) and ++(0,1) .. (-1.2,-2);
      \fill[fill=colour_diag1,opacity=.2]
        (-.2,-2) .. controls ++(0,1.5) and ++(0,1) .. (-1.2,-2);
      \draw[diag1=1,directed=.3] (.2,-2) .. controls ++(0,1.5) and ++(0,1) .. (1.2,-2);
      \fill[fill=colour_diag1,opacity=.2]
        (.2,-2) .. controls ++(0,1.5) and ++(0,1) .. (1.2,-2);
      \draw (.4,-1.2) to[out=110,in=0] (-1.5,.4);
      \draw (-.4,-1.2) to[out=180-110,in=180] (1.5,.4);
    }[scale=.6*\scl]
    \overset{\nc}{\longrightarrow}
    \tikzpic{
      \draw (0,0) circle (2cm);
      \clip (0,0) circle (2cm);
      \draw[diag2=1,rdirected] (-2,-2) .. controls ++(1,1) and ++(0,-1) .. (-.2,2);
      \draw[diag2=1,rdirected=.3] (-.2-.4,2) .. controls ++(0,-1) and ++(1,1) .. (-2-.4,-2);
      \fill[fill=colour_diag2,opacity=.2] 
        (-.2-.4,2) .. controls ++(0,-1) and ++(1,1) .. (-2-.4,-2)
        to (-2,-2) .. controls ++(1,1) and ++(0,-1) .. (-.2,2)
        to (-.2-.4,2);
      \draw[diag1=1,rdirected] (-2,-.5) .. controls ++(.5,.5) and ++(.5,-.5) .. (-1.5,2);
      \fill[fill=colour_diag1,opacity=.2] (-2,-.5) .. controls ++(.5,.5) and ++(.5,-.5) .. (-1.5,2) to (-2,2) to (-2,-.5);
      \draw[diag2=1,directed=.7] (2,-2) .. controls ++(-1,1) and ++(0,-1) .. (.2,2);
      \draw[diag2=1,directed] (.2+.4,2) .. controls ++(0,-1) and ++(-1,1) .. (2+.4,-2);
      \fill[fill=colour_diag2,opacity=.2] 
        (.2+.4,2) .. controls ++(0,-1) and ++(-1,1) .. (2+.4,-2)
        to (2,-2) .. controls ++(-1,1) and ++(0,-1) .. (.2,2)
        to (.2+.4,2);
      \draw[diag1=1,directed] (2,-.5) .. controls ++(-.5,.5) and ++(-.5,-.5) .. (1.5,2);
      \fill[fill=colour_diag1,opacity=.2] (2,-.5) .. controls ++(-.5,.5) and ++(-.5,-.5) .. (1.5,2) to (2,2) to (2,-.5);
      \draw[diag1=1,rdirected=.8] (-.2,-2) .. controls ++(0,.6) and ++(0,.6) .. (.2,-2);
      \draw[diag1=1,rdirected=.5] (1.2,-2) .. controls ++(0,1.5) and ++(0,1.5) .. (-1.2,-2);
      \fill[fill=colour_diag1,opacity=.2]
        (-.2,-2) .. controls ++(0,.6) and ++(0,.6) .. (.2,-2)
        to (1.2,-2) .. controls ++(0,1.5) and ++(0,1.5) .. (-1.2,-2)
        to (-.2,-2);
      \draw (-.9,-1.2) to[out=140,in=-60] (-1.65,0);
      \draw (.9,-1.2) to[out=180-140,in=180+60] (1.65,0);
    }[scale=.6*\scl]
    \mspace{50mu}
    \tikzpic{
      \draw (0,0) circle (2cm);
      \clip (0,0) circle (2cm);
      \draw[diag2=1,rdirected=.3] (.3,-2) to (.3,2);
      \draw[diag2=1,rdirected=.3] (-.3,2) to (-.3,-2);
      \fill[fill=colour_diag2,opacity=.2]
        (.3,-2) to (.3,2) to (-.3,2) to (-.3,-2) to (.3,-2);
      \draw[diag1=1,directed=.3] (.9,-2) .. controls ++(0,2) and ++(0,2).. (1.8,-2);
      \fill[fill=colour_diag1,opacity=.2]
        (.9,-2) .. controls ++(0,2) and ++(0,2).. (1.8,-2);
      \draw[diag1=1,rdirected=.3] (.9,2) .. controls ++(0,-2) and ++(0,-2).. (1.8,2);
      \fill[fill=colour_diag1,opacity=.2]
        (.9,2) .. controls ++(0,-2) and ++(0,-2).. (1.8,2);
      \draw[diag1=1,directed=.3] (-.9,2) .. controls ++(0,-2) and ++(0,-2).. (-1.8,2);
      \fill[fill=colour_diag1,opacity=.2]
        (-.9,2) .. controls ++(0,-2) and ++(0,-2).. (-1.8,2);
      \draw[diag1=1,rdirected=.3] (-.9,-2) .. controls ++(0,2) and ++(0,2).. (-1.8,-2);
      \fill[fill=colour_diag1,opacity=.2]
        (-.9,-2) .. controls ++(0,2) and ++(0,2).. (-1.8,-2);
      \draw (-1,-1) to (1,1);
      \draw (1,-1) to (-1,1);
    }[scale=.6*\scl]
    \overset{\nc}{\longrightarrow}
    \tikzpic{
      \draw (0,0) circle (2cm);
      \clip (0,0) circle (2cm);
      \draw[diag2=1,rdirected] (.3,-2) to (.3,2);
      \draw[diag2=1,rdirected] (-.3,2) to (-.3,-2);
      \fill[fill=colour_diag2,opacity=.2]
        (.3,-2) to (.3,2) to (-.3,2) to (-.3,-2) to (.3,-2);
      \draw[diag1=1,directed] (1,-2) to (1,2);
      \draw[diag1=1,directed] (1.7,2) to (1.7,-2);
      \fill[fill=colour_diag1,opacity=.2]
        (1,-2) to (1,2) to (1.7,2) to (1.7,-2) to (1,-2);
      \draw[diag1=1,rdirected] (-1,-2) to (-1,2);
      \draw[diag1=1,rdirected] (-1.7,2) to (-1.7,-2);
      \fill[fill=colour_diag1,opacity=.2]
        (-1,-2) to (-1,2) to (-1.7,2) to (-1.7,-2) to (-1,-2);
      \draw (-1,1) to (1,1);
      \draw (-1,-1) to (1,-1);
    }[scale=.6*\scl]
  \end{gather*}
  This concludes.
\end{proof}

\subsection{Addendum: another deformation of \texorpdfstring{$\glt$}{gl2}-foams}
\label{subsec:addendum_deformation_foam}

Given a presented graded-2-category, one can ask whether its defining relations can be ``deformed'' while keeping the same hom-basis.
Here we consider deformations of the form $a=\lambda b+\mu c$ where $a=b+c$ is a relation and $\lambda$ and $\mu$ are invertible scalars.
For instance, the graded-2-category $\gfoam_d$ is a deformation of the ``standard'' 2-category of $\glt$-foams, which is a posteriori recovered from $\gfoam_d$ by setting $X=Y=Z=1$.
In this addendum, we use rewriting to explore such ``deformations'' of $\glt$-foams, and show that under some conditions, $\gfoam_d$ and a variant $\gfoam_d'$ are the only two deformations of $\glt$-foams.

\medbreak

The structure of graded-2\nbd-cate\-gory on $\gfoam_d$ is already somehow universal. Indeed, the abelian group $\bZ^2$ is isomorphic to the abelian group presented by generators $D,\cup^L,\cup^R,\cap^L,\cap^R$ and relations
\begin{align*}
  \cap^R+\cup^R&=0\\
  \cap^L+\cup^L&=0\\
  \cup^R+\cap^L+D&=0\\
  \cup^L+\cap^R&=D\\
  \cup^L+\cup^R+\cap^L+\cap^R&=0
\end{align*}
obtained by taking the ``abelianization'' of the defining local relations of $\gfoam_d$.
Hence, if we assume that the grading is independent on the colours of the generators, then the $\bZ^2$-grading on $\gfoam_d$ is the most general one.
Moreover, a symmetric bilinear map $\bilfoam\colon\bZ^2\times\bZ^2\to\ringfoam^\times$ is determined by its values on the generators of $\bZ^2\times\bZ^2$, with relations:
\begin{align*}
  &\mu((1,0),(1,0))^2=\mu((0,1),(0,1))^2=1
  \quad\an\quad
  \mu((1,0),(0,1))\mu((0,1),(1,0))=1.
\end{align*}
This gives the parameters $X$, $Y$ and $Z$.

Let us now look how the defining relations could be deformed. For simplicity, we assume that the braid-like relations, pitchfork relations and dot slide (i.e.\ relations induced by $\foamX$) remain scalar-free.
Up to normalization, we can further assume that the scalars of the zigzag relations, the dot migration ($\dm$) and the evaluations of bubbles with external shadings ($\bbev$) keep the same scalars.

Pivotal $\foamE$-naturalities (\cref{lem:foam_pivotal_naturalities}) and the proof of \cref{lem:foam_T_congruent_iff_S_congruent} show that the neck-cutting cannot be modified.
Going over the remaining critical branchings, one finds exactly one extra possibility than $\gfoam_d$, where we have:
 \begin{IEEEeqnarray*}{cCc}
  \tikzpic{
    \funzip[0][1]\fzip\node at (.8+.5,-.3+1) {\scriptsize $i$};
  }[scale=.7]
  \mspace{5mu}\overset{}{\Rrightarrow}\mspace{5mu}
  XYZ
  \mspace{10mu}
  \tikzpic{
    \fdot\node at (.2,-.2) {\scriptsize $i$};
  }[][(0,0)]
  \mspace{10mu}
  \;+\;
  Z
  \mspace{10mu}
  \tikzpic{
    \fdot[0][0][2]\node at (.5,-.2) {\scriptsize $i+1$};
  }[][(0,0)]
  & \mspace{60mu} &
  \tikzpic{
    \flst[0][0]\frst[.5][0][2]\flst[1.5][0][2]\frst[2][0]
    \node[below=-2pt] at (0,0) {\scriptsize $i$};
    \node[below=-2pt] at (1.5,0) {\scriptsize $i+1$};
  }[yscale=4,scale=.35][(0,2*.35)]
  \Rrightarrow XYZ^{-1}\;
  \tikzpic{
    \clip (-1,0) rectangle (2,4);
    \fzip[0][3][2]
    \funzip[0][0][2]
    \begin{scope}[scale=2]
      \fcup[-.25][1]
      \fcap[-.25][0]
    \end{scope}
    \node[below=-2pt] at (-.5,0) {\scriptsize $i$};
    \node[below=-2pt] at (1,0) {\scriptsize $i+1$};
  }[scale=.35][(0,2*.35)]
\end{IEEEeqnarray*}
This defines an a priori distinct $(\bZ^2,\bilfoam)$-graded-2\nbd-cate\-gory $\gfoam_d'$. Renormalizing the rightward cap, we can also define $\gfoam_d'$ by only modifying the zigzag relations as follows:
\begin{equation*}
  \def\scl{.7}
  \tikzpic{
    \fcap[0][1]\flst[2][1]\flst\fzip[1][0]
    \node[left=-3pt] at (2,2) {\scriptsize $i$};
  }[scale=.5][(0,1*.5)]
  =
  \tikzpic{
    \flst[0][1]\flst\node[left=-3pt] at (0,2) {\scriptsize $i$};
  }[scale=.5][(0,1*.5)]
  \mspace{20mu}
  \tikzpic{
    \frst[0][1]\fcap[1][1]\fzip[0][0]\frst[2][0]
    \node[left=-3pt] at (0,2) {\scriptsize $i$};
  }[scale=.5][(0,1*.5)]
  = X
  \tikzpic{
    \frst[0][1]\frst\node[left=-3pt] at (0,2) {\scriptsize $i$};
  }[scale=.5][(0,1*.5)]
  \mspace{20mu}
  \tikzpic{
    \funzip[0][1]\frst[2][1]\frst\fcup[1][0]
    \node[left=-3pt] at (2,2) {\scriptsize $i$};
  }[scale=.5][(0,1*.5)]
  =XYZ^2
  \tikzpic{
    \frst[0][1]\frst\node[left=-3pt] at (0,2) {\scriptsize $i$};
  }[scale=.5][(0,1*.5)]
  \mspace{20mu}
  \tikzpic{
    \flst[0][1]\funzip[1][1]\fcup[0][0]\flst[2][0]
    \node[left=-3pt] at (0,2) {\scriptsize $i$};
  }[scale=.5][(0,1*.5)]
  = XZ^2\;
  \tikzpic{
    \flst[0][1]\flst\node[left=-3pt] at (0,2) {\scriptsize $i$};
  }[scale=.5][(0,1*.5)]
\end{equation*}
This other graded deformation (\cite[Definition~6.6.9]{Schelstraete_OddKhovanovHomology_2024}) verifies the same basis theorem (\cref{thm:foam_basis_theorem}), by exactly the same rewriting proof.
In particular, it verifies the same categorification theorem \cite[Theorem~2.29]{SV_OddKhovanovHomology_2023}.

\begin{remark}
  \label{rem:ladybugs}
  As noted in \cite[Remark~3.24]{SV_OddKhovanovHomology_2023}, the two variants $\gfoam_d$ and $\gfoam_d'$ correspond to the two variants of odd Khovanov homology, known as ``type X'' and ``type Y''. (Note that as homological invariants, these two variants are equivalent, while it is not clear to us whether $\gfoam_d$ and $\gfoam_d'$ are isomorphic.)
  The existence of these variants is explained by so-called ``ladybug squares'', certain squares in the hypercube that compose to zero.
  Interestingly, the same ladybug squares are at the heart of \emph{Khovanov homotopy type} \cite{LS_KhovanovStableHomotopy_2014}, a stable homotopy refinement of Khovanov homology. Most of the refinement is canonical, except on ladybug squares, for which a choice has to be made. This failure of canonicity is what allows Khovanov homotopy type to be a strictly stronger invariant than Khovanov homology
  \cite{
    LS_SteenrodSquareKhovanov_2014,
    Seed_ComputationsLipshitzSarkarSteenrod_2012,
  }.
\end{remark}

\printbibliography[heading=bibintoc]

\end{document}